\documentclass[11pt,a4paper,twoside,openright]{book}
\usepackage{epsfig}

\font\got=eufm10 at 11pt

\newtheorem{theorem}{Theorem}[section] 
\newtheorem{lemma}[theorem]{Lemma}
\newtheorem{proposition}[theorem]{Proposition}
\newtheorem{corollary}[theorem]{Corollary}
\newtheorem{definition}[theorem]{Definition}

\def\pagevide{\newpage{\pagestyle{empty}\cleardoublepage}}
\def\0{{\scriptscriptstyle 0}}
\def\Ad{\mathop{\rm Ad}}
\def\ad{\mathop{\rm ad}}
\def\a{{\mathcal A}}
\def\b{{\mathcal B}}
\def\build#1_#2^#3{\mathrel{\mathop{\kern 0pt#1}\limits_{#2}^{#3}}}
\def\bor{\mathop{\rm Bor}}
\def\c{{\mathcal C}}
\def\C{\hbox{\bf C}}
\def\card{\mathop{\rm Card}}
\def\diam{\mathop{\rm diam}}

\def\epsilon{\varepsilon}
\def\f{{\mathcal F}} 
\def\g{{\mathcal G}}
\def\h{{\mathcal H}}
\def\hol{\mathop{\rm hol}}
\def\id{\mathop{\rm Id}}
\def\i{{\mathcal I}}
\def\j{{\mathcal J}}
\def\kn{{\scriptstyle {k\over n}}}
\def\kpn{{\scriptstyle {{k+1}\over n}}}
\def\law{\mathop{\rm law}}

\def\l{{\mathcal L}}
\def\Log{\mathop{\rm Log}}
\def\lra{\longrightarrow}
\def\mod{\mathop{\rm mod}}
\def\n{{\mathcal N}}
\def\N{\hbox{\bf N}}
\def\R{\hbox{\bf R}}

\def\sR{\hbox{\scriptsize{\bf R}}}
\def\sigmaint{\sigma_{{\rm int}}}
\def\supp{\mathop{\rm Supp}}

\def\tr{\mathop{\rm tr}}
\def\v{{\mathcal V}}
\def\vol{\mathop{\rm vol}}
\def\wh{W\! H}
\def\Z{\hbox{\bf Z}}
\def\z{{\mathcal Z}}
\def\sZ{\hbox{\scriptsize{\bf Z}}}

\def\1{{\bf 1}}
\def\pf{\textsc{Proof}}
\def\qed{\hfill\hbox{\vrule\vbox to 2mm{\hrule width 2mm\vfill\hrule}\vrule}
  \\} 

\title{\Huge{Yang-Mills Measure on Compact Surfaces}}
\author{\Large{Thierry LEVY} \\  \\ Statistical Laboratory, Cambridge, UK \\
  \texttt{t.levy@statslab.cam.ac.uk}}
\date{January 2001}

\begin{document}

\frontmatter

\maketitle
\pagevide


\chapter*{Abstract} 
\thispagestyle{empty} We construct and study the Yang-Mills measure in two
dimensions. According to the informal description given by the physicists, it
is a probability measure on the space of connections modulo gauge
transformations on a principal bundle with compact structure group. We are
interested in the case where the base space of this bundle is a compact
orientable surface.

The construction of the measure in a discrete setting, where the base space of
the fiber bundle is replaced by a graph traced on a surface, is quite well
understood thanks to the work of E. Witten. In contrast, the continuum limit
of this construction, which should allow to put a genuine manifold as base
space, still remains problematic.

This work presents a complete and unified approach of the discrete theory and
of its continuum limit. We give a geometrically consistent definition of the
Yang-Mills measure, under the form of a random holonomy along a wide,
intrinsic and natural class of loops. This definition allows us to study
combinatorial properties of the measure, like its Markovian behaviour under
the surgery of surfaces, as well as properties specific to the continuous
setting, for example, some of its microscopic properties. In particular, we
clarify the links between the Yang-Mills measure and the white noise and show
that there is a major difference between the Abelian and semi-simple theories.
We prove that it is possible to construct a white noise using the measure as a
starting point and vice versa in the Abelian case but we show a result of
asymptotic independence in the semi-simple case which suggests that it is
impossible to extract a white noise from the measure.\\

\pagevide

\tableofcontents
\pagevide


\chapter*{Introduction}
\addcontentsline{toc}{chapter}{Introduction}
\markboth{INTRODUCTION}{INTRODUCTION}

This thesis is devoted to the construction and to the study of the Yang-Mills
measure in two dimensions. The quadridimensional and pseudo-Riemannian
equivalent of this measure is used by physicists in gauge theories such as
quantum electrodynamics and quantum chromodynamics, in order to describe the
fundamental interactions.  It appears in path integrals, which are known to be
often ill defined. Physicists describe the Yang-Mills measure as a probability
measure on the space of connections modulo gauge transformations on a
principal bundle with compact structure group. We are interested in the case
where the base space of this bundle is a compact orientable surface. The
informal expression of the measure is the following:
\begin{equation} \label{informal_expression}
d\mu(\omega)={1\over Z} \; e^{-{1\over 2}S(\omega)}
\;D\omega,
\end{equation}
where $S$ is the Yang-Mills action, that is, the $L^2$ norm of the
curvature. The constant $Z$ is a normalization constant and the measure
$D\omega$ should be a translation invariant measure on the space of
connections. 

The construction of the measure in a discrete setting, where the base space of
the fiber bundle is replaced by a graph traced on a surface, is quite well
understood thanks to the work of Witten \cite{Witten}. In contrast, the
continuum limit of this construction, which should allow to put a genuine
manifold as base space, still remains problematic. Several works
in this direction have led to substantial progress but not yet to an entirely
satisfactory solution.

This work presents a complete and unified approach of the discrete theory and
of its continuum limit. We give a geometrically consistent definition of the
Yang-Mills measure, under the form of a random holonomy along a wide,
intrinsic and natural class of loops. This definition allows us to study
combinatorial properties of the measure, like its Markovian behaviour under
the surgery of surfaces, as well as properties specific to the continuous
setting, for example, some of its microscopic properties. In particular, we
clarify the links between the Yang-Mills measure and the white noise and show
that there is a major difference between the Abelian and semi-simple theories.
We prove that it is possible to construct a white noise using the measure as a
starting point and vice versa in the Abelian case but we show a result of
asymptotic independence in the semi-simple case which suggests that it is
impossible to extract a white noise
from the measure.\\

\subsection*{Statement of the problem.} We are given a surface $M$, a Lie
group $G$ and a principal bundle $P$ over $M$. The surface $M$ is a
differentiable two-dimensional compact orientable manifold, with or without
boundary. It is endowed with a measure $\sigma$ which is equivalent to the
Lebesgue measure in any chart, with a positive smooth density. The group $G$
is a compact connected Lie group. In most examples, it will be either Abelian
or semi-simple. The fiber bundle $P$ is a principal fiber bundle over $M$ with
structure group $G$.

Recall that a connection on $P$ is a $G$-invariant choice of a horizontal
distribution in $P$ and that this choice can be represented by a {\got
  {g}}-valued 1-form $\omega$ on $P$, where {\got {g}} is the Lie algebra of
$G$. The curvature of the connection $\omega$ is the {\got{g}}-valued 2-form
$\Omega$ on $P$ defined by
\begin{equation} \label{courb}
\Omega(X,Y)=d\omega(X,Y)+[\omega(X),\omega(Y)].
\end{equation}
The curvature can be considered an $\ad P$-valued 2-form on $M$, where $\ad P$
is the fibre bundle associated with $P$ by the adjoint representation of $G$
on {\got{g}}. If we choose an orientation on $M$, $\Omega$ can be identified
with a section of $\ad P$. An $\ad$-invariant scalar product on {\got{g}}
allows to define a metric on $\ad P$ and hence the norm $\parallel \Omega
\parallel$ of the curvature. This norm does not depend on the choice of the
orientation of $M$ and the Yang-Mills action is defined on the space $\a$ of
connections on $P$ by
\begin{eqnarray*}
S: \a & \lra & \R_+ \\
\omega &\longmapsto & S(\omega)=\int_M \parallel \Omega \parallel^2 \;
d\sigma. 
\end{eqnarray*}

Our aim is to give a sense to the informal expression
(\ref{informal_expression}). The first problem is of course that there is no
translation invariant measure on the infinite dimensional affine space $\a$.
Another one is the invariance of the action $S$ under the action of a huge
group, that of gauge transformations of $P$. This group, denoted by $\j$, is
the group of diffeomorphisms of $P$ over the identity of $M$ that commute with
the action of $G$. It acts by pull-back on $\a$ and preserves $S$, since it
acts on the curvature by pointwise conjugation, which does not change the norm
in $\ad P$.  Because of this invariance, the constant $Z$ should be
proportional to the volume of $\j$, hence be infinite. To avoid this problem,
we try to construct the measure on the quotient space $\a/\j$ instead of $\a$.
This means that we will be able to integrate only gauge-invariant functions
against the Yang-Mills measure, in agreement with the physical principle
saying that observable quantities must be gauge-invariant. On the other hand,
this quotient space has a much more complicated structure than an affine
space.  This is why one usually tries to avoid to work directly on it and
prefer to work on a function space, as we explain below.  The work of D. Fine
\cite{Fine_S2,Fine_RS} is an exception from this point of view, since the
author analyzes the geometrical structure of $\a/\j$ in order to give sense to
(\ref{informal_expression}).

\subsection*{The Yang-Mills measure as random holonomy.} A starting point may
be to ask what functions we want to be able to integrate against the
Yang-Mills measure. Physicists' answer this question is that we must be able
to integrate Wilson loops.

A connection $\omega$ on $P$ defines a parallel transport along regular paths
on $M$. The parallel transport along a given path $c:[0,1]\lra M$ is a
$G$-equivariant diffeomorphism of the fiber $P_{c(0)}$ into the fiber
$P_{c(1)}$, denoted by $\hol(\omega,c)$. If $c$ is a loop and if we fix a
point $p$ in the fiber $P_{c(0)}$, this diffeomorphism can be represented by
the element $g$ of $G$ such that $\hol(\omega,c)(p)=pg$. If we choose another
point in $P_{c(0)}$, we find another element of $G$ conjugate to $g$.
So, for any representation $\rho$ of $G$ and any loop $l$, one defines the
Wilson loop $W_{l,\rho}$ by
$$W_{l,\rho}(\omega) = \tr \rho(\hol(\omega,l)).$$

The functions that we want to integrate are central functions of the holonomy
along loops. We just noticed that the holonomy along a loop determines a
conjugacy class in $G$. We must also take into account the action of the gauge
group, which conjugates by the same element of $G$ the holonomies along all
loops based at the same point. Let $LM$ denote the set of regular paths on $M$
and $\f(M,G)$, $\f(LM,G)$ the sets of $G$-valued functions on $M$, $LM$. An
element $j$ of the group $\f(M,G)$ acts on an element $f$ of $\f(LM,G)$ by:
$$j.f(l) = j(l(0))^{-1} f(l) j(l(0)).$$
It is possible to define a map from
$\a/\j$ into the quotient $\f(LM,G)/\f(M,G)$, mapping a connection to the
class of the holonomy that it determines along the elements of $LM$. An
argument of Sengupta \cite{Sengupta_S2} proves that this map is injective.
Thus, we change our point of view: we seek now a measure on the space
$\f(LM,G)/\f(M,G)$, viewing this space as a space of generalized connections
modulo gauge transformations. In fact, we shall construct a measure on
$\f(LM,G)$ and take the quotient of this measure by $\f(M,G)$. In other words,
we really want to construct a random holonomy instead of a random connection.
This will be easier because we can use the classical tools of probability to
construct a measure on a function
space. \\

At this point, it is necessary to characterize more precisely the Yang-Mills
measure. Either one tries to extract more information from the informal
expression of the measure or one looks for other description of this measure.
The last option is the one that we choosed in this work, using the
combinatorial description given by Migdal and Witten. The first one is based
on the Gaussian character of the measure and was used by Driver and Sengupta
\cite{Driver, Sengupta_S2,Sengupta_AMS}.

\subsection*{Gaussian interpretation: curvature of the random
  connection.} Assume that $G$ is abelian, for example $G=U(1)$. The relation
(\ref{courb}) between a connection and its curvature becomes linear. A formal
change of variables in (\ref{informal_expression}) gives
\begin{equation} \label{ccbb}
d\mu(\Omega)={1\over {Z'}} e^{-{1\over 2}\parallel \Omega \parallel^2} \;
D\Omega.
\end{equation}
Recall that the curvature may be seen as an $\ad P$-valued two-form, or as a
section of $\ad P$ if an orientation of $M$ is fixed. Since $G$ is Abelian,
its adjoint representation on its Lie algebra is trivial. Hence, the fiber
bundle $\ad P$, which is the vector bundle associated to $P$ for this adjoint
representation, is trivial, and may therefore be identified with $M \times$
{\got{g}}. We recognize in (\ref{ccbb}) the expression of a Gaussian measure
on the Hilbert space of square integrable {\got{g}}-valued functions on $M$.
This leads us to the main idea of the interpretation of
(\ref{informal_expression}): under the Yang-Mills measure, the random
curvature of a connection has a Gaussian distribution, it is a
{\got{g}}-valued white noise on $M$ with intensity $\sigma$.

This argument is of course specific to the abelian case, since in general,
$\Omega$ is a quadratic function of $\omega$. Nevertheless, the fact that $M$
is two-dimensional allows to overcome this problem: it is always possible to
get back to a situation similar to the abelian case by a gauge fixing
procedure. This requires a word of explanation.

If we choose a local trivialization of $P$ on an open subset $U\subset M$,
i.e. a local section $s: U \lra P$ of $P$, we can pull-back by $s$ all objects
living on $P$, in particular the connection and curvature forms. One denotes
usually $A=s^*\omega$ and $F=s^*\Omega$.
These forms on $M$ satisfy a structure equation $F=dA+[A,A]$ identical
to (\ref{courb}). Let $j$ be an element of the gauge group $\j$. This element
can act in two different ways in this situation, either by transforming the
section $s$ into $j\circ s$ or by transforming the forms $\omega$ and $\Omega$
into $j^*\omega$ and $j^*\Omega$. These two ways are indistinguishable from
the base space, since $(j\circ s)^*=s^*j^*$. So, we denote without ambiguity
$s^*j^*\omega=(j\circ s)^*\omega$ by $A^j$. 

What is specific to the two-dimensional case is that given a connection, it is
possible to choose $s$ in such a way that $[A,A]=0$. Choose $U$ small enough
to admit local coordinates $x,y$. Set $m=(0,0)$ in these coordinates and
choose $p$ in $P_m$. Then define $s$ along the $y$-axis by lifting it
horizontally, starting at $p$. Now, starting from each point $(0,y_0)$ in $U$,
define $s$ on the line through $(0,y_0)$ parallel to the $x$-axis by lifting
it horizontally, starting at $s(0,y_0)$. The section $s$ is smooth and
horizontal along all lines parallel to the $x$-axis. Thus, $A=s^* \omega$
takes the form $A=A_y dy$ and $[A,A]=0$.  When one looks at $P$ through $s$,
what one sees is similar to the abelian case, up to the fact that the section
through which the relation between connection and curvature should be linear
depends on the connection. 

\subsection*{From the curvature to the holonomy.} The next step is to define a
random holonomy using the random curvature. The method is based on
deterministic links between holonomy and curvature. Assume that $G=U(1)$ and
take $\R^2$ as base space, although it is not a compact surface. Given a
connection $\omega$ on the fiber bundle $\R^2\times G$ and a simple loop $l$
which bounds a domain $D$, the Stokes formula gives
$$\hol(\omega,l)=\exp i \oint_l A = \exp i \int_D dA = \exp i \int_D F= \exp i
(F,\1_D)_{L^2}.$$
This formulation is easily extended to the random case.
Indeed, pick a white noise $W$ on $\R^2$, i.e. an isometry from $L^2(\R^2)$
into a vector space of Gaussian random variables. One can replace $F$ by $W$
in the last expression and define a random holonomy along $l$ by
$$H_l=\exp i W(\1_D).$$
The construction that we present in the abelian case
in chapter \ref{abelian_theory} is an extension of this procedure to surfaces
whose topology is non trivial and where the interior of a loop is not well
defined.

It is possible here to understand better the difficulties of Driver and
Sengupta. They tried to use this method in the case of a non-Abelian structure
group. But in this case, the holonomy is not $\exp \int_l A$, but $P \exp
\oint_l A$, which is a compact notation for the solution of the differential
equation 
$$\cases{\dot h_t = A(\dot l(t)) h_t \cr h_0=1.}$$
The Stokes formula does not work in this frame. In some sense, one has to
choose in which order one multiplies the small elements of $G$ obtained by
integrating $F$ over small squares inside $D$. It is not surprising that
Driver and Sengupta had to use the coordinate system that allows to define a
section through which $[A,A]=0$ in order to determine this order. The problem
is that the class of loops along which they are able to define the random
holonomy  depends strongly on this choice of coordinates.

It should be noted that B.  Driver \cite{Driver_lassos} and L. Gross
\cite{Gross} introduced a new local object in order to replace the white noise
in this context and that this could lead to a way around the problem.

Although we do not treat this point in our work, we cannot conlude this
conclusion without mentioning the semi-classical limit of the Yang-Mills
measure. The remarkable fact is that, when the total surface of $M$ tends to
$0$, the measure concentrates on the set of flat connections over $M$ and
tends to the volume measure associated with the natural symplectic structure
on this space. There are a lot of references on this subject which is closely
related to the geometry of some moduli spaces
\cite{Forman,Becker_Sengupta,King_Sengupta,Sengupta_SCL,Liu,Atiyah_Bott}.

\subsection*{Combinatorial approach.} Our starting point is the combinatorial
approach initiated by A.A. Migdal in 1975 \cite{Migdal} and improved by
E.Witten in 1991 \cite{Witten}. The idea is to replace the base space of the
fiber bundle by a graph on a surface. This leads to a finite dimensional
problem, where we define the random holonomy only along the paths of a graph.
We also define conditional versions of the random holonomy. These conditional
versions will play a technical role in the continuous construction on surfaces
with boundary and lead also to the definition of very important objects, the
conditional partition functions. The main property of the discrete theory is
the invariance by subdivision. It explains that, up to some restrictions, the
law of the random holonomy is independent of the graph in which one works.

The next step towards the continuum limit is to take the projective limit of
the discrete measures associated with the graphs whose edges are piecewise
geodesic for some Riemannian metric on $M$. This allows us to define a random
holonomy along all piecewise geodesic paths on $M$. Then, we prove that this
random holonomy can be extended by continuity to the set $PM$ of piecewise
embedded paths on $M$, using a very natural approximation procedure. The law
of this new random holonomy is a measure on $\f(PM,G)$ which is proved to be
independent of all choices made during the construction. This measure is
pushed forward on $\f(LM,G)/\f(M,G)$ and then becomes what we call the
Yang-Mills measure. This measure is characterized by its consistence with the
discrete theory and a continuity property. It is multiplicative, as a random
holonomy is expected to be, and invariant by area-preserving diffeomorphisms.

All along the discrete construction, we study the special case $G=U(1)$.  This
analysis leads us to a second construction of a random holonomy, specific to
the Abelian case, based on the Gaussian character of the measure in two
dimensions. The characterization of the Yang-Mills measure given earlier
allows us to prove that this random holonomy has the same law as that defined
by the general procedure. Then we show that the holonomy along very small
loops can be used to construct a white noise on $M$, by means of a Wiener-like
integral.

It is then natural to try to adapt the extraction of the white noise to the
general case. We prove a result in the semi-simple case that strongly suggests
that this is impossible. Indeed, the $\sigma$-algebra generated by the
holonomies along very small loops seems to satisfy a zero-one law.

In the last part of this work, we study combinatorial properties of the
measure. We prove the Markov property of the random holonomy, extending to the
continuous setting a result that was proved in the discrete setting by C.
Becker and A. Sengupta \cite{Becker_Sengupta}.  Then, we study how it is
possible to glue together the Yang-Mills measures on two surfaces $M_1$ and
$M_2$ in order to get the measure on $M$, the surface obtained by gluing $M_1$
and $M_2$ together. We show that the measures on $M_1$ and $M_2$ do not
determine the measure on $M$. There is a lack of information that can be
parametrized by the centralizer of the holonomy along the common boundary of
$M_1$ and $M_2$.  Finally, we summarize the algebraic properties of the
conditional partition functions, whose importance had already been recognized
by Witten \cite{Witten}. We prove that a few of them generate all others by
algebraic transformations and identify these elementary functions. We also
show that the partition functions may be considered the transition functions
of the random holonomy as a Markov field and discuss to what extent they
determine the Yang-Mills measure.


\pagevide

\mainmatter

\chapter{Discrete Yang-Mills measure}

In this chapter we construct and we study the discrete Yang-Mills measure. It
is both the basis of the construction of the continuous Yang-Mills measure and
the frame in which computations are possible. The main results are the
invariance by subdivision of the discrete measure and the estimation of the
law of the random holonomy along small loops. 

\section{Notations}

Throughout this work, $M$ will denote a surface, i.e. a real differentiable
two-dimensional manifold, compact, connected, orientable, with or without
boundary. It is endowed with a Lebesguian measure $\sigma$, i.e. a measure
which has positive smooth density with respect to the Lebesgue measure in any
chart.

The boundary of $M$, if it is non empty, is the disjoint union of circles
$N_1,\ldots, N_p$. Let us make explicit what we call smooth objects on $M$ and
introduce the very useful notion of closure.

\begin{definition} \label{closure}
  A closure of $M$ is a triple $(i,M,M_1)$, where $M_1$ is a closed surface,
  i.e. a surface without boundary and $i :M \lra M_1$ is an embedding. If the
  complementary of $i(M)$ in $M_1$ is diffeomorphic to a disjoint union of
  disks, the closure is said to be minimal.
\end{definition}

Given two closures $(i_1,M,M_1)$ and $(i_2,M,M_2)$ of $M$, $i_1(M)$ and
$i_2(M)$ have diffeomorphic neighbourhoods in $M_1$ and $M_2$. So it makes
sense to say that an application (resp. a bundle, a section,...) is smooth on
$M$ if it is the restriction of a smooth application (resp. bundle,
section,...) defined on an open neighbourhood of $M$ in one of its closures.

The second basic object is $G$, a compact connected Lie group, that will be
chosen to be Abelian or semi-simple in most examples.

Let $P$ be a principal $G$-bundle over $M$. If $M$ has a boundary, $M$
retracts on a bunch of circles and $P$ is trivial. But if $M$ is closed, the
possible topological types for $P$ are classified by $\pi_1(G)$. A pleasant
way to see this is to cut $M$ along the boundary of a small disk. We get two
disjoint pieces. The restrictions of a bundle $P$ over $M$ to both pieces are
trivial and the topology of $P$ is completely determined by the transition
function along the boundary of the disk. This transition function is a map
$S^1 \lra G$ and it is a fact that two homotopic maps give rise to two
homeomorphic bundles. 

If $G=U(1)$, the element of $\pi_1(U(1))\simeq \Z$ determined by $P$
corresponds to the Chern class of the complex line bundle associated with $P$.
Note that when $G$ is semi-simple, $\pi_1(G)$ is finite.

A connection $\omega$ on $P$ is a choice at each point $p$ of $P$ of a
subspace in $T_pP$ supplementary to the vertical subspace of vectors tangent
to the action of $G$. Moreover, this distribution, called the horizontal
distribution, has to be invariant by the action of $G$.

Let $c:[0,1]\lra M$ be a regular path on $M$. A connection $\omega$ allows to
lift $c$ to a horizontal path in $P$ starting at any prescribed point in
$P_{c(0)}$. The function that maps a point $p$ of $P_{c(0)}$ to the end point
of the horizontal lift of $c$ starting at $p$ is called the parallel
transport or holonomy of $\omega$ along $c$. It is a $G$-equivariant map
$\hol(\omega,c):P_{c(0)} \lra P_{c(1)}$. If $c_1$ and $c_2$ are two paths such
that $c_1(1)=c_2(0)$, then the path $c_1 c_2$ exists and we have
$$\hol(\omega,c_1 c_2) = \hol(\omega,c_2)\circ \hol(\omega,c_1).$$

A gauge transformation is a diffeomorphism $j$ of $P$ over the identity of $M$
that commutes with the right action of $G$. Let $\omega$ be a connection on
$P$. Let $c :[0,1] \lra M$ be a piecewise $C^1$ path. A gauge transformation
$j$ allows to define a new connection $j^*\omega$ whose holonomy is related to
that of $\omega$ through the relation:
$$\hol(j^*\omega,c)=j(c(1))^{-1} \circ \hol(\omega,c) \circ j(c(0)).$$
Remark
that these holonomies are conjugate if $c$ is a loop. More detailed
presentations of the theory of fiber bundles and connections can be found for
example in \cite{Kobayashi_Nomizu_I,Bleecker}.

\section{Graphs on $M$}

In order to reduce to a discrete setting, we will replace $M$ by a graph drawn
on $M$ and adapt the notions of fibre bundle, connection and gauge
transformation. 

\subsection{Pregraphs}

We say that an application $c :[0,1]\lra M$ is smooth (resp. an embedding) if
it is the restriction of a smooth application (resp. embedding) defined on an
open interval containing $[0,1]$.

\begin{definition} A parametrized path on $M$ is a mapping $c
  :[0,1]\lra M$ which is the concatenation of a finite number of smooth
  embeddings.
\end{definition}

Two parametrized paths are said to be equivalent if they differ by an
increasing reparametrization.

\begin{lemma} The equivalence of parametrized paths preserves their
  orientation, image, end points, injectivity, injectivity on $(0,1)$.
\end{lemma}

Equivalence classes of parametrized paths are called simply paths. The set of
paths on $M$ is denoted by $PM$. 

A path whose end points are equal is called a loop and a loop which is
injective on $(0,1)$ is said to be simple. Given a path $a$, we denote by
$a^{-1}$ the path obtained by reversing the orientation of $a$. An edge is an
injective path $a$ such that $a([0,1]) \cap \partial M$ is empty or a finite
union of segments.

\begin{definition} A pregraph on $M$ is a set $\Gamma=\{a_1,\ldots,a_r\}$ of
  edges that meet each other only at their end points, i.e. such that for each
  distinct $i$ and $j$ between $1$ and $r$, one has
$$a_i([0,1]) \cap a_j([0,1]) =a_i(\{0,1\}) \cap a_j(\{0,1\}).$$
\end{definition}

We call support of a pregraph $\Gamma$ the union of the images of its edges. A
pregraph $\Gamma$ is said to be connected if its support $\supp(\Gamma)$ is
connected.

We call faces of a pregraph $\Gamma$ the connected components of $M\backslash
\supp(\Gamma)$. We denote by $\f(\Gamma)$ the set of these faces.

\begin{proposition} \label{disks_imply_H1} Let $\Gamma$ be a connected
  pregraph on $M$. Suppose that every face of $\Gamma$ is diffeomorphic to a
  disk. Then the map $H_1(\supp(\Gamma);\Z) \lra H_1(M;\Z)$ induced by the
  inclusion is surjective.
\end{proposition}

\pf. Let $c:[0,1] \lra M$ be a loop. There exists on each face of $\Gamma$ a
point that is not in the image of $c$. Let us fix such a point in each face
and remove it from $M$. The remaining open set $U$ retracts on the support of
$\Gamma$ because each face with a point removed retracts on its boundary. This
retraction induces a homotopy from $c$ to a loop whose image is included in
$\supp(\Gamma)$. So each loop of $M$ is homotopic, thus homologous to a loop
of $\supp(\Gamma)$. This proves the result. \qed

\subsection{Graphs}

Given a pregraph $\Gamma$, we call path in $\Gamma$ a concatenation of edges
of $\Gamma$, with natural or reverse orientation. We denote by $\Gamma^*$ the
set of these paths.

\begin{definition} A graph on $M$ is a connected pregraph $\Gamma$ whose
  faces are diffeomorphic to disks and such that for each component $N_i$ of
  $\partial M$, there exists an element of $\Gamma^*$ whose image is equal to
  $N_i$.
\end{definition}

The reason for which we choose these properties is that we want a graph to
take the whole topology of $M$ into account, including its boundary.

\begin{definition} \label{finess}
  Let $\Gamma_1$ and $\Gamma_2$ be two pregraphs. We say that $\Gamma_2$ is
  finer than $\Gamma_1$ and write $\Gamma_1 < \Gamma_2$ if each edge of
  $\Gamma_1$ is a path of $\Gamma_2^*$.
\end{definition}

\begin{proposition} \label{complete_a_pregraph} 
  Let $\Gamma_0$ be a pregraph on $M$. There exists a graph $\Gamma$ which is
  finer than $\Gamma_0$. Moreover, if $\supp(\Gamma_0)$ is contained in an
  open set diffeomorphic to a disk, it is possible to construct $\Gamma$ in
  such a way that it has the same number of faces as $\Gamma_0$.
\end{proposition}

\begin{lemma} \label{disks_simply_connected} A pregraph whose faces are
  diffeomorphic to disks is necessarily connected.
\end{lemma}

\pf. Let $\Gamma$ be a pregraph with faces diffeomorphic to disks. Suppose
that $\Gamma=\Gamma' \cup \Gamma''$, where $\Gamma'$ and $\Gamma''$ have
disjoint supports. Each face of $\Gamma$ is a disk, so it has a connected
boundary, which is included either in the support of $\Gamma'$ or in that of
$\Gamma''$. The closures of the unions of faces whose boundary lies in
$\supp(\Gamma')$ (resp. $\supp(\Gamma'')$) form a partition of $M$ into two
closed sets, which is in contradiction with connectedness of $M$. \qed

\begin{figure}[h]
\begin{center}
\scalebox{0.7}{\input{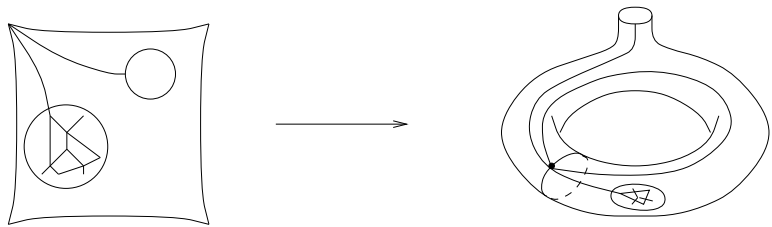}}
\end{center}
\caption{Fundamental domain in the universal covering for a torus with one
  hole.
\label{univ}}
\end{figure}

Before to prove proposition \ref{complete_a_pregraph}, let us recall some
classical facts about the topology of $M$. If $M$ has no boundary and is not a
sphere, its universal covering is diffeomorphic to a plane. In this plane, it
is always possible to choose a polygonal fundamental domain for the covering
map, namely a $4g$-gonal domain if $g$ is the genus of $M$. This means that it
is possible to see topologically $M$ as the result of the identification of
some edges of a polygon. If $M$ is a shpere, it can be seen as a disk whose
upper and lower half of the boundary have been identified. If $M$ has a
boundary, there are holes in the universal covering, one for each boundary
component in a fundamental domain. It is possible to choose a fundamental
domain such that the holes are in its interior. Thus, it is possible to
represent a surface with boundary by a picture like picture \ref{univ}.\\

\pf {\textsl{of proposition} \ref{complete_a_pregraph}}. By definition of an
edge, the set $\supp(\Gamma_0)\cap \partial M$ is a finite union of segments.
Cutting some edges of $\Gamma_0$ in several pieces if necessary, we can assume
that these segments are exactly images of edges. Then, it is possible to add
to $\Gamma_0$ edges in such a way that $\supp(\Gamma_0)\cap \partial
M=\partial M$. If $\supp(\Gamma_0)$ did not meet a component of $\partial M$
initially, it is necessary to add at least two edges on this component. So we
can construct a pregraph $\Gamma_1$ which is finer than $\Gamma_0$, whose
support contains $\partial M$.

Each face of $\Gamma_1$ is homeomorphic to the interior of a compact surface
with boundary. On any such surface, there exists a graph, for example a
triangulation. We add to $\Gamma_1$ the edges that are necessary to transform
it into a graph on each face which is not diffeomorphic to a disk. All faces
of the resulting pregraph $\Gamma_2$ are diffeomorphic to disks. By the
preceding lemma, it is connected. Thus, it is a graph.  The first part of the
proposition is proved.

If $\supp(\Gamma_0)$ is contained in a disk, it is possible to move this disk
by a diffeomorphism of $M$ into any prescribed disk. It is possible to get the
situation described by the picture, in the universal covering, where the disk
is in the interior of a fundamental domain. Then it is easy to complete
$\Gamma_0$ into a graph with the same number of faces and get back to the
initial situation by the inverse diffeomorphism. \qed

\section{Discrete holonomy and gauge transformations}
\label{discrete_holonomy}

Choose a graph $\Gamma=\{a_1,\ldots,a_r\}$ on $M$ and denote by
$i:\supp(\Gamma) \lra M$ the canonical injection. Consider the fiber bundle
$i^*P$ on $\supp(\Gamma)$. It is a trivial bundle, whatever the topology of
$P$ was. We identify $i^*P$ with $\supp(\Gamma)\times G$. Let $\v(\Gamma)$
denote the set of vertices of $\Gamma$.

\begin{lemma} Let $\omega_1$ and $\omega_2$ be two connections on
  $i^*P$. Suppose that $\omega_1$ and $\omega_2$ have the same holonomy along
  each edge of $\Gamma$. Then there exists a gauge transformation $j$ of
  $i^*P$ that leaves the fibers over the points of $\v(\Gamma)$ invariant and
  such that $j^*\omega_1=\omega_2$.
\end{lemma}

\pf. Choose a parametrization of $a_1\in\Gamma$. Set
$$j(a_1(t))=\hol(\omega_1,{a_1}_{|[0,t]}) \circ
\hol(\omega_2,({a_1}_{|[0,t]})^{-1}).$$
Then 
$$\hol(\omega_2,{a_1}_{|[0,t]})= j(a_1(t))^{-1} \circ
\hol(\omega_1,{a_1}_{|[0,t]})= \hol(j^*\omega_1,{a_1}_{|[0,t]}).$$
The assumption about $\omega_1$ and $\omega_2$ makes sure that it is
possible to extend the construction of $j$ to the whole graph and that
$j(m)=\id_{P_m}$ for all vertex $m$ of $\Gamma$. 

This implies that $\omega_2$ and $j^*\omega_1$ have the same holonomy, hence
the same horizontal paths, so they are equal.\qed

In the discrete setting, we expect to be able to compute the holonomy only
along edges of the graph. So we identify a connection with the holonomy that
it determines along the edges of $\Gamma$ and, according to the preceding
lemma, consider gauge transformations that act only on the fibers over the
vertices of $\Gamma$. Finally, using the identification $i^*P =
\supp(\Gamma)\times G$, we can identify holonomies and gauge transformations
with elements of $G$, defining for example $h_c(\omega)$ and $j(m)$ by
$\hol(\omega,c)(c(0),1)=(c(1),h_c(\omega))$ and $j(m,1)=(m,j(m))$.

This leads us to the following definitions.

\begin{definition} A discrete connection on $\Gamma$ is a map from $\Gamma$
  into $G$. A discrete gauge transformation is a map from the set $\v(\Gamma)$
  of vertices of $\Gamma$ into $G$. 
\end{definition}

A discrete gauge transformation $j:\v(\Gamma) \lra G$ acts on a discrete
connection $g=(g_1,\ldots,g_r)$ by :
$$j\cdot g=(j(a_1(1))^{-1} g_1 j(a_1(0)),\ldots,j(a_r(1))^{-1} g_r
j(a_r(0))).$$

A discrete connection $g=(g_1,\ldots,g_r)$ determines a multiplicative
application from $\Gamma^*$ into $G$. Given a path $c=a_{i_1}^{\epsilon_1}
\ldots a_{i_n}^{\epsilon_n}$, $\epsilon_i=\pm 1$, we can compute
$g_{i_n}^{\epsilon_n} \ldots g_{i_1}^{\epsilon_1}$ (with reversed order!). In
other words, any path $c$ of $\Gamma^*$ gives rise to a map $h_c$ from
$G^\Gamma$ to $G$ defined by
$$h_c(g_1,\ldots,g_r)=g_{i_n}^{\epsilon_n} \ldots g_{i_1}^{\epsilon_1}.$$
This map is well defined because there is only one way to decompose a path in
product of edges. To see why, it is enough to consider the times
at which a path crosses a vertex of $\Gamma$.

\begin{proposition} Let $c_1$ and $c_2$ be two paths such that $c_1 c_2$ is
  also a path. Then 
$$h_{c_1 c_2}=h_{c_2}h_{c_1}.$$
\end{proposition}

This basic property of the discrete holonomy will be refered to as the
multiplicativity of the holonomy.

A gauge transformation $j$ transforms $h_c$ in $h_c \circ j$, with 
$$h_c \circ j=j(c(1))^{-1} h_c j(c(0)).$$

At this stage, it may be noticed that central functions of the holonomy along
loops are invariant under gauge transformations. This is why they will play
such a major role in the sequel.

\section{Discrete Yang-Mills measure}
\label{discrete_YM}

We keep a graph $\Gamma=\{a_1,\ldots,a_r\}$ on $M$. In the discrete setting, a
probability measure on the quotient space $\a /\j$ of connections modulo gauge
transformations is represented by a probability measure on $G^\Gamma$
invariant under the action of $G^{\v(\Gamma)}$.

The basic example of such an invariant measure is the product of Haar
measures. We shall construct the discrete Yang-Mills on $G^\Gamma$ as :
$$dP={{dP}\over{dg}} \; dg,$$
where $dg=dg_1\otimes \ldots \otimes dg_r$. The
density ${{dP}\over {dg}}$ will be a product of central functions of
holonomies along loops, a feature that makes $P$ invariant. Recall that a
function $p$ is said to be central on $G$ if $p(xy)=p(yx)$, or equivalently
$p(y^{-1}xy)=p(x)$ for all $x,y$ in $G$. 

Each face $F$ of $\Gamma$ has a boundary which is the image of a path defined
up to the choice of an origin and an orientation. So the function $h_{\partial
  F} : G^\Gamma \lra G$ is defined up to conjugation and inversion. For any
central function $p$ invariant by inversion, the function $p\circ h_{\partial
  F}$ is well defined.

Let us denote by $(p_t)_{t>0}$ the fundamental solution of the heat equation
on $G$ endowed with its biinvariant Riemannian metric, normalized to have
total volume equal to $1$. It satisfies
$$(\partial_t - {1\over 2}\Delta) p_t = 0 \;\;\; {\rm on} \;\;\; \R^*_+ \times
G,$$
and for any function $f$ continuous on $G$,
$$\int_G f(g)p_t(g)\; dg \build{\lra}_{t \to 0}^{} f(1).$$
For any positive $t$, $p_t$ is a positive central function, invariant by
inversion, such that $\int_G p_t(g)\; dg =1$. For the moment, the choice of
the heat kernel may seem to be quite arbitrary. We shall discuss this point
at the end of section \ref{inv-by-sub}.

For each face $F$ of $\Gamma$, the function $p_{\sigma(F)}(h_{\partial F})
:G^\Gamma \lra \R^*_+$ is well defined, with $\sigma$ denoting the surface
measure on $M$. Set
\begin{equation}
D=\prod_{F\in \f(\Gamma)} p_{\sigma(F)}(h_{\partial F}) :G^\Gamma \lra
\R^*_+,
\end{equation}
\begin{equation}
Z=\int_{G^\Gamma} D\; dg.
\end{equation}
Now define $P$ on $(G^\Gamma,\bor(G^\Gamma))$ by
\begin{equation}
dP={1\over Z} \; D \; dg.
\end{equation}
Given $n$ paths $c_1,\ldots,c_n$ in $\Gamma^*$, we define the law of the
discrete holonomy along $c_1,\ldots,c_n$ as the joint law of the $n$-uple
$(h_{c_1},\ldots,h_{c_n})$ under $P$.

\section{Conditional Yang-Mills measure}
\label{conditionning}

When $M$ has a boundary, it is natural to want to impose the holonomy along
the components of $\partial M$. It may also be useful to be able to impose
holonomy along some other loops even if $M$ has no boundary.

\subsection{Conditional Haar measure}

\begin{proposition} \label{conditioned_Haar_measure} Let $n$ be a positive
  integer. Let $x$ be an element of $G$. There exists on $G^n$ a measure
  $\nu_x^n$ such that $g_n\ldots g_1=x$ $\nu_x^n$-a.s. and such that for any
  function $f$ continuous on $G^n$ and any $i$ between $1$ and $n$,
  \begin{eqnarray*}
  && \hskip -.5cm \nu_x^n(f)=\int_{G^{n-1}}f(g_1,\ldots,g_{i-1},(g_n\ldots
  g_{i+1})^{-1}x(g_{i-1}\ldots g_1)^{-1},g_{i+1},\ldots,g_n)\\
  && \hskip 9cm  dg_1\ldots \widehat{dg_i} \ldots dg_n.
  \end{eqnarray*}
  Moreover, one has
  $$\nu_x^n(f)=\lim_{t\to 0} \int_{G^n} f(g_1,\ldots,g_n) \; p_t(g_n \ldots
  g_1 x^{-1})\; dg_1\ldots dg_n.$$
  Finally, $\int_G \nu_x^n \; dx = dg$, as
  measures on $G^n$.
\end{proposition}

\pf. Pick $i$ between $1$ and $n$, and $t>0$. By centrality of $p_t$ and then
right invariance of $dg_i$, one has
\begin{eqnarray*}
&& \hskip -.5cm \int_{G^n} f(g_1,\ldots,g_n) p_t(g_n\ldots g_1 x^{-1}) \; dg =
\\ 
&=&  \int_{G^n} f(g_1,\ldots,g_n)\; p_t(g_i(g_{i-1}\ldots g_1)x^{-1}(g_n \ldots g_{i+1})) \;dg \\
&=&  \int_{G^n} f(g_1,\ldots,g_{i-1},g_i(g_n\ldots
g_{i+1})^{-1}x(g_{i-1}\ldots 
g_1)^{-1},g_{i+1},\ldots,g_n)\; p_t(g_i) \\ 
&& \hskip 8.2cm dg_i dg_1 \ldots \widehat{dg_i}\ldots dg_n \\
& \build{\lra}_{t\to 0}^{}& \int_{G^{n-1}}f(g_1,\ldots,g_{i-1},(g_n\ldots
g_{i+1})^{-1}x(g_{i-1}\ldots g_1)^{-1},g_{i+1},\ldots,g_n) \\ 
&& \hskip 8.6cm dg_1\ldots \widehat{dg_i} \ldots dg_n. 
\end{eqnarray*}
Thus the limit exists and the last expression
does not depend on $i$. It defines a probability measure on $G^n$.

If $f$ vanishes on the hypersurface $\{g_1\ldots g_n=x\}$, then
$\nu_x^n(f)=0$. 
So we do have $g_n\ldots g_1=x$ $\nu_x^n$-a.s.

Finally, since $\int_G p_t(g)\; dg =1$,
$$\int_G\int_{G^n} f(g_1,\ldots,g_n) \; p_t(g_n\ldots g_1 x^{-1}) \; dx dg
= \int_{G^n} f(g_1,\ldots,g_n) \; dg,$$
which implies the last statement when $t$ tends to zero.  \qed

\subsection{Conditional Yang-Mills measure}

Let $L_1,\ldots,L_q$ be disjoint simple loops of $\Gamma^*$ whose image is
either a component of $\partial M$ or contained in the interior of $M$. We
want to choose the law of $(h_{L_1},\ldots,h_{L_q})$. For this, it is enough
to be able to impose a deterministic value to each $h_{L_i}$. Let
$(x_1,\ldots,x_q)$ be an element of $G^q$. Let $\Gamma'$ denote the set of
edges of $\Gamma$ that do not appear in the decomposition of any $L_i$. We
denote by $dg'$ the product of Haar measures on $G^{\Gamma'}$. The fact that
the conditional Haar measure is not invariant by permutation of the factors on
$G^n$ leads either to a very heavy or to an elliptic notation. We will choose
the second option, except during a few lines. Suppose that
$L_1=a^{\epsilon^1_1}_{i^1_1} \ldots a^{\epsilon^1_{n_1}}_{i^1_{n_1}}, \ldots,
L_q=a^{\epsilon^q_1}_{i^q_1} \ldots a^{\epsilon^q_{n_q}}_{i^q_{n_q}}$, with
$\epsilon^i_j=\pm 1$. We denote by $d\nu_{x_1} \ldots d\nu_{x_q}dg'$ the
following measure on $G^\Gamma$:
$$d\nu_{x_1}^{n_1}(g^{\epsilon^1_1}_{i^1_1},\ldots,
g^{\epsilon^1_{n_1}}_{i^1_{n_1}}) \ldots
d\nu_{x_q}^{n_q}(g^{\epsilon^q_1}_{i^q_1},\ldots,
g^{\epsilon^q_{n_q}}_{i^q_{n_q}}) dg'.$$
With this notation, set:
\begin{equation}
Z(x_1,\ldots,x_q) = \int_{G^\Gamma} D \; d\nu_{x_1}\ldots
  d\nu_{x_q} dg',
\end{equation}
\begin{equation}
dP(x_1,\ldots,x_q)={1\over{Z(x_1,\ldots,x_q)}}  D \; d\nu_{x_1}\ldots
d\nu_{x_q} dg'.
\end{equation}

The function $Z(x_1,\ldots,x_q)$ is called conditional partition function on
$M$ with respect to $L_1,\ldots,L_q$. 

\begin{proposition} \label{partition_law}
  Choose $r\leq q$ and $x_{r+1},\ldots,x_q \in G$.  The distribution of
  $(h_{L_1},\ldots,h_{L_{r}})$ under $P(x_{r+1},\ldots,x_q)$ is equal to:
  $${{Z(x_1,\ldots,x_q)}\over {Z(x_{r+1},\ldots,x_q)}} \; dx_1 \ldots
  dx_{r},$$
  where each element $x_i$ corresponds to the loop $L_i$. In
  particular, the distribution of $(h_{L_1},\ldots,h_{L_{q}})$ under $P$ is ${1\over Z}
  Z(x_1,\ldots,x_q)\; dx_1 \ldots dx_q$.
\end{proposition}

\pf: The last part of the statement is just the case $r=q$.  Let $f$ be a
continuous function on $G^r$. We have:
\begin{eqnarray*}
&& \hskip -1cm \int_{G^\Gamma} f(h_{L_1},\ldots,h_{L_{r}})\;
dP(x_{r+1},\ldots,x_q) \; =\\ 
&=& {1\over {Z(x_{r+1},\ldots,x_q)}} \int_{G^\Gamma}
f(h_{L_1},\ldots,h_{L_{r}}) D \; d\nu_{x_{r+1}} \ldots d\nu_{x_q} dg' \\
&=& {1\over {Z(x_{r+1},\ldots,x_q)}} \int_{G^{r}}
\int_{G^\Gamma} 
f(h_{L_1},\ldots,h_{L_{r}}) D 
\; d\nu_{x_1} \ldots d\nu_{x_q} dg' \; dx_1\ldots dx_{r} \\
&=& \int_{G^{r}} f(x_1,\ldots,x_{r}) \left[ {1\over
    {Z(x_{r+1},\ldots,x_q)}} \int_{G^\Gamma} D 
\; d\nu_{x_1} \ldots d\nu_{x_q} dg' \right] \; dx_1\ldots dx_{r} \\
&=& \int_{G^{r}} f(x_1,\ldots,x_{r}) {{Z(x_1,\ldots,x_q)}\over {Z(x_{r+1},\ldots,x_q)}}\;
dx_1\ldots dx_{r} .
\end{eqnarray*}
\vskip-.8cm \qed

\begin{corollary} \label{discrete_disintegration}
The map $(x_1,\ldots,x_q) \mapsto
P(x_1,\ldots,x_q)$ is a disintegration of the measure $P$ with respect
to the random variable $(h_{L_1},\ldots,h_{L_q})$. This means that \\
\indent 1. $(h_{L_1},\ldots,h_{L_q})=(x_1,\ldots,x_q)$
 $P(x_1,\ldots,x_q)$-a.s. \\
\indent 2. denoting by $\eta$ the law of $(h_{L_1},\ldots,h_{L_q})$ under
$P$, we have
$$P=\int_{G^q} P(x_1,\ldots,x_q) \; d\eta(x_1,\ldots,x_q).$$
\end{corollary}

\pf. The first part is a direct consequence of the definition of
$P(x_1,\ldots,x_q)$. A simple computation proves the second one:
\begin{eqnarray*}
\int_{G^q} P(x_1,\ldots,x_q)\; d\eta(x_1,\ldots,x_q) &=& {1\over Z}\int_{G^q}
P(x_1,\ldots,x_q) Z(x_1,\ldots,x_q) \; dx_1 \ldots dx_q \\
&=& {1\over Z}\int_{G^q} D \; \nu_{x_1} \ldots \nu_{x_q} \;dg' \; dx_1 \ldots
dx_q \\ 
&=& P.
\end{eqnarray*}
\vskip -.8cm \qed

This corollary says that the $P(x_1,\ldots,x_q)$ are really what we expected
them to be. Now, given a measure $\beta$ on $G^q$, we can choose the law of
$(h_{L_1},\ldots,h_{L_q})$ to be $\beta$ by putting the measure
$P_\beta=\int_{G^q} P(x_1,\ldots,x_q) \; d\beta(x_1,\ldots,x_q)$ on
$G^\Gamma$. 

\subsection{Gauge transformations}

Let us compute how the conditional measure $P(x_1,\ldots,x_q)$ is transformed by a gauge
transformation.

\begin{lemma} \label{gt_of_P}
Let $j$ be a discrete gauge transformation. The following
  equality holds :
$$j_*P(x_1,\ldots,x_q)=P(y_1 x_1 y_1^{-1},\ldots,y_q x_q y_q^{-1}),$$
where $y_i=j(L_i(0))$.
\end{lemma}

\pf. For sake of simplicity, let us write the proof in the case $q=1$, the
general case being exactly similar, only with heavier notations. Suppose that
$L_1=a_1\ldots a_m$. Let $f$ be a continuous function on $G^\Gamma$.
\begin{eqnarray*}
j_*P(x_1)(f) &=& {1\over{Z(x_1)}} \int_{G^\Gamma} f\circ j \prod_{F\in
  \f(\Gamma)} p_{\sigma(F)}(h_{\partial F}) \; \nu_{x_1}^m(g_1,\ldots,g_m)
  dg_{m+1}\ldots dg_r \\
&=& {1\over{Z(x_1)}} \int_{G^\Gamma} f \prod_{F\in
  \f(\Gamma)} p_{\sigma(F)}(h_{\partial F}\circ j^{-1}) \;
 \nu_{x_1}^m(j(a_1(1))g_1j(a_1(0))^{-1},\ldots \\ 
&& \hskip 3.6cm \ldots,j(a_m(1))g_mj(a_m(0))^{-1})
  \; dg_{m+1}\ldots dg_r \\
&=& {1\over{Z(x_1)}} \int_{G^\Gamma} f \prod_{F\in
  \f(\Gamma)} p_{\sigma(F)}(h_{\partial F}) \;
  \nu_{j(L_1(0)) x_1
  j(L_1(0))^{-1}}^m(g_1,\ldots,g_m)\\
&& \hskip 8cm  dg_{m+1}\ldots dg_r \\
&=& {{Z(y_1 x_1 y_1^{-1})}\over{Z(x_1)}} P(y_1 x_1
  y_1^{-1}) (f),
\end{eqnarray*}

with $y_1=j(L_1(0))$. Setting $f$ to be identically equal to $1$, we get
$Z(y_1 x_1 y_1^{-1})=Z(x_1)$. \qed

Let us state the invariance property of the partition function that we just proved:

\begin{proposition} \label{inv_part}
For any $y_1,\ldots, y_q$ in $G$, one has 
$$Z(y_1^{-1} x_1 y_1, \ldots, y_q^{-1} x_q y_q) = Z(x_1,\ldots,x_q).$$
\end{proposition}

According to this result, the conditional partition function can also be
viewed as a function on $(G/\Ad)^q$, where $Ad$ is the adjoint action on $G$
given by $\Ad(y)x=y^{-1}xy$. We will use this point of view in the next
paragraph.

It is clear now that $P(x_1,\ldots,x_q)$ is not gauge invariant in
general. We will explain how to overcome this problem. 

Let $t$ be an element of the quotient $G/\Ad$, that is, a conjugacy class in
$G$. Let $x$ be an element of this class. The measure $\int_G
\delta_{yxy^{-1}} \; dy$ does not depend on the choice of $x$ in $t$. We shall
denote it by $\delta_t$. Similarly, we denote $\delta_{t_1}\otimes \ldots
\otimes \delta_{t_q}$ by $\delta_{t_1,\ldots,t_q}$. Set
$$P(t_1,\ldots,t_q)=\int_{G^q} P(x_1,\ldots,x_q) \;
\delta_{t_1,\ldots,t_q}(x_1,\ldots,x_q).$$

\begin{proposition} 
The measure $P(t_1,\ldots,t_q)$ is gauge invariant.
\end{proposition}

\pf: Let $j$ be a discrete gauge transformation and set $y_i=j(L_i(0))$.
According to the lemma \ref{gt_of_P}, we have:
\begin{eqnarray*}
j_* P(t_1,\ldots,t_q) &=& \int_{G^q} j_* P(x_1,\ldots,x_q) \; \delta_{t_1,\ldots,t_q}(x_1,\ldots,x_q)\\
&=& \int_{G^q}  P(y_1 x_1 y_1^{-1},\ldots, y_q x_q y_q^{-1}) \; \delta_{t_1,\ldots,t_q}(x_1,\ldots,x_q)\\
&=& \int_{G^q}  P(y_1 z_1 x_1^0 z_1^{-1} y_1^{-1},\ldots,y_q z_q
x_q^0 z_q^{-1} y_q^{-1})\; dz_1 \ldots dz_q \\
&=& \int_{G^q} P(z_1^{-1} x_1^0 z_1,\ldots, z_q^{-1} x_q^0 z_q) \;
dz_1 \ldots dz_q \\
&=& \int_{G^q} P(x_1,\ldots,x_q) \;
\delta_{t_1,\ldots,t_q}(x_1,\ldots,x_q) \\
&=&P(t_1,\ldots,t_q),
\end{eqnarray*}
where each $x_i^0$ was an arbitrary element of $t_i$.\qed

It will emerge later that the measures $P(t_1,\ldots,t_q)$ are in fact more
natural than the $P(x_1,\ldots,x_q)$. 

\section{Invariance by subdivision}
\label{inv-by-sub}

The invariance by subdivision is the main feature of the discrete theory. It
allows to prove that the law of the discrete holonomy along given loops does
not depend on the graph in which one computes it.

The fact that the heat kernel $(p_t)_{t>0}$ is a convolution semi-group will
play a central role in the proof. This means that for any $x\in G$ and any
$s,t$ such that $0<s<t$,
$$\int_G p_s(xy^{-1}) p_t(y) \; dy = p_t(x).$$

Let $\Gamma_1$ and $\Gamma_2$ be two graphs on $M$. Suppose that $\Gamma_2$ is
finer than $\Gamma_1$ and set $\Gamma_1=\{a_1.\ldots,a_r\}$. By definition,
each edge $a_i$ of $\Gamma_1$ is a path in $\Gamma_2^*$ and it gives rise to a
function $h_{a_i}: G^{\Gamma_2} \lra G$. The $r$-uple of those functions
constitutes a single function $(h_{a_1},\ldots,h_{a_r}):G^{\Gamma_2}\lra
G^r=G^{\Gamma_1}$ that we denote by $f_{\Gamma_1 \Gamma_2}$. the invariance by
subdivision is expressed by the following result :

\begin{theorem} \label{invariance_by_subdivision} Let $\Gamma_1$ and
  $\Gamma_2$ be two graphs on $M$ such that $\Gamma_2$ is finer than
  $\Gamma_1$. Let $L_1,\ldots,L_q$ be disjoint simple loops of
  $\Gamma_1^*$. Let $x_1,\ldots,x_q$ be elements of $G$. Then \\
\indent  1. The map $f_{\Gamma_1 \Gamma_2}:G^{\Gamma_2} \lra G^{\Gamma_1}$ is
  surjective.\\
\indent  2. This map satisfies : $\left(f_{\Gamma_1
      \Gamma_2}\right)_*P^{\Gamma_2}(x_1,\ldots,x_q)=P^{\Gamma_1}(x_1,\ldots,x_q).$
\end{theorem}

From now on, it will be sometimes necessary to write explicitely the graph in
which we consider objects such as $P$, $Z$, $D$.

We begin by proving that it is always possible to go from one graph to a finer
graph by a finite sequence of elementary transformations. 

\begin{lemma} Let $\Gamma$ and $\Gamma'$ be two graphs such that $\Gamma <
  \Gamma'$. There exist an increasing sequence of graphs $\Gamma_0=\Gamma
  <\Gamma_1<\ldots<\Gamma_n<\ldots$, stationnary of limit $\Gamma'$ and such
  that for any nonnegative $n$, one can transform $\Gamma_n$ into
  $\Gamma_{n+1}$ by one of the two following elementary operations:\\
(V) Add a vertex to $\Gamma_n$, i.e. replace an edge $a$ by two edges $b$ and
  $c$ such that $a=bc$,\\
(E) Add an edge to $\Gamma_n$, this new edge joining two vertices of
  $\Gamma_n$.
\end{lemma}

\pf. We proceed by induction on $n$. $\Gamma_0$ is given, equal to $\Gamma$.
Suppose $\Gamma_n$ given, with $\Gamma_n<\Gamma'$. Recall that $\v(\Gamma)$
denotes the set of vertices of $\Gamma$.\\
$\bullet$ We have $\v(\Gamma_n)\subset \v(\Gamma')\cap\supp(\Gamma_n)$. If
this inclusion is a strict one, pick an element of
$(\v(\Gamma')\cap\supp(\Gamma_n))\backslash \v(\Gamma_n)$. It is a vertex of
$\Gamma'$ which is on an edge of $\Gamma_n$ whithout being one of its end
points. By an operation (V), we add this vertex to $\Gamma_n$ and get
$\Gamma_{n+1}$ which is still finer than $\Gamma'$. Note that
$\card(\Gamma_{n+1})=\card(\Gamma_n)+1$. \\
$\bullet$ If $\v(\Gamma_n)=\v(\Gamma')\cap\supp(\Gamma_n)$, then each edge of
$\Gamma_n$ is an edge of $\Gamma'$. In other words, $\Gamma_n \subset
\Gamma'$. If this inclusion is a strict one, there exists an edge of $\Gamma'$
which is not an edge of $\Gamma_n$ and by connectedness of $\Gamma'$ we may
assume that this edge has at least one of its end points on $\supp(\Gamma_n)$.
By an operation (E), we add this edge to $\Gamma_n$ and get $\Gamma_{n+1}$
which is still connected and finer than $\Gamma'$. The pregraph $\Gamma_{n+1}$
is a graph. Indeed, we just noticed that it is connected and it is finer than
$\Gamma_0$, so that its support contains $\partial M$. It can happen that the
operation (E) cuts a face in two pieces, but they are still diffeomorphic to
disks. We also have $\card(\Gamma_{n+1})=\card(\Gamma_n)+1$.\\
$\bullet$ If $\Gamma_n=\Gamma'$, just set $\Gamma_{n+1}=\Gamma_n$.\\

At each step, the fact that $\Gamma_n$ is a graph implies that $\Gamma_{n+1}$
is also a graph : connectedness is preserved, as well as boundary
properties. The faces of a graph are not modified by an operation (V) and it
can happen that an operation (E) cuts a face into two pieces, which are still
diffeomorphic to disks.

For each $n$, $\Gamma_n<\Gamma'$ implies
$\card(\Gamma_n)\leq\card(\Gamma')$. On the other hand, elementary operations
increase strictly the cardinal of the graph. Thus, there is necessarily only a
finite number of such operations before the sequence becomes stationnary. \qed

\begin{figure}[hbtp]
\begin{center}
\input{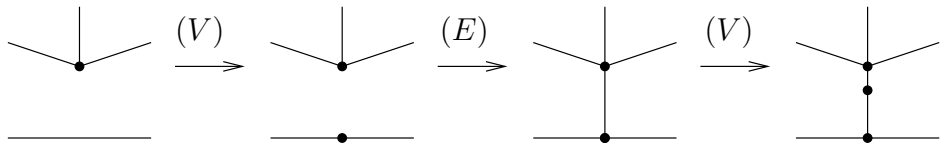}
\end{center}
\caption{Examples of elementary transformations of a graph.
\label{elem}}
\end{figure}

\begin{lemma} Let $\Gamma_1<\Gamma_2<\Gamma_3$ be three graphs. Then
  $$f_{\Gamma_1 \Gamma_3}=f_{\Gamma_1 \Gamma_2} \circ f_{\Gamma_2 \Gamma_3}.$$
\end{lemma}

\pf. It is the associativity of the product in $G$. \qed

This lemma shows that it is enough to prove the theorem
\ref{invariance_by_subdivision} when $\Gamma_2$ can be deduced from $\Gamma_1$
by an elementary operation. One recovers the general case by composition.

During the proof, we set $\beta=\delta_{(x_1,\ldots,x_q)}$ and use the
notations $P_\beta=P(x_1,\ldots,x_q)$ and $Z_\beta=Z(x_1,\ldots,x_q)$ in order
to make the expressions shorter.

\pf {\textsl{of theorem}} \ref{invariance_by_subdivision} : $1.$ We prove that
$f_{\Gamma_1 \Gamma_2}$ is surjective. If $\Gamma_2$ can be deduced from
$\Gamma_1$ by an operation (E), $f_{\Gamma_1 \Gamma_2}$ is just the projection
that forgets the factor associated with the new edge. It is of course
surjective. In the case of an operation (V), $f_{\Gamma_1 \Gamma_2}$ preserves
all factors except those associated with the two new edges, that are
multiplied. It is also surjective.

$2.$ Let us begin by the case of an operation of type (E). We fix some
notations. Set $\Gamma_1=\{a_1,\ldots,a_r\}$ and
$\Gamma_2=\{a_1,\ldots,a_r,b\}$. The new edge $b$ is located in a face $F_0$
of $\Gamma_1$. Two situations are possible : either $b$ has one end point on
$\partial F_0$ or it has both. In the first case, $F_0$ is still a face of
$\Gamma_2$, with a new factor $bb^{-1}$ in its boundary. In the second case,
$F_0$ is cut into two faces $F_1$ and $F_2$ by $b$. Let us consider this
second case. The boundaries of $F_0$, $F_1$ and $F_2$ can be written
respectively $\partial F_0=c_1 c_2$, $\partial F_1=c_1 b^{-1}$ and $\partial
F_2 = b c_2$. Let $f$ be a continuous function on $G^{\Gamma_1}$.
\begin{eqnarray*}
&& \hskip -.5cm \int_{G^{\Gamma_1}} f \; d\left(\left(f_{\Gamma_1 \Gamma_2}\right)_*
  P_\beta^{\Gamma_2}\right) = {1\over {Z_\beta^{\Gamma_2}}}
  \int_{G^{\Gamma_2}} f(g_1,\ldots,g_r) p_{\sigma(F_1)}(g_{r+1}^{-1}h_{c_1})
  p_{\sigma(F_2)}(h_{c_2}g_{r+1}) \\
&& \hskip 5.7cm \prod_{F\in\f(\Gamma_1)\backslash F_0}
  p_{\sigma(F)}(h_{\partial F}) \; 
  d\nu_{x_1} \ldots d\nu_{x_q} dg',
\end{eqnarray*}
where $g_{r+1}$ is the element associated with $b$. Since the $L_i$'s are
paths in $\Gamma_1$, the new edge $b$ is not involved in their decomposition.
Thus we can isolate $dg_{r+1}$ in $dg'$ and integrate against it. We use the
fact that the heat kernel is a convolution semi-group. We get the following
expression :
\begin{eqnarray*}
&& \hskip -.7cm ={1\over {Z_\beta^{\Gamma_2}}} \int_{G^{\Gamma_1}}
f(g_1,\ldots,g_r) 
p_{\sigma(F_1)+\sigma(F_2)}(h_{c_1}h_{c_2}) \prod_{F\in\f(\Gamma_1)\backslash
  F_0} p_{\sigma(F)}(h_{\partial F}) \; d\nu_{x_1} \ldots d\nu_{x_q} dg'\\
&& \hskip -.7cm = {{Z^{\Gamma_1}_\beta}\over{Z^{\Gamma_2}_\beta}}
\int_{G^{\Gamma_1}} f \; 
dP^{\Gamma_1}_\beta.
\end{eqnarray*}
Setting $f=1$, we get $Z^{\Gamma_1}_\beta=Z^{\Gamma_2}_\beta$ and the result. 

The case where the new edge does not cut $F_0$ into two faces is even simpler:
the factor $bb^{-1}$ vanishes in all computations, because $f$ does not depend
on the factor associated to $b$.

Now consider the case of an operation (V). Set $\Gamma_1=\{ a_1,\ldots,a_r\}$
and $\Gamma_2=\{b,c,a_2,\ldots,a_r\}$, with $a_1=bc$. The edge $a_1$ can be on
the boundary of one or two faces, depending on the fact that it is on
$\partial M$ or not. 
\begin{eqnarray*}
&& \hskip -.5cm \int_{G^{\Gamma_1}} f \; d\left(\left(f_{\Gamma_1
  \Gamma_2}\right)_* 
  P_\beta^{\Gamma_2}\right) = {1\over{Z_\beta^{\Gamma_2}}} \int_{G^{\Gamma_2}}
f(g_c g_b,g_2,\ldots,g_r) D^{\Gamma_1}(g_c g_b,g_2,\ldots,g_r)\\
&& \hskip 9.5cm  \nu_{x_1}
\ldots \nu_{x_q} dg',
\end{eqnarray*}
where $(g_b,g_c,g_2,\ldots,g_r)$ denotes the generic element of
$G^{\Gamma_2}$. We have to discuss two cases : either $a_1$ is involved in the
decomposition of one of the $L_i$'s, say $L_1$, or it is not. 
If it is not, we can isolate $dg_b dg_c$ in the $dg'$ term. By integrating
against $dg_b$, the dependence in $g_c$ disappears by right invariance of
$dg_b$ and we get
$$\ldots={1\over{Z_\beta^{\Gamma_2}}} \int_{G^{\Gamma_1}} f(g,g_2,\ldots,g_r)
D^{\Gamma_1}(g,g_2,\ldots,g_r) \; \nu_{x_1} \ldots \nu_{x_q} dg' =
{{Z^{\Gamma_1}_\beta}\over{Z^{\Gamma_2}_\beta}} \int_{G^{\Gamma_1}} f \;
dP^{\Gamma_1}_\beta,$$
and we conclude as before.  If $a_1$ is involved in the
decomposition of $L_1$, we can suppose that $L_1=a_1\ldots a_m$, with $m\geq
2$. We write $d\nu_{x_1}$ in a convenient way, putting the contioning on
$g_m$, which is necessarily distinct from $g_1$. We get :
\begin{eqnarray*}
&& \hskip -.5cm {1\over{Z_\beta^{\Gamma_2}}} \int_{G^{\Gamma_1}} f(g_c g_b,
g_2, \ldots,
\widetilde{g_m},\ldots,g_r) D^{\Gamma_2}(g_c g_b,g_2,\ldots,\widetilde{g_m},
\ldots,g_r) \\
&& \hskip 6.2cm  dg_b dg_c dg_2 \ldots dg_{m-1} \nu_{x_2}\ldots \nu_{x_q} dg',
\end{eqnarray*}
where $\widetilde{g_m}=x_1(g_{m-1}\ldots g_2 g_c g_b)^{-1}$. This is equal to 
\begin{eqnarray*}
&& \hskip -.5cm {1\over{Z_\beta^{\Gamma_2}}} \int_{G^{\Gamma_1}} f(g, g_2,
\ldots, 
\widetilde{g_m},\ldots,g_r) D^{\Gamma_2}(g,g_2,\ldots,\widetilde{g_m},
\ldots,g_r) \\
&& \hskip 7cm  dg dg_2 \ldots dg_{m-1} \nu_{x_2}\ldots \nu_{x_q} dg',
\end{eqnarray*}
with $\widetilde{g_m}$ now
equal to $x_1(g_{s_1}\ldots g_2 g)^{-1}$. This is one more time equal to
$${{Z^{\Gamma_1}_\beta}\over{Z^{\Gamma_2}_\beta}} \int_{G^{\Gamma_1}} f \;
dP^{\Gamma_1}_\beta$$
and we get the result. \qed

\begin{corollary} Let $c_1,\ldots,c_n$ be paths that are simultaneously
  elements of $\Gamma_1^*$ and $\Gamma_2^*$, where $\Gamma_1$ and $\Gamma_2$
  are two graphs such that $\Gamma_1 < \Gamma_2$. Then the law of the discrete
  holonomy along $c_1,\ldots,c_n$ is the same on $\Gamma_1$ and $\Gamma_2$. In
  other words, the law of $(h^{\Gamma_1}_{c_1},\ldots,h^{\Gamma_1}_{c_n})$ on
  $(G^{\Gamma_1},P_\beta^{\Gamma_1})$ and the law of
  $(h^{\Gamma_2}_{c_1},\ldots,h^{\Gamma_2}_{c_n})$ on
  $(G^{\Gamma_2},P_\beta^{\Gamma_2})$ are equal.
\end{corollary}

\pf. It is enough to verify that
$(h^{\Gamma_1}_{c_1},\ldots,h^{\Gamma_1}_{c_n}) \circ f_{\Gamma_1 \Gamma_2} = 
(h^{\Gamma_2}_{c_1},\ldots,h^{\Gamma_2}_{c_n})$. This is true if the $c_i$'s
are edges of $\Gamma_1$, thus it is true in general by multiplicativity. \qed

During the proof of the theorem, we also proved an important result about
the conditional partition function:

\begin{proposition} \label{Z_constant}
  Let $\Gamma_1$ and $\Gamma_2$ be two graphs such that $\Gamma_1 <\Gamma_2$.
  Take $L_1,\ldots, L_q$ and $x_1,\ldots,x_q$ as usual. Then
  $$Z^{\Gamma_1}(x_1,\ldots,x_q)=Z^{\Gamma_2}(x_1,\ldots,x_q).$$
\end{proposition}

Let us discuss briefly the choice of the heat kernel in the definition of
$P$. This choice is the key of the physical relevance of the theory. It is a
physicist, A. Migdal \cite{Migdal}, who suggested first to use the heat kernel
in the mathematical formulation of the theory. Nevertheless, it is possible to
construct a discrete theory using any other convolution semigroup. For
example, Albeverio, H{\o}egh Krohn and Holden investigated some properties of
the random fields obtained this way \cite{Albeverio}. But it would probably be
much more difficult to construct a continuous theory without the nice
regularity properties that characterize the heat kernel among all other
convolution semigroups. 

\section{Invariance by area-preserving diffeomorphisms}

The manifold $M$ is given with its differentiable structure and the Lebesguian
surface measure $\sigma$. Let $\Gamma$ be a graph on $M$ and $\phi:M \lra M$ a
diffeomorphism such that $\phi_*\sigma=\sigma$. Then $\phi$ transforms
$\Gamma$ into a graph $\phi(\Gamma)$ and induces a bijection between faces of
$\Gamma$ and $\phi(\Gamma)$ that preserves the surface. Thus, the natural
bijection induced between $G^\Gamma$ and $G^{\phi(\Gamma)}$ preserves the
discrete Yang-Mills measure.  Let us state in a slightly more general way this
invariance property.

\begin{proposition} \label{inv_by_apd}
  Let $(M_1,\sigma_1)$ and $(M_2,\sigma_2)$ be two surfaces and $\phi:M_1 \lra
  M_2$ a diffeomorphism such that $\phi_* \sigma_1= \sigma_2$. Let $\Gamma_1$
  be a graph on $M_1$ and $\Gamma_2=\phi(\Gamma_1)$ the corresponding graph on
  $M_2$. Still denoting by $\phi:G^{\Gamma_1} \lra G^{\Gamma_2}$ the induced
  bijection, one has
$$\phi_* P^{\Gamma_1}(x_1,\ldots,x_q) = P^{\Gamma_2}(x_1,\ldots,x_q).$$
\end{proposition}

Thus, for each family $(c_1,\ldots,c_n)$ of paths in $\Gamma_1^*$, the
discrete holonomy along $(c_1,\ldots,c_n)$ and the
discrete holonomy along $(\phi(c_1),\ldots,\phi(c_n))$ have the same
distribution.

\section{Examples} 

In this section, we will compute the law of the discrete holonomy in two basic
situations. 

\subsection{Holonomy along an open path}
\label{hol_edge}
Let $\Gamma$ be a graph and $c \in \Gamma^*$ be an open path, i.e. a path such
that $c(0)\neq c(1)$. Let $L_1,\ldots,L_q$ be disjoint
simple loops of $\Gamma^*$ and $t_1,\ldots,t_q$ be elements of $G/\Ad$. Let us
compute the law of $h_c$ under $P(t_1,\ldots,t_q)$. We will use the gauge
invariance of $P(t_1,\ldots,t_q)$.

Let $f$ be a continuous function on $G$ and $j$ a discrete gauge
transformation. Recall from the proposition \ref{gt_of_P} that
$j_*P(t_1,\ldots,t_q)=P(t_1,\ldots,t_q)$, so that:
\begin{eqnarray*}
\int_{G^\Gamma} f(h_c) \; dP(t_1,\ldots,t_q) &=& \int_{G^\Gamma} f(h_c\circ j)\;
dP(t_1,\ldots,t_q)\\ 
&=& \int_{G^\Gamma} f(j(c(1))^{-1} h_a j(c(0)))\; dP(t_1,\ldots,t_q).
\end{eqnarray*}
Thus the law of $h_c$ is right and left invariant on $G$ : it is the Haar
measure. 

\subsection{Holonomy along the boundary of a small disk}
\label{small_disk}
Let $\Gamma$ be a graph on $M$, $L_1,\ldots,L_q$ be disjoint simple loops of
$\Gamma^*$. Let $l$ be a loop of $\Gamma^*$ which is the boundary of a disk
$D$ such that $L_i([0,1])$ is not constained in $\overline{D}$ for each $i$.
We will estimate the law of $h_l$ on $(G^\Gamma,P_\beta)$.

Let $\rho$ be the function defined on $G$ by $\rho(x)=d(1,x)$, where $d$ is
the biinvariant Riemannian distance. We want to estimate the size of
$\rho(h_l)$.

$$\int_{G^\Gamma} \rho(h_l)\; dP_\beta = {1\over {Z_\beta}} \int_{G^\Gamma}
\rho(h_l) \prod_{F\in \f(\Gamma)} p_{\sigma(F)}(h_{\partial F}) \;
d\nu_{x_1}\ldots d\nu_{x_q} dg'.$$
We need a result about graphs.

\begin{lemma} Let $M$ be a surface and $\Gamma$ a graph on $M$. There exists a
  subgraph of $\Gamma$ which has only one face.
\end{lemma}

\pf. If a graph has more than one face, there exists an edge which is on the
boundary of two different faces. If we remove this edge, the resulting
pregraph is still a graph. Indeed, this removed edge was necessarily in the
interior of $M$, thus boundary properties are preserved. The faces of the new
pregraph are those of the old one, except two faces that were glued along a
segment. So the new face is diffeomorphic to a disk and by lemma
\ref{disks_simply_connected}, we know that the new pregraph is connected. In a
finite number of such steps, one gets a subgraph of $\Gamma$ with only one
face. \qed

The pregraph constituted by the $L_i$'s cuts $M$ into several pieces
homeomorphic to surfaces with boundary $M_1,\ldots,M_k$. The graph $\Gamma$
induces a graph on each $M_i$. By the preceding lemma, there exists a subgraph
$\Gamma'$ of $\Gamma$ that has exactly one face on each $M_i$ and such that
$L_1,\ldots,L_q \in \Gamma'^*$. Now add to $\Gamma'$ the edges required to
form $l$ and, if necessary, a simple path connecting $l$ with
$\supp(\Gamma')$. Finally, the assumption that the image of $L_i$ is never
contained in $\overline{D}$ shows that it is possible, maybe by adding some
vertices to $\Gamma'$, to be sure that each $L_i$ has an edge outside $D$. We
get a graph $\Gamma''$ which is included in $\Gamma$ and in which we will
compute, using the invariance by subdivision. We use the notation $P_\beta$
and $Z_\beta$.
\begin{eqnarray*}
\int_{G^\Gamma} \rho(h_l) \; dP_\beta^\Gamma &=& \int_{G^{\Gamma''}}
\rho(h_l)\; dP_\beta^{\Gamma''} \\
&=& {1\over {Z_\beta}} \int_{G^{\Gamma''}} \rho(h_l) \prod_{F\in\f(\Gamma'')}
p_{\sigma(F)}(h_{\partial F}) \; d\nu_{x_1}\ldots d\nu_{x_q} dg'.
\end{eqnarray*}

\begin{lemma} The function $t \lra \parallel p_t \parallel_\infty$ is
  decreasing on $(0,\infty)$.
\end{lemma}

\pf. Let $0<s<t$ be two positive times. Let $x$ be an element of $G$. We can
estimate $p_t(x)$ in the following way, keeping in mind that $p_t$ is a
positive function on $G$ :
$$p_t(x) = \int_G p_s(xy^{-1}) p_{t-s}(y) \; dy \leq \parallel p_s
\parallel_\infty \int_G p_{t-s}(y) \; dy \leq \parallel p_s
\parallel_\infty.$$ 
\vskip -.9cm \qed

Recall that each $L_i$ has at least an edge outside $\overline{D}$. For
each $L_i$, we put the conditioning of $d\nu_{x_i}$ on one of these edges. 

On the other hand, each face $F$ of $\Gamma''$ which is not included in $D$ is
included in a face of $\Gamma'$, i.e. in a $M_i$: its surface is greater than
$\sigma(M_i)-\sigma(D)$.  We assume that $\sigma(D) \leq {1\over 2}
\inf_i \sigma(M_i)$.  Then
$$
\int_{G^{\Gamma''}} \rho(h_l)\; dP_\beta \leq {1\over {Z_\beta}} \prod_{F
  \not\subset D} \sup_{x\in G} \left| p_{\sigma(F)}(x) \right|
\int_{G^{\Gamma''}} \rho(h_l) \prod_{F\subset D} p_{\sigma(F)}(h_{\partial
  F})\; dg' $$
The last integral is nothing but an integral against the
discrete Yang-Mills measure on $D$, which is a surface with boundary. Using
the invariance by subdivision inside $D$, we can replace the graph induced by
$\Gamma''$ by a very simple graph whose support is just $\partial D=l([0,1])$.
This leads to
\begin{equation} \label{estim0}
\int_{G^\Gamma} \rho(h_l) \; dP_\beta^\Gamma \leq {1\over {Z_\beta}} \prod_{i=1}^k \parallel p_{{\sigma(M_i)}\over 2}
\parallel_\infty  \int_G \rho(g) p_{\sigma(D)}(g)\; dg.
\end{equation}

We are led to a problem of estimation of the heat kernel at small time. 

\begin{lemma} \label{estimate_heat_kernel} The following estimates hold :
  $$\int_G \rho(g)^4 p_t(g) \; dg = O(t^2),\;\;\;\;\; \int_G \rho(g)^2 p_t(g)
  \; dg = O(t),\;\;\;\;\; \int_G \rho(g) p_t(g) \; dg = O(\sqrt{t}).$$
\end{lemma}

Let $d$ denote the dimension of $G$. We use the following result proved in
\cite{Varopoulos}(V.4.3): 

\begin{proposition} There exists a positive constant $C$ such that for all $t
  \in (0,1)$, all $g \in G$, 
  $${1\over {C}} t^{-{d\over 2}} e^{-{{C\rho(g)^2}\over t}} \leq p_t(g) \leq
  C t^{-{d\over 2}} e^{-{{\rho^2(g)}\over {Ct}}}.$$
\end{proposition}

\pf. The first estimate implies both others. We use normal coordinates at the
identity of $G$. Let $D_R$ be a geodesic disk of radius $R$ around $1$, with
$R$ such that $\exp$ is a diffeomorphism  from $B(0,R) \subset T_1G$ onto
$D_R$. We cut the integral according to $G=D_R \cup D_R^c$. One $D_R^c$, we
have:
$$\int_{D_R^c} \rho(g)^4 p_t(g) \; dg \leq C t^{-{d\over 2}} \diam(G)^4
e^{-{{R^2}\over {Ct}}} \leq C_1 t^{-{d\over 2}} e^{-{{R^2}\over {Ct}}}.$$
For
the part corresponding to $D_R$, we use spherical coordinates $(r,\theta)$ on
$\exp^{-1}(D_R)=B(0,R)$. Note that on $B(0,R)$, the image of the Haar measure
by $\exp^{-1}$ can be compared to the Lebesgue measure, so:
\begin{eqnarray*}
\int_{D_R} \rho(g)^4 p_t(g) \; dg &\leq& C_2 \int_{[0,R]\times S^{d-1}}
  r^4 p_t(\exp(r,\theta)) r^{d-1}\;dr d\theta \\
&\leq & C_3 t^{-{d\over 2}} \int_0^R r^{d+3} e^{-{{r^2}\over
  {Ct}}} \; dr\\
&\leq & C_3 t^{2+{d\over 2}} t^{-{d\over 2}} \int_0^\infty {{r^{d+3}}\over
  {t^{{{d+3}\over 2}}}} e^{-{{r^2}\over {Ct}}} \; {{dr}\over {t^{1\over 2}}}
  \\
&\leq& C_4 t^2.
\end{eqnarray*}
This estimation remains true if we replace $R$ by $R'<R$. Thus, for $t$ small
enough, we have:
\begin{eqnarray*}
\int_G \rho(g)^4 p_t(g) \; dg &=& \int_{D_{t^{1/4}}} \rho(g)^4 p_t(g) \; dg
+\int_{D_{t^{1/4}}^c} \rho(g)^4 p_t(g) \; dg\\
&\leq &  C_4 t^2 + C_1 t^{-{d\over 2}} e^{-{{1\over {Ct^{1/2}}}}} = O(t^2).
\end{eqnarray*}
\vskip -.8cm \qed

Finally, we deduce from relation \ref{estim0} and the preceding lemma the
following proposition:
\begin{proposition} \label{fundamental_estimate}
  Let $\Gamma$ be a graph on $M$. Let $L_1,\ldots,L_q$ be disjoint simple
  loops of $\Gamma^*$ and $x_1,\ldots,x_q$ be elements of $G$. Let $l$ be the
  boundary of a disk $D$ such that none of the $L_i$'s has its image contained
  in $\overline{D}$. There exist two positive constants $s$ and $C$ depending
  on the $L_i$'s but not on the $x_i$'s such that if $\sigma(D)\leq s$, then
$$\int_{G^\Gamma} \rho(h_l) \; dP(x_1,\ldots,x_q) \leq C \sqrt{\sigma(D)}.$$
\end{proposition}

This regularity property will play an essential role in the construction of
the continuous measure.

\section{Discrete Abelian theory}
\label{Abelian_theory_1}

\subsection{Decomposition of cycles}
\label{decomposition_of_cycles}

Until now, we only used the compactness of $G$. We will finish this first
chapter with a detailed study of the case $G=U(1)$. All results could be
extended without conceptual problems to the case $G=U(1)^n$, i.e. the general
compact Abelian case, but this would also make the notations much heavier.

We fix $M$, $\sigma$ as usual and a graph $\Gamma$ on $M$. Our aim is to
analyze the law of the family $(h_c)_{c\in \Gamma^*}$. Set
$\{a_1,\ldots,a_r\}=\Gamma$.  Since $G$ is Abelian, the function $h_c:G^\Gamma
\lra G$ associated with a path $c$ depends only on the number of occurences of
each $a_i$ in the decomposition of $c$, not on the order of the edges in this
decomposition. In other words, the function $h_c$ depends only on the image of
$c$ by the natural morphism of monoids $\Gamma^* \lra \Z^\Gamma$ which sends
$a_i$ to $(0,\ldots,1,\ldots,0)$ with a $1$ at the $i$-th place. Conversely,
each element of $\Z^\Gamma$ determines without ambiguity a function from
$U(1)^\Gamma$ into $U(1)$.

So, the natural index space in this context is $\Z^\Gamma$ instead of
$\Gamma^*$ and this allows to consider linear combination of paths.  Let us
denote by $C\Gamma \subset \Z^\Gamma$ the set of linear combination of loops,
also called cycles. We are especially interested in the law of $(h_c)_{c \in
  C\Gamma}$. The reason for which we consider only loops will become clear at
the end of chapter 2. Basically, it is because for an arbitrary $G$, the
holonomy along an open path is not a gauge-invariant function of a connection.

Let us recall a classical result about the homology of $M$. 

\begin{theorem} \label{description_of_H1} Let $g$ be the genus of $M$ and $p$
  the number of connected components of $\partial M$. Then
$$H_1(M;\Z)\simeq \cases{\Z^{2g} \;\;{\rm if}\;\; p=0 \cr
     \Z^{2g+p-1} \;\; {\rm if }\;\; p>0.}$$
\end{theorem}

If $p>0$, one can construct a system of loops representing a basis of $H_1(M)$
by taking $p-1$ components of $\partial M$ and $2g$ loops of $M$ that generate
the $H_1$ of a minimal closure of $M$, i.e. a surface obtained from $M$ by
gluing a disk along each boundary component.

So, let us choose such a system composed by $\ell_1,\ldots,\ell_{2g}$ in
$\Gamma^*$ and $p-1$ loops $N_1,\ldots, N_{p-1}$ that we denote just as the
corresponding boundary components, with an abuse of notation. We can obtain
the $\ell_i$'s by deforming an arbitrary system of generators using the same
technique as in the proof of the proposition \ref{disks_imply_H1}.

Now let $c$ be a cycle in $C\Gamma$. There is an unique decomposition
$$c=\lambda_1 \ell_1 +\ldots + \lambda_{2g} \ell_{2g} + \nu_1 N_1 +\ldots +
\nu_{p-1} N_{p-1} + c^\perp,$$
with $\lambda_i, \nu_j \in \Z$ and $c^\perp \in
C\Gamma$ a cycle homologous to zero. Let us denote by $C_0\Gamma$ the
submodule of $C\Gamma$ spanned by the cycles homologous to zero.

\begin{proposition} \label{faces_are_basis} If $\partial M$ is empty
  (resp. non empty), the boundaries of all faces except one chosen arbitrarily
  (resp. of all faces) form a basis of the submodule $C_0\Gamma$ of
  $C\Gamma$. 
\end{proposition}

We will prove this proposition very soon. Set $\f(\Gamma)=\{F_1,\ldots,F_n\}$
and choose for each $F_i$ a cycle $\partial F_i$ whose image is the boundary
of $F_i$. We can write :
\begin{equation}\label{decomp}
c=\lambda_1 \ell_1 +\ldots + \lambda_{2g} \ell_{2g} + \nu_1 N_1 +\ldots +
\nu_{p-1} N_{p-1} + \mu_1 \partial F_1 +\ldots + \mu_n \partial F_n,
\end{equation}
the decomposition being non unique if $M$ is closed. The relation
\ref{decomp}, together with the multiplicativity of the holonomy, shows that
the law of the family $(h_c)_{c\in C\Gamma}$ is completely determined by the
law of what we will call a fundamental system :
$$(h_{\ell_1},\ldots,h_{\ell_{2g}},h_{N_1},\ldots,h_{N_{p-1}},h_{\partial
  F_1},\ldots,h_{\partial F_n}).$$

\pf \textsl{of proposition} \ref{faces_are_basis} : To begin with, suppose
that $M$ has no boundary. We proceed by induction on $n=\card \f(\Gamma)$. If
$n=1$, the only loop in $C\Gamma$ is $\partial F$ and it is homologically
trivial.

Now suppose that the result is true for a graph with $n-1$ faces. Let $\Gamma$
be a graph with $n$ faces. There is an edge of $\Gamma$, say $a_r$, which is
on the boundary of two distinct faces, say $F_{n-1}$ and $F_n$. Let
$\Gamma'=\{a_1,\ldots,a_{r-1}\}$ be the graph obtained by removing $a_r$. It
has $n-1$ faces $F_1,\ldots,F_{n-2},F_{n-1}\cup F_n$. Let $c$ be a cycle of
$C_0\Gamma$. We can decompose it uniquely in $c=c_0+p a_r$ with $p \in \Z$ and
$c_0 \in C\Gamma'$. We can also write $\partial F_{n-1} = a_r +b$ with $b\in
C\Gamma'$. So , we have $c=(c_0-pb)+p\partial F_{n-1}$.  By induction,
$c_0-pb$, which is homologous to zero in $\Gamma'$, is a linear combination of
$\partial F_1,\ldots,\partial F_{n-2}$. Thus, $\partial F_1,\ldots,\partial
F_{n-1}$ generate the submodule of homologically trivial cycles in $C\Gamma$.
On the other hand, $\partial F_1,\ldots,\partial F_{n-2}$ are linearly
indendent by induction and $\partial F_{n-1}$ is independent of the submodule
that they generate, because it contains the edge $a_r$. This gives the result
when $M$ is closed.

If $M$ has a boundary, consider a minimal closure $i_1: M \lra M_1$ of $M$ and
identify $M$ with $i_1(M)$. Let $c$ be a cycle homologous to zero in $M$. It
is also homologous to zero in $M_1$ and can be decomposed using the result on
$M_1$ into :
$$c=\sum_{F_i\in \f(\Gamma), F_i\subset M} \mu_i \partial F_i + \nu_1 N_1
+\ldots + \nu_{p-1} N_{p-1},$$
because the $N_i$'s are the boundaries of the
faces of $\Gamma$ on $M_1 - M$. This decomposition gives, in $H_1(M)$,
$$[c]=0=\nu_1 [N_1]+\ldots + \nu_{p-1}[N_{p-1}],$$
implying
$\nu_1=\ldots=\nu_{p-1}=0$ and $c=\mu_1 \partial F_1+\ldots+\mu_n\partial
F_n$. The independence of the $\partial F_i$'s on $M_1$ implies their
independence on $M$. \qed

\subsection{Study of a fundamental system}

We want to study the discrete Yang-Mills measure conditioned by the holonomies
along the boundary components of $M$. Let $x_1,\ldots,x_p$ be elements of
$U(1)$. Under $P_\beta=P(x_1,\ldots,x_p)$, the law of $(h_{N_1},\ldots,h_{N_{p-1}})$
is deterministic, equal to $\delta_{(x_1,\ldots,x_{p-1})}$.

\begin{proposition} Under the measure $\nu_{x_1} \otimes \ldots \otimes
  \nu_{x_q}\otimes dg'$ on $G^\Gamma$, the variables
  $h_{\ell_1},\ldots,h_{\ell_{2g}}$, $h_{\partial F_1},\ldots,h_{\partial F_{n-1}}$
  are uniform and independent on $U(1)$.
\end{proposition}

\pf. We compute the characteristic function of
$(h_{\ell_1},\ldots,h_{\ell_{2g}},h_{\partial F_1},\ldots,h_{\partial
  F_{n-1}})$, seen as a $\C^{2g+n-1}$-valued random variable. In order to
simplify the notations, we choose an orientation of $M$ and assume that each
$N_i\subset \partial M$ and each $\partial F_i$ is oriented according to the
usual convention. Let $\lambda_1,\ldots,\lambda_{2g},\mu_1,\ldots,\mu_{n-1}$
be integers.
\begin{eqnarray*}
F(\lambda_1,\ldots,\lambda_{2g},\mu_1,\ldots,\mu_{n-1})&=& \int_{U(1)^\Gamma}
h_{\ell_1}^{\lambda_1}\ldots h_{\ell_{2g}}^{\lambda_{2g}}h_{\partial
  F_1}^{\mu_1}\ldots h_{\partial F_{n-1}}^{\mu_{n-1}}\; dP_\beta\\
&=& \int_{U(1)^\Gamma} h_{\lambda_1 \ell_1 + \ldots +\lambda_{2g}\ell_{2g}+\mu_1
  \partial F_1+\ldots+\mu_{n-1} \partial F_{n-1}} \;dP_\beta\\
&=& \int_{U(1)^\Gamma} h_{a_1}^{\alpha_1}\ldots h_{a_r}^{\alpha_r} \;
dP_\beta, 
\end{eqnarray*}
where $\sum_i \lambda_i \ell_i + \sum_j \mu_j \partial F_j = \sum_k \alpha_k
a_k$. Suppose that the $a_i$'s are labeled in such a way that $N_1=a_1\ldots
a_{i_1},\ldots,N_p=a_{i_{p-1}+1}\ldots a_{i_p}$ with $1<a_1<\ldots<a_{i_p}$.
\begin{eqnarray*}
F(\lambda_i,\mu_j) &=& \int_{U(1)^\Gamma} (h_{a_1}^{\alpha_1}\ldots
h_{a_{i_1}}^{\alpha_{i_1}}) 
\ldots (h_{a_{i_{p-1}+1}}^{\alpha_{i_{p-1}+1}} \ldots
h_{a_{i_p}}^{\alpha_{i_p}}) h_{a_{i_{p}+1}}^{\alpha_{i_{p}+1}} \ldots 
h_{a_r}^{\alpha_r} \; dP_\beta\\
&=& \int_{U(1)^\Gamma} g_1^{\alpha_1}\ldots g_{i_1}^{\alpha_{i_1}} \;
d\nu_{x_1} \ldots \int_{U(1)^\Gamma} g_{i_{p-1}+1}^{\alpha_{i_{p-1}+1}} \ldots
g_{i_p}^{\alpha_{i_p}} \; d\nu_{x_p}\\
&& \hskip 5.2cm  \int_{U(1)^\Gamma} g^{\alpha_{i_{p}+1}}\; dg \ldots
\int_{U(1)^\Gamma} g^{\alpha_r} \; dg.
\end{eqnarray*}

This product is zero if one of the $\alpha_k$'s with $k\geq i_p+1$ is nonzero.
Otherwise, $\partial M$ is non empty and the cycle $\sum \lambda_i \ell_i + \sum
\mu_j \partial F_j$ has all its edges on $\partial M$. Thus, we have an
equality
$$\sum_{i=1}^{2g} \lambda_i \ell_i + \sum_{j=1}^{n-1} \mu_j \partial
F_j=\sum_{k=1}^p \nu_k N_k,$$ 
which, in $H_1(M)$, implies $\nu_p[N_p] = \sum \lambda_i [\ell_i] - \sum_{k<p}
\nu_k [N_k]$. Since $[N_p]=-[N_1]-\ldots-[N_{p-1}]$, this implies
$\lambda_i=0$ for all $i$ and $\nu_k=\nu_p$ for all $k$. We get
$$\sum_{j=1}^{n-1} \mu_j \partial F_j = \nu \sum_{i=1}^p N_i=\nu \sum_{j=1}^n
\partial F_j.$$
Since $\partial F_1,\ldots,\partial F_n$ are independent, the
comparison of the coefficients of the $\partial F_n$'s gives $\nu=0$ and then
$\mu_j=0$.  Finally, the cycle $\sum \lambda_i \ell_i + \sum \mu_j \partial F_j$
is equal to zero and $F(\lambda_i,\mu_j)=1$. Thus, we proved that
$F(\lambda_i,\mu_j)$ is equal to zero, except if all $\lambda_i$'s and
$\mu_j$'s are zero, in which case it is equal to $1$. This proves the result.
\qed

The last element to study in the fundamental system is $h_{\partial F_n}$. We
have
$$\partial F_n=\sum_{i=1}^{p} N_i -\sum_{j=1}^{n-1} \partial F_j,$$
so that $h_{\partial F_n}=x_1\ldots x_p h_{\partial F_1}^{-1}\ldots
h_{\partial F_{n-1}}^{-1}$ under $\nu_{x_1}\otimes\ldots\otimes
\nu_{x_p}\otimes dg'$.

\begin{proposition} \label{law_1} Set $x=x_1\ldots x_p$ if $M$ has a boundary
  and $x=1$ if $M$ is closed. For any function $f$ continuous on $G^{2g+n+p}$,
\begin{eqnarray*}
&& \hskip -.5cm
\int_{G^\Gamma}f(h_{\ell_1},\ldots,h_{\ell_{2g}},h_{N_1},\ldots,h_{N_{p-1}}, 
h_{\partial F_1},\ldots,h_{\partial F_n}) \; dP(x_1,\ldots,x_p) = \\
&& \hskip -.5cm \int_{G^{2g+n}} f(u_1,\ldots,u_{2g},x_1,\ldots,x_{p-1},v_1,\ldots,v_n)
p_{\sigma(F_1)}(v_1)\ldots p_{\sigma(F_n)}(v_n) \\
&&  \hskip 7cm du_1 \ldots du_{2g}\; d\nu_x^n(v_1,\ldots,v_n).
\end{eqnarray*}
\end{proposition}
Note that, in the Abelian setting, the measure $\nu_x^n$ is invariant by
permutations of the factors in $U(1)^n$. 

\subsection{Gaussian aspect of the Abelian theory}

We proved that the law of the whole family $(h_c)_{c\in C\Gamma}$ is
determined by the law of a fundamental system and we just described this law.
So we could consider that the proposition \ref{law_1} is the answer to our
question. In fact, it is possible to be much more explicit by taking the
gaussian character of the Abelian theory into account. The crucial part of the
law of a fundamental system is of course that of $(h_{\partial
  F_1},\ldots,h_{\partial F_n})$. We will concentrate on this law.

\begin{proposition} \label{gaussian_representation} Let $Y_1,\ldots,Y_n$ be
  independent centered real gaussian random variables with $Y_i\sim
  \n(0,\sigma(F_i))$.  Let $S=Y_1+\ldots +Y_n$ be their sum. For each
  $i=1,\ldots,n$, set
  $$X_i=Y_i-{{\sigma(F_i)}\over {\sigma(M)}} S.$$
  Let $T$ be a real random
  variable, independent of the $Y_i$'s, with the following discrete law :
  $$P(T=t)=\cases{\left(\displaystyle \sum_{s,e^{is}=x}
      e^{-{{s^2}\over{2\sigma(M)}}}\right)^{-1} e^{-{{t^2}\over {2\sigma(M)}}}
    \;\;{\rm if}\;\; e^{it}=x \cr 0\;\;{\rm otherwise},}$$
  where, as before,
  $x=x_1\ldots x_p$ if $M$ has a boundary and $x=1$ if $M$ is closed. Then,
  for any function $f$ continuous on $G^n$,
\begin{equation}
\label{eqn1}
\int_{G^\Gamma} f(h_{\partial F_1},\ldots,h_{\partial F_n})\; dP_\beta = E\;
f\left(e^{i\left(X_1+{{\sigma(F_1)}\over{\sigma(M)}}
    T\right)},\ldots,e^{i\left(X_n+{{\sigma(F_n)}\over{\sigma(M)}} T\right)}
\right) 
\end{equation}
\end{proposition}

The law of $T$ described in this theorem is just that of a
$\n(0,\sigma(M))$ random variable conditioned to take its values in
$\exp^{-1}(x)$, where $\exp(t)=e^{it}$. We shall discuss the meaning of this
variable in section \ref{meaning_of_T}. \\

\pf. In this proof, we set $\sigma_i=\sigma(F_i)$ and
$\sigma_M=\sigma(M)$. One easily computes 
$$E X_i X_j =\delta_{ij}\sigma_i -{{\sigma_i \sigma_j}\over {\sigma_M}}$$
and
$\sum X_i = 0$ a.s. The law of $(X_1,\ldots,X_n)$ has no density with respect
to Lebesgue measure on $\R^n$, but that of $(X_1,\ldots,X_{n-1})$ does, on
$\R^{n-1}$. Denote by $C$ the $(n-1)\times (n-1)$ covariance matrix of
$(X_1,\ldots,X_{n-1})$. One easily checks that $C^{-1}$ is given by
$$(C^{-1})_{ij}={{\delta_{ij}\over {\sigma_i}}}+{1\over {\sigma_n}}.$$
So the density of the law of $(X_1,\ldots,X_n)$ is :
$$d\eta(t_1,\ldots,t_{n-1})= {1\over Z} \exp  -{1\over 2}
\left(\sum_{i=1}^{n-1} {{t_i^2}\over \sigma_i} + \sum_{i,j=1}^{n-1}
  {{t_it_j}\over {\sigma_n}}  \right)  \; dt_1 \ldots dt_n.$$
Let us fix a number $t_0$ such that $e^{it_0}=x$. We also set
$t_n=-t_1-\ldots-t_{n-1}$. We can compute the right term of (\ref{eqn1}): it
is equal to
\begin{eqnarray*}
&&{1\over Z} \sum_{k \in \sZ} \int_{\sR^{n-1}}
f\left(e^{it_1+{{\sigma_1}\over{\sigma_M}}(t_0+2k\pi)},\ldots,e^{it_{n-1}+{{\sigma_{n-1}}\over{\sigma_M}}(t_0+2k\pi)},e^{it_n+{{\sigma_n}
      \over{\sigma_M}}(t_0+2k\pi)}\right) \\
&& \hskip 2cm  \exp {-{1\over
    2}\left(\sum_{i=1}^{n-1}{{t_i}^2\over{\sigma_i}}  + 
    \sum_{i,j=1}^{n-1} {{t_it_j}\over {\sigma_n}}\right)}
\exp {-{(t_0+2k\pi)^2}\over{2\sigma_M}}\; dt_1 \ldots dt_{n-1}\\
&=&{1\over Z}\sum_{k,q_1,\ldots,q_{n-1} \in \sZ}\int_{[0,2\pi]^{n-1}}
f(e^{it_1},\ldots,e^{it_{n-1}},e^{i(t_n+t_0)}) \\
&& \hskip 0cm \exp \left(
-{1\over 2} \sum_{i=1}^{n-1} {1\over {\sigma_i}}(t_i+2\pi
q_i-{{\sigma_i}\over{\sigma_M}}(t_0+2k\pi))^2 + \right. \\
&& \hskip 0cm\left. + {1\over {2\sigma_n}} \left(
  \sum_{i=1}^{n-1} t_i+2\pi q_i-{{\sigma_i}\over
    {\sigma_M}}(t_0+2k\pi)\right)^2 \right)
\exp {-{{(t_0+2k\pi)^2}\over{2\sigma_M}}} \; dt_1 \ldots dt_{n-1}.
\end{eqnarray*}

We do not care about normalization constants, since two probability measures
with proportional densities are equal. Now we compute the left hand side of
(\ref{eqn1}) using the following expression of the heat kernel :
$$p_s(e^{it})={1\over Z}\sum_{p\in\sZ} e^{-{{(t-2p\pi)^2}\over {2s}}},$$
which
is just the image by the exponential map of the heat kernel on $\R$. We get
\begin{eqnarray*}
&& \hskip -1.25cm \int_{G^\Gamma} f(h_{\partial F_1},\ldots,h_{\partial F_n})\; dP_\beta =
\int_{G^n} f(v_1,\ldots,v_n) p_{\sigma_1}(v_1)\ldots p_{\sigma_n}(v_n)
\nu_x^n(v_1,\ldots,v_n)\\
&=&\sum_{p_1,\ldots, p_n\in \sZ} \int_{[0,2\pi]^{n-1}}
f(e^{it_1},\ldots,e^{it_{n-1}},e^{i(t_n+t_0)}) \exp \left(-{1\over 2}
\sum_{i=1}^{n-1} {{(t_i-2p_i\pi)^2}\over{\sigma_i}}\right) \\
&& \hskip 4.5cm  \exp \left(-{1\over 2}
{{(t_n+t_0-2p_n\pi)^2}\over {\sigma_n}}\right) \; dt_1 \ldots dt_{n-1}.
\end{eqnarray*}

The result will be a consequence of the following equality:
\begin{eqnarray*}
&&\hskip -1.25cm \sum_{p_1,\ldots,p_n \in\sZ} \exp \left( -{1\over 2}
\sum_{i=1}^{n-1} {{(t_i-2p_i\pi)^2}\over{\sigma_i}} -{1\over 2}
{{(t_n+t_0-2p_n\pi)^2}\over {\sigma_n}}\right) = \\
&&\hskip -1.25cm\sum_{k,q_1,\ldots,q_{n-1}\in \sZ} \exp
  {-{(t_0+2k\pi)^2}\over{2\sigma_M}} -{1\over 2} \sum_{i=1}^{n-1} {1\over
    {\sigma_i}}(t_i+2\pi q_i-{{\sigma_i}\over{\sigma_M}}(t_0+2k\pi))^2
  -\\
&& \hskip 5cm -{1\over {2\sigma_n}} \left( \sum_{i=1}^{n-1} t_i+2\pi
    q_i-{{\sigma_i}\over {\sigma_M}}(t_0+2k\pi)\right)^2
\end{eqnarray*}

Setting $q_n=-q_1-\ldots-q_{n-1}$, we have
\begin{eqnarray*}
&&\hskip -.5cm \sum_{i=1}^{n-1} {1\over {\sigma_i}}(t_i+2\pi
  q_i-{{\sigma_i}\over 
  {\sigma_M}}(t_0+2k\pi))^2 +{1\over {\sigma_n}}(\sum_{i=1}^{n-1}t_i+2\pi q_i
  -{{\sigma_i}\over {\sigma_M}}(t_0+2k\pi))^2 = \\
&=& \sum_{i=1}^{n-1} {1\over {\sigma_i}} ([t_i+2\pi q_i]-{{\sigma_i}\over
  {\sigma_M}}(t_0+2k\pi))^2+ \\
&& \hskip 3cm + {1\over {\sigma_n}}\left ([-t_n-2\pi
  q_n-(t_0+2k\pi)] + {{\sigma_n}\over {\sigma_M}}(t_0+2k\pi)\right)^2 \\
&=& \sum_{i=1}^{n-1} {1\over {\sigma_i}} \left[(t_i+2\pi q_i)^2 -{{2\sigma_i}\over
  {\sigma_M}}(t_i+2\pi q_i)(t_0+2k\pi)+{{\sigma_i^2}\over
  {\sigma_M^2}}(t_0+2k\pi)^2\right]+ \\
&&\hskip 2.5cm+ {1\over {\sigma_n}}\left[(-t_n-2\pi q_n -(t_0+2k\pi))^2 +{{\sigma_n^2}\over
  {\sigma_M^2}}(t_0+2k\pi)^2+ \right. \\
&&\hskip 4.5cm \left. +{{2\sigma_n}\over{\sigma_M}}(-t_n-2\pi q_n
  -(t_0+2k\pi))(t_0+2k\pi)\right]\\
&=& \sum_{i=1}^{n-1} {1\over {\sigma_i}}(t_i+2\pi q_i)^2+ {1\over
  {\sigma_M}}(t_0+2k\pi)^2\left(\sum_{i=1}^{n-1} {{\sigma_i}\over {\sigma_M}}
  + {{\sigma_n}\over {\sigma_M}}\right) \\
&&\hskip 1cm - {2\over {\sigma_M}}\sum_{i=1}^{n-1}(t_i+2\pi q_i)(t_0+2k\pi)
  +{2\over
  {\sigma_M}}\sum_{i=1}^{n-1}(t_i+2\pi q_i)(t_0+2k\pi)\\
&&\hskip 3.5cm - {2\over {\sigma_M}} (t_0+2k\pi)^2 +{1\over
  {\sigma_n}}(-t_n-2\pi q_n 
  -(t_0+2k\pi))^2 \\
&=& \sum_{i=1}^{n-1} {1\over {\sigma_i}} (t_i+2\pi q_i)^2 + {1\over
  {\sigma_n}} (t_n+t_0 +2\pi q_n+2k\pi)^2-{1\over {\sigma_M}}(t_0+2k\pi)^2.
\end{eqnarray*}
Setting $p_1=q_1,\ldots,p_{n-1}=q_{n-1}$ and $p_n=-q_n-k$, we get the
result. \qed

\subsection{The double layer potential}

To go further, we would like to represent isometrically the vector
$(X_1,\ldots,X_n)$ by a vector of functions of $L^2(M,\sigma)$ naturally
associated with $F_1,\ldots,F_n$. To begin with, remark that the vector
$(\1_{F_1},\ldots,\1_{F_n})$ has the same covariance as $(Y_1,\ldots,Y_n)$.
Now set
$$u_i=\1_{F_i}-{{\sigma(F_i)}\over {\sigma(M)}}.$$
The vector
$(u_1,\ldots,u_n)$ has the same covariance as $(X_1,\ldots,X_n)$. Each $u_i$
can be seen as the orthogonal projection of $\1_{F_i}$ on the hyperplane
$L^2_0(M,\sigma)$ of functions whose mean is equal to zero. In fact, the
$u_i$'s are the most natural generalizations on $M$ of the classical index of
a loop around a point in the plane. We will give a more direct definition of
$u_i$.

To do this, we endow $M$ with a Riemannian metric whose volume coincides with
$\sigma$. There exist a lot of such metrics, because $\sigma$ is equivalent to
the Lebesgue measure in any chart with a smooth density, as well as the
Riemannian volume of any Riemannian metric on $M$.

The choice of a compatible metric on $M$ gives rise to a Laplace operator
$\Delta$ and to a Hodge operator $*$ on $\bigwedge^1(T^*M)$. There exists on
$M$ a Green function $G : M\times M \lra \R_+$ defined outside the diagonal
which is symmetric, smooth and such that
\begin{equation}
\label{def_G}
\cases{\Delta G(x,\cdot) = \delta_x - {1\over {\sigma(M)}} \;\; \forall x\in
  M \cr
\int_M G(x,y) \; d\sigma(y)=0 \;\; \forall x\in M \;\; {\rm when} \;\;
  \partial M = \emptyset \cr
*dG_x=0 \;\; {\rm on} \;\; \partial M \;\; \forall x\in M \;\; {\rm when}\;\;
\partial M \neq \emptyset,}
\end{equation}
where $G_x$ denotes the function $G(x,\cdot)$. A proof of this fact can be
found in \cite{Aubin}. Note that when $M$ has a boundary, there exists a
solution to $\Delta \widetilde G_x = \delta_x$.  Nevertheless, this choice
would be incompatible with the condition $*dG_x=0$ on $\partial M$ which
implies $\int_M \Delta G_x =0$.

\begin{definition} \label{def_dlp}
Let $c$ be a path on $M$. We call double layer potential of
  $c$ the function $u_c$ defined on $M$ outside the image of $c$ by :
$$u_c(x)=\int_c *dG_x.$$
\end{definition}

Note that the double layer potential is additive: if $c_1$ and $c_2$ are two
cycles of $C\Gamma$, then $u_{c_1+c_2}=u_{c_1}+u_{c_2}$ $\sigma$-a.e. on $M$. 

\begin{proposition} \label{compute_dlp} Let $c$ be a simple loop which
is the boundary of a subset $U$ of $M$. Set $V=U^c$. Then
$$u_c(x)={{\sigma(V)}\over {\sigma(M)}} \1_U(x)-{{\sigma(U)}\over
  {\sigma(M)}}\1_V(x)=\1_U - {{\sigma(U)}\over {\sigma(M)}}.$$
In particular, $u_c \in L^2(M,\sigma)$ and $\parallel u_c \parallel_2
=\left( {{\sigma(U)\sigma(V)}\over {\sigma(M)}}\right)^{{1\over 2}}$.
\end{proposition}

\pf. Let $x$ be in $U$. Since $*dG_x=0$ on $\partial M$, we have:
$$ u_c(x)= \int_{\partial U} *dG_x = -\int_{\partial V} *dG_x
= -\int_V \delta_x -{1\over {\sigma(M)}} = {{\sigma(V)}\over {\sigma(M)}}.$$
Now let $x$ be in $V$.
$$u_c(x)=\int_{\partial U} *dG_x = \int_U \delta_x - {1\over {\sigma(M)}}\\
= -{{\sigma(U)}\over {\sigma(V)}}.$$
The last part of the statement follows
easily.\qed

\begin{corollary}\label{id_dlp} The vector $(u_1,\ldots,u_n)$ is equal to
  $(u_{\partial F_1},\ldots,u_{\partial F_n})$.
\end{corollary}

To go from functions $u_{\partial F_i}$ to random variables $X_i$, we need an
isometry of $L^2(M,\sigma)$ into a gaussian space, in other words a white
noise on $(M,\sigma)$. Let us consider a white noise
\begin{eqnarray*}
W:L^2(M,\sigma)&\lra& \g \\
u &\longmapsto & W(u)
\end{eqnarray*}
such that for any $u,v\in L^2(M)$, $W(u)$ and $W(v)$ are real centered
gaussian random variables such that $E[W(u)W(v)]=(u,v)_{L^2}$. The proposition
\ref{gaussian_representation} can be rewritten in the following form:

\begin{proposition} \label{wn_representation} The following equality holds in
  law:  
  $$(h_{\partial F_1},\ldots,h_{\partial F_n}) \build{=}_{}^{\law} \left(e^{i
      \left(W(u_{\partial F_1}) + {{\sigma(F_1)}\over
          {\sigma(M)}}T\right)},\ldots,e^{i \left(W(u_{\partial F_n}) +
        {{\sigma(F_n)}\over {\sigma(M)}}T\right)} \right).$$
\end{proposition}
 
We would like to extend this result to arbitrary cycles homologous to
zero. Let $c_1,\ldots c_k$ be cycles of $C_0\Gamma$. For each $i$, $c_i$ is a
linear combination of the $\partial F_i$ thus $u_{c_i}$ is well defined and is
in $L^2(M)$. So, $W(u_{c_i})$ is well defined. 

We have to generalize the term ${{\sigma(F_i)}\over {\sigma(M)}}$. Since $c_i$
is homologous to zero, it is the boundary of a two-chain denoted by
$\alpha$. If $M$ has a boundary, $H_2(M)=0$ and $\alpha$ is well defined by
$\partial \alpha=c$. We identify $\sigma$ with a $2$-form on $M$ and set
$\sigma(\alpha)=|\langle \sigma,\alpha \rangle |$, using the natural pairing between
$2$-forms and $2$-chains. So the number
$$\sigmaint(c)={{\sigma(\alpha)}\over {\sigma(M)}}$$
is well defined. If $M$
is closed, $H_2(M)\simeq \Z$ and $\alpha$ is defined up to a multiple of
$[M]$. So the number $\sigmaint(c)$ is only defined modulo $1$. But in
this case, $T$ takes its values in $\exp^{-1}(1)=2\pi \Z$ so that $\exp i
\sigmaint(c)T$ is well defined.

\begin{proposition} \label{wn_rep_loops}
  Let $(c_1,\ldots,c_k)$ be cycles of $C_0\Gamma$. Then the following equality
  in law holds:
  $$(h_{c_1},\ldots,h_{c_k}) \build{=}_{}^{\law} \left(e^{i \left(W(u_{c_1}) +
        \sigmaint(c_1) T\right)},\ldots,e^{i \left(W(u_{c_k}) + \sigmaint(c_k)
        T\right)} \right).$$
\end{proposition}

\pf. By proposition \ref{wn_representation}, the result is true when
$(c_1,\ldots,c_k)=(\partial F_1,\ldots, \partial F_n)$. Since the boundaries
of the faces constitute a basis of $C_0\Gamma$, it is sufficient to show that
the new set of variables defined using the white noise satisfy the same
multiplicativity property as $(h_c)_{c\in C_0\Gamma}$.

On one hand, $W$ is linear and the double layer potential is additive, so that
$\exp i W(u_{c_1+ c_2})= \exp i W(u_{c_1}) \exp i W(u_{c_2})$. On the other
hand, $c_1=\partial \alpha_1$ and $c_2=\partial \alpha_2$ imply
$c_1+c_2=\partial (\alpha_1 + \alpha_2)$, so $\sigmaint$ is also additive.
This proves the result. \qed

The results that we proved in this section are the starting point of the more
detailed investigation that will be done in chapter \ref{abelian_theory},
after the continuous Yang-Mills measure has been constructed. 

Some properties of the double layer potential will be proved in the next
chapter, using a favorable technical context. Nevertheless, it is necessary to
state here a fundamental property that will be proved at the end of the
chapter \ref{abelian_theory}.

\begin{theorem} \label{dlp_l2}
For any path $c$ of $PM$, the function $u_c$ is in
  $L^2(M,\sigma)$. 
\end{theorem}


\pagevide
\chapter{Continuous Yang-Mills measure}
\label{Continuous_theory}

In chapter 1, we defined a random holonomy along the paths of a graph on $M$.
Our aim in this chapter is to extend this definition to all paths on $M$. The
problem is that there are families of paths that cannot be realized as
subfamilies of any $\Gamma^*$, $\Gamma$ being a graph on $M$: even two smooth
paths can cross each other an infinity of times an give rise to an infinity of
connected components on $M$, a situation in which we are unable to write the
joint law of their holonomies using the tools of the preceding chapter.

The two properties of the discrete theory that are essential to our purpose
are the invariance by subdivision, expressed in theorem
\ref{invariance_by_subdivision} and the regularity property of the proposition
\ref{fundamental_estimate}. Our basic idea is to cover $M$ with finer and
finer graphs and to prove that the discrete measures on these graphs converge
in some sense to a continuous object that will be called continuous Yang-Mills
measure.

\section{Projective systems}

Let $(\Gamma_\lambda)_{\lambda \in \Lambda}$ be a family of graphs on $M$,
which approximates correctly $M$, whatever this means exactly. In sections
\ref{discrete_YM} and \ref{conditionning}, we explained how to construct a
family of probability spaces corresponding to this family of graphs. We will
now consider the projective limit of this family of probability spaces.

\begin{definition} Let $\Lambda$ be an ordered set such that for all
  $\lambda,\mu \in \Lambda$, there exists $\nu\in \Lambda$ such that $\lambda
  <\nu$ and $\mu < \nu$. A projective family of probability spaces indexed by
  $\Lambda$ is a family $(\Omega_\lambda,P_\lambda)$ of probability spaces
  together with a family of measurable maps $f_{\lambda \mu}: \Omega_\mu \lra
  \Omega_\lambda$ defined for all $\lambda <\mu$, such that \\
  \indent  1. $f_{\lambda \lambda}=\id_{\Omega_\lambda}$, \\
  \indent 2. $f_{\lambda \nu} = f_{\lambda \mu} \circ f_{\mu \nu} $ if
  $\lambda <\mu
  <\nu $, \\
  \indent  3. $f_{\lambda \mu} (\Omega_\mu) = \Omega_\lambda$, \\
  \indent 4. $(f_{\lambda \mu})_* P_\mu = P_\lambda$.
\end{definition}

The projective limit of such a system is by definition the set 
$$\build{\lim}_{\longleftarrow}^{} \Omega_\lambda=\Omega=\{
(\omega_\lambda)_{\lambda\in\Lambda} \in \prod_{\lambda\in\Lambda}
\Omega_\lambda \; |\; \forall \lambda<\mu, f_{\lambda
  \mu}(\omega_\mu)=\omega_\lambda \} .$$

For each $\lambda\in \Lambda$, the projection on the $\lambda$-th coordinate
defines a map $f_{\lambda}: \Omega \lra \Omega_{\lambda}$. The main result is
that, under certain conditions, there exist a $\sigma$-algebra and a
probability measure on $\Omega$ which are consistent with all $P_\lambda$ via
the maps $f_{\lambda}$.

\begin{theorem} [\cite{Choksi} 2.2] \label{choksi} If all $\Omega_\lambda$ are
  compact Borel probability spaces, then there exist a $\sigma$-algebra and a
  probability measure $P$ on $\Omega$ such that, for all $\lambda$ in
  $\Lambda$, ${f_{\lambda}}_*P=P_\lambda$. The space $(\Omega,P)$ is called
  projective limit of the family $(\Omega_\lambda,P_\lambda;f_{\lambda\mu})$.
\end{theorem}

It would be appealing to take the family of all graphs as index set $\Lambda$,
with the order defined in \ref{finess}. The problem is that given two graphs,
it is not always true that there exists a third graph which is finer than both
others, just because two edges belonging to two different graphs can intersect
very badly. So, the first assumption about the ordering on $\Lambda$ would not
be satisfied.

\section{Piecewise geodesic graphs}

We are led to consider a family of graphs small enough for that problem not to
occur. A convenient family is that of graphs with edges piecewise geodesic for
some Riemannian metric on $M$. Another possibility is to consider graphs with
piecewise real analytic edges for some complex structure on $M$. This has been
investigated by Ashtekar and Lewandowski \cite{Ashtekar_Lewandowski}.

We fix a surface $(M,\sigma)$. For technical reasons, we suppose that $M$ is
closed, until the section \ref{boundary_case} where we shall derive the
general case from the case without boundary.

Let us choose $q$ disjoint simple loops $L_1,\ldots,L_q$ on $M$ whose image is
a smooth submanifold of $M$, in other words, $q$ disjoint embeddings of $S^1$
into $M$. We will sometimes think of these loops as the boundary of a
submanifold of $M$ or just as loops along which we want to contition the
holonomy.

\begin{proposition} \label{metric_Moser}
  There exists a Riemannian metric on $M$ whose Riemannian volume coincides
  with $\sigma$ and such that $L_1,\ldots,L_q$ are geodesics.
\end{proposition}

\pf. If $q$=0, let us choose an arbitrary metric on $M$. Its Riemannian volume
is equivalent to $\sigma$, with a smooth density. Multiplying the metric by an
appropriate smooth positive function, we get a new metric, conformal to the
first one, whose Riemannian volume is exactly $\sigma$.

If $q>0$, the proof is more complicated. Let us first construct a metric for
which all $L_i$'s are geodesic. Each $L_i$ has a tubular neighbourhood in $M$
which is diffeomorphic to a cylinder $S^1 \times (-1,1) \ni (\theta,t)$, with
$L_i=\{t=0\}$. In these coordinates, $L_i$ is certainly geodesic for the
metric $d\theta^2 + dt^2$. If the tubular neighbourhoods were chosen small
enough to be disjoint, this procedure defines a Riemannian metric on the
reunion of these cylinders and we extend it arbitrarily to a metric $g_0$ on
the whole surface $M$.  The loops $L_1,\ldots, L_q$ cut $M$ into submanifolds
with boundary $M_1,\ldots,M_k$. Multiplying $g_0$ by a good positive
function which is identically equal to $1$ in a neighbourhood of each
$L_i$, we can obtain a new metric $g_1$ for which all $L_i$ are still
geodesics and also such that
\begin{equation}
\label{rel-vol}
\vol {}_{g_1}(M_i)=\sigma(M_i) \;\;\; \forall i=1\ldots k.
\end{equation}

Now we must redistribute the surface inside each $M_i$. We adapt a proof of
Moser's theorem, which is for example proved in \cite{Berger_Gostiaux}.

\begin{theorem}[Moser] Let $\alpha$ and $\beta$ be two volume $2$-forms on a
  closed compact surface $M$ such that
$$\int_M \alpha=\int_M \beta.$$
There exists a diffeomorphism $\phi: M \lra M$ such that $\phi^*\beta=
\alpha$.
\end{theorem} 

We choose two forms $\alpha$ and $\beta$ representing respectively $\sigma$
and $\vol_{g_1}$. We know more about them than what is needed in Moser's
theorem, but we also want to prove more: we would like to find a
diffeomorphism of $M$ that sends $\alpha$ to $\beta$ and also that preserves
the $L_i$'s, so that they remain geodesic after pulling back $g_1$ by the
diffeomorphism.

For each $i=1\ldots q$, let $j_i: L_i([0,1]) \lra M$ denote the canonical
injection. The fact that $\int_M \beta-\alpha=0$ implies that there exists a
form $\gamma \in \Lambda^1(T^*M)$ such that $d\gamma=\beta-\alpha$. We show
that $\gamma$ can be chosen such that $j_i^*\gamma=0$ for all $i$, or in other
words such that $\gamma(X)=0$ for any vector $X$ tangent to a $L_i$.  Pick
$\gamma$ such that $d\gamma =\beta - \alpha$, consider the element
$(\int_{L_1} \gamma,\ldots,\int_{L_q} \gamma)$ of $\R^q$ and the $q$-uple
$([L_1],\ldots,[L_q])$ of vectors of $H_1(M;\Z)$.  Suppose that
$n_1,\ldots,n_q$ are integers such that $\sum_i n_i [L_i]=0$.  Then $\sum_i
n_i [L_i]$ is the boundary of a $2$-chain in $M$, which is necessarily a
linear combination of $M_1,\ldots,M_k$. Thanks to (\ref{rel-vol}), this
implies $\int_{\sum n_i L_i} \gamma = 0 =\sum_{i=1}^q n_i \int_{L_i} \gamma.$
We proved that a relation $\sum n_i [L_i] =0$ implies $\sum n_i \int_{L_i}
\gamma =0$. Thus there exists a linear form $\zeta$ on $H_1(M;\Z)$ such that
$(\zeta,[L_i])=\int_{L_i}\gamma$. This linear form can be represented by an
element of $H^1(M;\R)$ and this element can be represented by a closed
$1$-form on $M$ that we still denote by $\zeta$. The form $\gamma-\zeta$
satisfies:
$$
\cases{ d(\gamma-\zeta)=d\gamma=\beta-\alpha \cr \int_{L_i}\gamma-\zeta=0
  \;\;\; \forall i=1\ldots q.}$$
This last relation proves that for each $i$,
$j_i^*(\gamma-\zeta)$ is exact on $L_i$ and can be written $du_i$ with $u_i
\in C^\infty(L_i)$. Let $u$ be a smooth function on $M$ such that
$u_{|L_i}=u_i$ for each $i$. The form $\gamma-\zeta-du$ satisfies:
$$\cases{ d(\gamma-\zeta-du)=\beta-\alpha \cr 
 j_i^*(\gamma-\zeta-du)=0 \;\;\; \forall i=1\ldots q.}$$

We proved that it is possible to choose $\gamma$ such that $\gamma(X)=0$ for
each vector $X$ tangent to a $L_i$ and we choose $\gamma$ in that way. The end
of the proof is similar to that of Moser's theorem.

For each $t\in[0,1]$, set $\alpha_t=(1-t)\alpha + t\beta$ and define the
vector field $X_t$ on $M$ by $i_{X_t}\alpha_t=-\gamma$.  The field $X_t$
depends smoothly on $t$ and induces a flow $(\phi_t)_{t\in[0,1]}$. We compute
the derivative of $\phi_t^*\alpha_t$. For any $t_0 \in (0,1)$,
$${d\over {dt}}\biggl|_{t=t_0} (\phi_t^*\alpha_t) = {d\over
  {dt}}\biggl|_{t=t_0} (\phi_t^*\alpha_{t_0}) +
{d\over{dt}}\biggl|_{t=t_0}(\phi_{t_0}^*\alpha_t).$$
The second term of the
r.h.s. is equal to
$$\phi_{t_0}^*\left({d\over{dt}}\biggl|_{t=t_0}
  \alpha_t\right)=\phi_{t_0}^*(\beta-\alpha).$$
We denote by $\l_{X_t}$ the Lie derivative with respect
to the field $X_t$ and use Cartan's relation $\l_{X_t}=d\circ i_{X_t} +
i_{X_t}\circ d$. We find that the first term is equal to
$$\l_{X_{t_0}} \phi_{t_0}^*\alpha_{t_0} = d\left(
  i_{X_{t_0}}\phi_{t_0}^*\alpha_{t_0}\right) = d(\phi_{t_0}^*
i_{X_{t_0}}\alpha_{t_0})=
-d(\phi_{t_0}^*\gamma)=-\phi_{t_0}^*(d\gamma)=-\phi_{t_0}^* (\beta-\alpha).$$
Thus, ${d\over {dt}} \phi_t^*\alpha_t=0$, so that
$\phi_1^*\beta=\phi_1^*\alpha_1=\phi_0^*\alpha_0=\alpha.$ For any vector $X$
tangent to a $L_i$, the equality
$\alpha_t(X_t,X)=i_{X_t}\alpha_t(X)=\gamma(X)=0$ proves that the field $X_t$
is tangent to $L_i$. So the flow $\phi_t$ preserves the $L_i$'s.

The Riemannian volume of the metric $g=\phi_1^*g_1$ is
$\phi^*\vol_{g_1}=\sigma$. Moreover, $\phi_1$ is an isometry from $(M,g)$ into
$(M,g_1)$ that preserves the $L_i$'s. Since they are geodesics for $g_1$, they
are also geodesics for $g$. \qed

From now on, we fix on $M$ a metric given by the last proposition. Let us
recall a classical result that summarizes the main property of the geodesics
that we will use. A proof of a local version of this theorem can be found for
example in \cite{do_Carmo} (proposition 3.4.2). The compactness of $M$ allows
to globalize the result. 

\begin{theorem} \label{injectivity_radius} There exists a positive
  real number $R_M$, called convexity radius of $M$, such that if $x$ and $y$
  are two points of $M$ contained in a ball of radius smaller than $R_M$, they
  are joined by a unique piece of minimizing geodesic and this piece of
  geodesic stays inside the ball.
\end{theorem}

This theorem implies in particular the following result:

\begin{proposition} \label{geod-intersection}
  Let $\zeta_1$ and $\zeta_2$ be two finite pieces of geodesics. The
  intersection of $\zeta_1$ and $\zeta_2$ is the union of a finite number of
  isolated points and at most two segments.
\end{proposition}

\pf. If $\zeta_1$ and $\zeta_2$ meet at an infinity of points, it is easy to
check that there exists a couple $(t_1,t_2)$ of times such that
$\zeta_1(t_1)=\zeta_2(t_2)$ and $\dot\zeta_1(t_1)=\dot\zeta_2(t_2)$. So they
are two pieces of the same infinite geodesic. If this geodesic is periodic,
$\zeta_1$ and $\zeta_2$ can intersect along one or two segments. Otherwise,
they have one segment in common plus a finite number of isolated points. \qed

We denote by $\g$ the set of graphs whose edges are piecewise geodesic and
such that $L_1,\ldots,L_q$ are in $\Gamma^*$. The set $\g$ is ordered by the
relation $<$.

\begin{proposition} \label{finer}
  Given two graphs $\Gamma_1,\Gamma_2$ in $\g$, there exists $\Gamma_3$ in
  $\g$ such that $\Gamma_1<\Gamma_3$ and $\Gamma_2<\Gamma_3$.
\end{proposition}

\pf. The idea is to superpose $\Gamma_1$ and $\Gamma_2$. Given an edge $a$ of
$\Gamma_2$, we know that $a([0,1])\cap \supp(\Gamma_1)$ is a finite reunion of
segments and points. So, it is possible to add a finite number of new vertices
and new edges to $\Gamma_1$ in such a fashion that $a$ becomes a path in the
new graph. Repeating this procedure for each edge of $\Gamma_2$ gives the
result. \qed

Let us fix an element $(x_1,\ldots,x_q)$ of $G^q$. With each graph $\Gamma\in
\g$, we associated a space $(G^\Gamma,P^\Gamma(x_1,\ldots,x_q))$. The last
proposition states the last property that was missing for the family
$(G^\Gamma,P^\Gamma(x_1,\ldots,x_q)),\Gamma\in\g$ to be a projective family
of probability spaces. Each $G^\Gamma$ is compact, so theorem \ref{choksi}
asserts that this projetive family has a projective limit
$(\Omega,\a,P(x_1,\ldots,x_q))$ which is a probability space endowed with
functions $f_\Gamma: \Omega \lra G^\Gamma$ such that ${f_\Gamma}_*
P(x_1,\ldots,x_q) =P^\Gamma(x_1,\ldots,x_q)$.  This space contains in itself
the same information that is contained in all spaces
$(G^\Gamma,P^\Gamma(x_1,\ldots,x_q))$: each random variable $h_\zeta: G^\Gamma
\lra G$ gives rise to a random variable
$$H_\zeta: \Omega \build{\lra}_{}^{f_\Gamma} G^\Gamma
\build{\lra}_{}^{h_\zeta} G$$
and the law of a $n$-uple
$(h_{\zeta_1},\ldots,h_{\zeta_n})$ computed in any graph of $\g$ is always
equal to that $(H_{\zeta_1},\ldots,H_{\zeta_n})$ under $P$. Moreover, since
any piecewise geodesic path can be seen as a path in a graph, there is a
well-defined random variable $H_\zeta$ on $\Omega$ associated with any such
path $\zeta$. Remark that the multiplicativity property is preserved:

\begin{proposition} Let $\zeta_1$ and $\zeta_2$ be two piecewise geodesic
  paths such that $\zeta_1(1)=\zeta_2(0)$. Then $H_{\zeta_1
  \zeta_2}=H_{\zeta_2}H_{\zeta_1}$ a.s.
\end{proposition}

From now on, we use greek letters to denote the piecewise geodesic paths and
denote by $PGM$ the set of these paths.

\section{Preliminary results}

\subsection{Lassos}

\begin{definition} A lasso is a simple loop or a path of the form
  $l=cbc^{-1}$, where $c$ is an injective path and $b$ a simple loop which
  meets $c$ only at its base point. The loop $b$ is detetermined by $l$ and is
  called the buckle of the lasso $l$.
\end{definition}

A notion of lasso close to this one has already been used by Driver in
\cite{Driver_lassos}. In \cite{Gross_King_Sengupta}, Gross, King and Sengupta
also suggested that the use of lassos might be helpful in this construction.

Lassos are useful at least for two reasons: the first one is that it is easy
to compute the law of their holonomy and the second one is that any reasonable
loop can be decomposed in some sense into a product of lassos.  Let us begin
with this second point.

There is a natural equivalence relation between paths, which is the
following:
\begin{definition} Two paths are said to be basically equivalent if one of
  them can be written $c_1 c_2 c_2^{-1} c_3$ and the other one $c_1 c_3$,
  where $c_1,c_2,c_3 \in PM$. Two paths $c$ and $c'$ are equivalent, and we
  denote $c\simeq c'$, if there exists a finite chain $c=c_0, \ldots,c_n=c'$
  such that any two successive terms of this chain are basically equivalent.
\end{definition}

\begin{lemma} Let $\zeta_1$ and $\zeta_2$ be two paths of $PGM$ and
suppose that $\zeta_1 \simeq \zeta_2$. Then $H_{\zeta_1}=H_{\zeta_2}$
$P(x_1,\ldots,x_q)$-a.s.
\end{lemma}

\pf. This is a consequence of the multiplicativity of the random
holonomy. \qed

Let us define the class of paths that can be decomposed into a product
of lassos.

\begin{definition} \label{finite-self-intersection}
A path $c\in PM$ is said to have finite
self-intersection if there exists a graph $\Gamma$ such that $c\in
\Gamma^*$. 
\end{definition}
Remark that this definition is not the usual one of finite
self-intersection. Indeed, our definition allows for a path a finite number of
points and also a finite number of segments as auto-intersection set. In
particular, the proposition \ref{geod-intersection} shows that a piecewise
geodesic path, which is a concatenation of injective pieces of different
geodesics, has finite self-intersection in the sense of
\ref{finite-self-intersection}. 

The Riemannian metric chosen on $M$ allows us to compute the length of
a path $c$, that we denote by $\ell(c)$.

\begin{proposition} \label{decomposition_in_lassos} Let $c$ be a path
  with finite self-intersection. \\
  \indent 1. If $c(0) \neq c(1)$ then $c$ is equivalent to a unique product $c
  \simeq l_1\ldots l_p c'$, where the $l_i$'s are lassos which are non
  equivalent to a constant loop and $c'$ is an injective path joining $c(0)$
  to $c(1)$.  Moreover, if $b_i$ denotes the buckle of the lasso $l_i$ for
  each $i$, the following inequality holds:
  $$\ell(c) \geq \sum_{i=1}^p \ell(b_i) + \ell(c').$$
  \indent 2. If
  $c(0)=c(1)$, the result remains true after removing $c'$ of all expressions.
\end{proposition}

\pf. We proceed by induction on the number of edges in a decomposition of $c$
as a path in a graph. If $c$ is an edge, we are in the first case and the
result is true. Suppose that $c=a_1 \ldots a_r$. If $c$ is an injective path
or a simple loop, the result is true. Otherwise, the idea is to trace $c$ out
until the first time it intersects itself. Let $i$ be the smaller integer such
that there exists $1 \leq j <i$ verifying $a_j(0)=a_i(1)$. Such an $i$ exists,
and $(i,j) \neq (r,1)$. We have:
$$c \simeq a_1 \ldots a_{j-1}\cdot a_j \ldots a_i \cdot (a_1 \ldots
a_{j-1})^{-1} \cdot a_1 \ldots a_{j-1} \cdot a_{i+1} \ldots a_r .$$
It can
happen that the first piece or the last piece of this decomposition are empty,
respectively if $j=1$ or $i=r$, but these two situations cannot coexist. If
$j=i-1$, it is possible that $a_i=a_j^{-1}$ so that $a_j a_i$ is equivalent to
a constant path. This cannot happen if $j< i-1$, in which case $a_j \ldots
a_i$ is a genuine simple path. Thus, the product of the three first terms is
either equivalent to a constant path, or is a simple loop (if $j=1$), or a
lasso. The product of the two last terms is the product of a number of edges
which is positive and strictly less than $p$. So by induction, this path
$\widetilde c$ is equivalent to $l_1 \ldots l_q c'$, or to $l_1 \ldots l_q$ if
$c$ is a loop.  Note that $\ell(c)=\ell(a_1)+\ldots +\ell(a_r)=\ell(\widetilde
c)+\ell(a_j)+\ldots +\ell(a_i)$. So, by induction, $\ell(c)\geq
\ell(b_1)+\ldots+\ell(b_q)+\ell(c')+\ell(a_j)+\ldots + \ell(a_i)$. In the case
$j=i-1$ and $a_j=a_i^{-1}$, we have $c \simeq \widetilde c$ and the result is
true with a strict inequality. Otherwise, there exists a lasso $l_0$ such that
$c \simeq l_0 \widetilde c$ and the length of the buckle of this lasso $l_0$
is exactly $\ell(a_j)+\ldots +\ell(a_i)$. \qed

\begin{figure}[h]
\begin{center}
\input{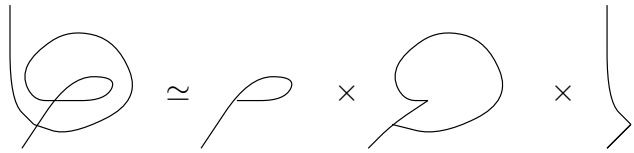}
\end{center}
\caption{An example of decomposition.
\label{dec}}
\end{figure}

\subsection{Holonomy along small piecewise geodesic loops}

In order to estimate the holonomy along a small lasso, we need, according to
the proposition \ref{fundamental_estimate}, to know the area enclosed by its
buckle. This area can be controlled by the length of the buckle using an
isoperimetric inequality. We recall a classical fact about open covering of
metric compact sets. A proof can be found in \cite{Massey}.

\begin{lemma} Let $M$ be a metric compact set. Let $(O_i)_{i\in I}$ be
an open covering of $M$. There exists a positive real number
$R$ called Lebesgue number of this covering, such that
for any ball $B$ of radius smaller than $R$ in $M$, there exists an
$i \in I$ such that $B\subset O_i$.
\end{lemma}

\begin{proposition} \label{isop} There exist $R>0$ and $K>0$
  such that any simple loop $l$ contained in a ball $B$ of radius smaller than
  $R$ is the boundary of an open set $U \subset B$ such that
  $$\sigma(U)\leq K \ell(l)^2.$$
\end{proposition}

\pf. Let $R_0$ be such that any closed geodesic ball of $M$ of radius smaller
than $R_0$ is diffeomorphic to a disk. Let $B_1,\ldots,B_n$ be a covering of
$M$ by open balls of radius $R_0$. Let us denote by $g$ the metric on $M$ and
$g_0$ the euclidean metric on $\R^2$. For each $i$, there exists a
diffeomorphism $\phi_i:\overline{B_i} \lra \overline{D}(0,1) \subset \R^2.$
Since the metrics $g$ and $\phi^*g_0$ can be compared on $\overline{B_i}$, the
usual isoperimetric inequality on $\overline{D}(0,1)$ gives rise to an
inequality on $\overline{B_i}$, with some constant $K_i$. Let $K$ be the
supremum of $K_1,\ldots,K_n$. Let $R$ be a Lebesgue number of the covering
$B_1,\ldots,B_n$. Then the statement holds with this choice of $K$ and
$R$. \qed 

Now we can estimate the holonomy along a small lasso:

\begin{proposition} \label{small_lassos} There exist $L_0>0$ and $K>0$ such
  that if $\lambda$ is a piecewise geodesic lasso whose buckle $\beta$ has a
  length smaller than $L_0$, then
  $$E\rho(H_\lambda) = E d(H_\lambda,1) \leq K \ell(\beta).$$
\end{proposition}

\pf. The lasso $\lambda$ can be written $\sigma \beta \sigma^{-1}$, so thanks
to invariance by conjugation of the distance on $G$, we have $E
\rho(H_\lambda) =E\rho(H_\beta)$. Let $L_0$ be shorter than the shortest
length of a loop non homotopic to a point and also shorter than the radius $R$
given py the proposition \ref{isop}. Then the hypothesis $\ell(\beta) \leq
L_0$ implies that $\beta$ is the boundary of a small disk $D$. Using
proposition \ref{fundamental_estimate}, we get
$$E \rho(H_\beta) \leq C \sqrt{\sigma(D)} \leq K \ell(\beta).$$
\vskip -1cm \qed

This result suggests that it will be possible to prove regularity results for
the random holonomy using the following distance between $G$-valued random
variables:

\begin{definition} Let $X$ and $Y$ be two $G$-valued random variables defined
  on the same probability space. The distance $d_P(X,Y)$ is defined by:
  $$d_P(X,Y)=E d(X,Y),$$
  where $d$ is the biinvariant Riemannian distance on $G$.
\end{definition}

The first example of such regularity results is the following one: 

\begin{proposition} \label{small_loops}
Let $\zeta$ be a piecewise geodesic loop of length smaller than $L_0$. Then
$$d_P(H_\zeta,1) \leq K \ell(\zeta).$$
\end{proposition}

\pf. Since $\zeta$ is piecewise geodesic, it has finite
self-intersection. So it is equivalent to a product of piecewise
geodesic lassos: $\zeta\simeq\lambda_1\ldots \lambda_p$. This gives:
\begin{eqnarray*}
d_P(H_\zeta,1) &=& d_P(H_{\lambda_1}\ldots H_{\lambda_p},1) \\
&\leq& d_P(H_{\lambda_1}\ldots H_{\lambda_p},H_{\lambda_2}\ldots
H_{\lambda_p})+\ldots +d_P(H_{\lambda_p},1)\\
&\leq& d_P(H_{\lambda_1},1)+\ldots +d_P(H_{\lambda_p},1) \\
&\leq& K (\ell(\beta_1) + \ldots \ell(\beta_p)) \\
&\leq& K \ell(\zeta).
\end{eqnarray*}
\vskip -1cm \qed

\subsection{Double layer potential of small piecewise geodesic loops}

Using the same techniques as in the preceding pararaph, we will
estimate the double layer potential of a small loop. This is the first step in
the proof of the proposition \ref{continuity_loops_dlp}, that will play an
important role in the study of the Abelian theory. 

Recall that the definition of the potential (see \ref{def_dlp}) depends on a
Riemannian metric on $M$ whose Riemannian volume is equal to $\sigma$. For the
moment, we only know that the potential of any element of $PM$ is in
$L^2(M,\sigma)$ (see theorem \ref{dlp_l2}).

\begin{proposition} Let $l$ be a lasso with buckle $b$. Suppose that
  $\ell(b)\leq L_0$, where $L_0$ is the length given by \ref{small_lassos}.
  Then
$$\parallel u_l \parallel_{L^2} \leq K \ell(b).$$
\end{proposition}
 
\pf. The length $L_0$ is such that $b$ is necessarily the boundary of a disk
$D$ whose area satisfies $\sigma(D) \leq K \ell(b)^2$. Thus, by proposition
\ref{compute_dlp},
$$\parallel u_b \parallel_{L^2} =\left( {{\sigma(D) \sigma(D^c)}\over
{\sigma(M)}} \right)^{{1\over 2}} \leq \sigma(D)^{{1\over 2}} \leq K
\ell(b).$$
Since $u_l=u_b$ a.e., we have the result. \qed

As in the preceding paragraph, this result can be extended to loops
with finite self-intersection. 

\begin{proposition} \label{small_loops_dlp}
  Let $l$ be a loop with finite self-intersection and of length smaller than
  $L_0$. Then
$$\parallel u_l \parallel_{L^2} \leq K \ell(l).$$
\end{proposition}

\pf. Let us write that $c$ is equivalent to a product of lassos: $c\simeq
l_1\ldots l_p$. Two paths that are equivalent have the same double layer
potential almost everywhere, so that
\begin{eqnarray*}
\parallel u_c \parallel_{L^2} &\leq& \parallel u_{l_1} + \ldots +
u_{l_p} \parallel_{L^2} \\
&\leq& \parallel u_{l_1} \parallel_{L^2} + \ldots + \parallel u_{l_p}
\parallel_{L^2} \\
&\leq& K ( \ell(b_1)+\ldots + \ell(b_p) ) \\
&\leq& K \ell(c).
\end{eqnarray*}
\vskip -1cm \qed

\label{correspondance}

The fact that the propositions \ref{small_loops} and \ref{small_loops_dlp} are
very similar will allow us later, in proposition \ref{dlp_reg} for example, to
transpose directly some regularity results about the random variables $H$ to
the double layer potential.

\subsection{Topology on the space of paths}

According to the proposition \ref{fundamental_estimate}, it seems to be
necessary to control the surface left between two loops in order to control
the distance between their holonomies. Given a Riemannian metric on $M$, the
uniform distance defined as follows:
$$d_\infty(c_1,c_2)=\inf \sup_{t\in [0,1]} d(c_1(t),c_2(t))$$
allows to control this surface, where the infimum is taken over all
reparametrizations of $c_1$ and $c_2$. 

In the paper \cite{Becker}, C. Becker says that the double layer potential of
a loop depends continuously of this loop in norm $L^2$ when the set of loops
on $M$ is endowed with the topology induced by the uniform norm. His
proposition $3.2$ depends on the validity of this assertion, which is probably
true if one restricts to simple loops, but not if one allows loops to have a
self-intersection, even a finite one. Let us describe a counterexample.
Becker stated his result on $\R^2$, but this does not change the situation
very much.  Let $M$ be the sphere $S^2$ embedded as usual in $\R^3$, endowed
with the standard metric. Let us consider the pencil of planes $(P_t)_{t\in
  [0,4\pi)}$ containing the horizontal line $z=0, y=-1$, indexed in the
following way: denoting by $C_t$ the intersection of $S^2$ with the lower
half-space bounded by $P_t$, we have $\sigma(C_t)=t$. For any $t\in [0,4\pi)$,
denote by $c_t$ a loop based at $(0,-1,0)$ whose image is the intersection of
$P_t$ with $S^2$, oriented negatively with respect to the $z$ axis.  Let
$0<t_1<\ldots<t_n<4\pi$ be $n$ distinct times and set $c=c_{t_1} \ldots
c_{t_n}$. For each $i$, $u_{c_i}=\1_{C_{t_i}}-{{t_i}\over {4\pi}}$. Thus,

\begin{eqnarray*}  
 \parallel u_c \parallel_{L^2}^2 &=&\parallel \1_{C_{t_1}} +\ldots + 
\1_{C_{t_n}} \parallel_{L^2}^2 + {{(t_1+\ldots + t_n)^2}\over {16\pi^2}} \\
&& \hskip 4cm -2\left({{t_1+\ldots + t_n}\over {4\pi}},\1_{C_{t_1}} + \ldots +
  \1_{C_{t_n}}\right)_{L^2} \\
&\geq& n^2 t_1 + {{n^2 t_1^2}\over {16\pi^2}} - {{n^2 t_n^2}\over
{2\pi}}.
\end{eqnarray*}

Suppose that $t_1={1 \over n}$ and $t_n={2\over n}$. Then $\parallel u_c
\parallel_{L^2}^2$ is of the order of $n$, so it grows to infinity when $n$
tends to infinity. But at the same time, the loop $c$ tends to the constant
loop equal to $(0,-1,0)$ in the topology induced by the distance $d_\infty$.
The potential of this constant loop being equal to zero, this contradicts the
continuity.

Therefore, it is necessary to endow the space of paths with a topology finer
than that induced by $d_\infty$ if we expect some kind of continuity of the
double layer potential and of the random holonomy. It has emerged in the last
paragraphs that the length plays a role in the continuity results.

\begin{definition} On the set of paths $PM$, we define the distance
$d_1$ by
$$d_1(c_1,c_2)=d_\infty(c_1,c_2) + |\ell(c_1)-\ell(c_2)|.$$
\end{definition}

\begin{proposition} The topology induced on $PM$ by the distance $d_1$
does not depend on the Riemannian metric chosen on $M$.
\end{proposition}

\pf. By compactness of $M$, two different metrics induce two equivalent
Riemannian distances on $M$ and thus two equivalent distances $d_1$ on $PM$.
\qed

\section{Approximation of embedded paths}
\label{embedded}
We want to extend the definition of the random holonomy to all paths in $PM$,
by approximation. Since any path of $PM$ is, by definition, a concatenation of
embedded submanifolds of $M$, it is natural to begin with those paths who are
embedded submanifolds themselves.

\subsection{Tubular neighbourhoods and Fermi coordinates}

These paths have the following nice property: they possess a tubular
neighbourhood that can be described using Fermi coordinates. Let us
fix a path $c$ which is an embedded submanifold. The proof of the
following result can be found in \cite{Gray}. Let us fix a
parametrization of $c$ and a vector field $N$ along $c$, unitary and
normal to $c$. 

\begin{proposition} There exists a positive real number $r$ such that
the mapping
\begin{eqnarray*}
\psi:[0,1] \times [-r,r] &\lra& M \\
(t,s)\hskip .8cm &\longmapsto& \exp_{c(t)}(s N_{c(t)})
\end{eqnarray*}
is a diffeomorphism onto its image, which is called tubular neighbourhood of
$c$ or tube around $c$. The coordinates $(t,s)$ are
called Fermi coordinates. They satisfy:\\
\indent 1. For any fixed $t_0$, the curve $s\mapsto \psi(t_0,s)$ is a piece of
geodesic normal to $c$. \\
\indent 2. For any couple $(t,s)\in [0,1]\times [-r,r]$,
$d(\psi(t,s),c([0,1]))=s.$
\end{proposition}

We shall always assume that the radius of the tubular neighbourhoods that we
consider are smaller than the convexity radius $R_M$ of $M$, defined
in \ref{injectivity_radius}. 

\subsection{Piecewise geodesic approximation}

The path $c$ is fixed until the end of the next section, together with a
tubular neighbourhood of radius $r$. 
\begin{proposition} \label{pg_approximation}
Let $x,y,z$ be three real numbers such that $0\leq
x<y<z\leq r$. There exists a piecewise geodesic path $\sigma$ such that
\\
\indent 1. $\sigma(0)=\psi(0,y)$ and $\sigma(1)=\psi(1,y)$, \\
\indent 2. $\sigma((0,1)) \subset \psi((0,1)\times (x,z))$, \\
\indent 3. $\sigma$ is injective.
\end{proposition}

We construct $\sigma$ as an approximation of the path $c_y:t\mapsto
\psi(t,y)$, in the same way as one would approximate a curve in $\R^2$ by
piecewise linear paths. 

\begin{lemma} \label{velocity_bounded} 
  Set $\displaystyle \delta_n(c_y)=\sup_{0\leq k \leq n-1}
  d\left(c_y\left({k\over {n}}\right),c_y\left({{k+1}\over
        {n}}\right)\right).$ Then $$\delta_n(c_y) \build{\lra}_{n\to
  \infty}^{} 0.$$
\end{lemma}

\pf : The norm of the velocity of $c$ is bounded. \qed

For $n$ large enough and for each $k=1\ldots n$, the points
$c_y({{k-1}\over n})$ and $c_y({{k}\over n})$ are close enough to be
joined by a unique minimizing geodesic $\zeta_{n,k}$ that stays at a distance
smaller than $\delta_n$ of $c_y({k\over n})$. We will always assume
that $n$ is large enough for this property to be true and set
$\zeta_n=\zeta_{n,1}\ldots \zeta_{n,n}$. 

\begin{lemma} For $n$ large enough, $\zeta_n$ is the graph of a
continuous function in Fermi coordinates. More precisely, there
exists a continuous function $\varphi: [0,1] \lra (-r,r)$ such that for
each $t\in[0,1]$, $\zeta_n(t)=\psi(t,\varphi(t))$.
\end{lemma}

\pf. It is sufficient to prove that each $\zeta_{n,k}$ is the graph of a
continuous function defined on $[{{k-1}\over n},{k\over n}]$ and that these
functions can be put together to form $\varphi$. 
Let $n$ be such that $\delta_n < {{r-y}\over 2}$. We show
that $\zeta_n$ stays inside the tubular neighbourhood of $c$. 

The first point is that $\zeta_{n,k}$ cannot meet the horizontal
boundary $\psi([0,1]\times \{-r,r\})$, because this boundary is at distance
$r$ of $c$ and $\zeta_{n,k}$ stays at distance smaller than
$\delta_n+y<r$. 

The vertical part of the boundary $\psi(\{0,1\}\times [-r,r])$ is made of two
pieces of minimizing geodesics, so that $\zeta_{n,k}$, which is also
minimizing, cannot meet twice one of these pieces without belonging to the
same infinite geodesic. This is impossible because the geodesics supporting
the vertical boundary meet $c_y$ only once.

The only way $\zeta_{n,k}$ could exit the tube around $c$ would be to exit
through one piece of the vertical boundary and get back through the other.
Suppose that $n$ is large enough for $\delta_n$ being smaller than ${1\over 2}
d(\psi(\{0\}\times [-r,r]),\psi(\{1\}\times [-r,r]))$. Then the situation
described above cannot happen, since any two points of the image of
$\zeta_{n,k}$ are at distance smaller than $\delta_n$. So, $\zeta_{n,k}$ stays
inside the tube.

Each vertical slice $t=t_0$ of the tube is a minimizing piece of a geodesic
that meets $c$ only once inside the tube, so that it meets $\zeta_{n,k}$ at
most one time. Thus, $\zeta_{n,k}$ is the graph of a smooth function defined
on the segment $[{{k-1}\over n},{k\over n}]$, equal to $y$ at both end points
of this segment. All these functions can be put together to make $\varphi$,
which is continuous. \qed

\pf \textsl{of proposition \ref{pg_approximation}}. Choose $n$ be large enough
for $\zeta_n$ to be the graph of a function in Fermi coordinates and such that
$\delta_n<\inf(z-y,y-x)$. As a graph, $\zeta_n$ is necessary injective, which
is statement $3$. The inequality $|\varphi(t)-y|=d(\zeta_n(t),c_y) \leq
\delta_n$ shows that $\zeta_n$ stays in $\psi([0,1]\times (x,z))$.  Together
with the fact that $\zeta_n$ meets at most once each vertical boundary, this
gives statement $2$. Statement $1$ is a direct consequence of the definition
of $\zeta_n$. So $\sigma=\zeta_n$ has all the properties required.\qed

\section{Random holonomy along embedded paths}

We suppose that the surface of the tube is smaller than the constant $s$ given
by the proposition \ref{fundamental_estimate}. We prove that the random
holonomy along a piecewise geodesic approximation of $c$ converges in
probability to a random variable and that this limit does not depend on the
particular choice of the approximation.

\subsection{Existence of a limit random holonomy}

For $n\geq 0$, set $x_n={r\over {2^{n+1}}}, \; y_n={3\over
  2}{r\over{2^{n+1}}}, \; z_n={r\over {2^{n}}}$ and let $\sigma_n$ be a path
given by the proposition \ref{pg_approximation}.  For each $n\geq 0$, let
$\lambda_n$ denote the vertical segment joining $(0,0)$ to $(0,y_n)$ and
$\rho_n$ the vertical segment joining $(1,y_n)$ to $(1,0)$. Finally, set
$\alpha_n=\lambda_n \sigma_n \rho_n$.

\begin{proposition} The sequence of random variables
$(H_{\alpha_n})_{n\geq 0}$ is a Cauchy sequence with respect to the
distance $d_P$.
\end{proposition}

\pf. Let $m\geq n$ be two integers. We want to estimate
$d_P(H_{\alpha_m},H_{\alpha_n})=d_P(H_{\alpha_n^{-1}\alpha_m},1)$. But
$\alpha_n^{-1}\alpha_m $ is equivalent to a simple loop which is the boundary
of an open set contained in $\psi([0,1]\times [0,{r\over {2^n}}])$. Thus, the
assumption on the surface of the tube allows us to apply proposition
\ref{fundamental_estimate}. We get:
$$ d_P(H_{\alpha_m},H_{\alpha_n}) \leq C \sigma\left(\psi\left([0,1]\times
\left[0,{r\over {2^n}}\right]\right)\right) \leq {{C'}\over {2^n}}.$$
This proves the result. \qed

\begin{figure}[h]
\begin{center}
\input{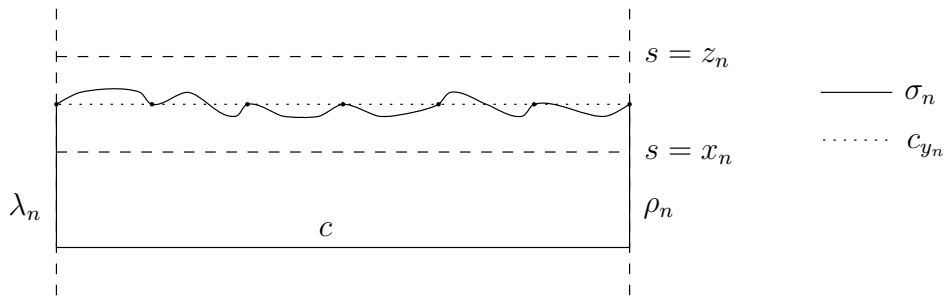}
\end{center}
\caption{Definition of the sequence $(\alpha_n)$.
\label{app1}}
\end{figure}

The space of $G$-valued random variables on $(\Omega,\a,P_\beta)$ endowed with
the distance $d_P$ is complete: it can be isometrically embedded in a $L^1$
space by embedding $G$ isometrically in some $\R^n$. So the sequence
$(H_{\alpha_n})_{n\geq 0}$ has a limit that we denote by $H_c$, anticipating
the fact that this limit does not depend on the choice of the sequence
$(\alpha_n)$.

\subsection{Unicity of the limit random holonomy}

\begin{lemma} \label{injective_paths} For all $\epsilon>0$, there
  exists $\delta>0$ such that if $\alpha$ is an {\rm injective} piecewise
  geodesic path with the same end points as $c$, such that $\alpha((0,1))
  \subset \psi((0,1) \times (-r,r))$ and such that
  $d_\infty(c,\alpha)<\delta$, then $d_P(H_c,H_\alpha)<\epsilon$.
\end{lemma}

In this statement, it is not necessary to control $|\ell(\alpha)-\ell(c)|$
because $\alpha$ is assumed to be injective.\\

\pf. Let $C$ be the constant given by the proposition
\ref{fundamental_estimate}. Let $n$ be such that
$d_P(H_c,H_{\alpha_n})<\epsilon/2$ and $C\sigma(\psi([0,1]\times [-{r\over
  {2^n}},{r\over{2^n}}]))<\epsilon/2$. Set $\delta={1\over {2^{n+1}}}$ and
suppose that $d_\infty(c,\alpha)<\delta$. Then $\alpha$ meets $\alpha_n$ only
at its end points. Thus $\alpha_n\alpha^{-1}$ is the boundary of an open set
included in $\psi([0,1]\times [-{r\over {2^n}},{r\over {2^n}}])$, so that
\begin{eqnarray*}
d_P(H_c,H_\alpha) &\leq& d_P(H_c,H_{\alpha_n}) + d_P(H_{\alpha_n},H_\alpha) \\
&\leq& {\epsilon \over 2} + C\sigma\left(\psi\left([0,1]\times \left[-{r\over
        {2^n}},{r\over 
  {2^n}}\right]\right)\right)\\
&\leq & \epsilon.
\end{eqnarray*}
\vskip -.8cm\qed

The control of the length of $\alpha$ allows to drop all restrictive
conditions on $\alpha$, unless those concerning end points. The main
result of this section is the following:

\begin{proposition} \label{pg_continuity}
For all $\epsilon>0$, there exists $\delta>0$ such that if $\alpha$
is a piecewise geodesic path with the same end points as $c$ and such
that $d_1(c,\alpha)<\delta$, then $d_P(H_c,H_\alpha)<\epsilon$.
\end{proposition}

\begin{lemma}
\label{golab}
For all $\delta>0$, there exists $\eta>0$ such that if $c'$ is another path
such that $d_\infty(c,c')<\eta$, then $\ell(c') > \ell(c) - \delta.$
\end{lemma}

Note that this result is not symmetric in $c$ and $c'$. Indeed,
it is true that there exists $\eta'$ such that $d_\infty(c,c')<\eta'$
implies $\ell(c)>\ell(c')-\delta$, but $\eta'$ may be much smaller than
$\eta$ (consider for $c'$ a zigzag approximating a straight line for
example). One could reformulate this result by saying that for any sequence
$c_n$ converging uniformly to
$c$, $\liminf \ell(c_n) \geq \ell(c)$.  \\

\pf : Let $n$ be large enough for the following inequality to hold:
$$\ell(c)-\sum_{k=0}^{n-1} d(c(\kn),c(\kpn)) < {\delta \over 2}.$$
Let $c'$ be a path with fixed parametrization such that $d_\infty(c,c')<
{\delta \over {4n}}$. Then, on one hand, $\displaystyle \ell(c') \geq \sum
d(c'(\kn),c'(\kpn)).$ On the other hand,
\begin{eqnarray*}
d(c(\kn),c(((k+1))/n)) &\leq & d(c(\kn),c'(\kn)) + d(c'(\kn),c'(\kpn)) +
d(c'(\kpn),c(\kpn)) \\ 
&\leq& {\delta \over {2n}} + d(c'(\kn),c'(\kpn)). 
\end{eqnarray*}
Thus, $\displaystyle \ell(c') \geq \sum_{k=0}^{n-1} d(c(\kn),c(\kpn)) -
{\delta \over 2} \geq \ell(c) - \delta$. We see that $\eta={\delta \over
  {4n}}$ is a possible choice. \qed

\pf \textsl{of proposition \ref{pg_continuity}.} Denote by $\delta_0$ the
distance between $c(0)$ and $c(1)$. Assume that $d_\infty(\sigma,c)$ is
smaller than $\inf(r,\delta_0/5)$. Recall that $r$ is assumed to be smaller
than the convexity radius of $M$ (see \ref{injectivity_radius}).

The points $\alpha(0)$ and $\alpha(1)$ are respectively in the balls
$B_0=B(c(0),2d_\infty(\alpha,c))$ and $B_1=B(c(1),2d_\infty(\alpha,c))$. These
balls are disjoint, hence there exists a last time $\tau_0$ at which $\alpha$
exits $B_0$ and a first time $\tau_1$ at which it enters $B_1$. The points
$\alpha(\tau_0)$ and $\alpha(\tau_1)$ are necessarily inside the tube, for the
points of $M$ at distance smaller than $r$ of $c$ are inside the tube or in
$B_0 \cup B_1$.  In Fermi coordinates, we can write $\alpha(\tau_0)=(t_0,s_0)$
and $\alpha(\tau_1)=(t_1,s_1)$. Note that $t_0>0$ and $t_1<1$: otherwise, we
would have $|s_0|$ or $|s_1|$ equal to $2d_\infty(\alpha,c)$.

Let $\gamma_0$ be the path that follows $c$ from time 0 to $t_0$ and then the
geodesic normal to $c$ from $(t_0,0)$ to $(t_0,s_0)$. Similarly, let
$\gamma_1$ be the path that follows the normal geodesic from $(t_1,s_1)$ to
$(t_1,0)$ and then $c$ from time $t_1$ to 1.  We write $\alpha$ in the
following way:
$$\alpha \simeq \alpha_{|[0,\tau_0]} \gamma_0^{-1} \cdot \gamma_0
\alpha_{|[\tau_0,\tau_1]} \gamma_1 \cdot \gamma_1^{-1} \alpha_{|[\tau_1,1]}.$$

The first and the third terms are small loops that we shall study later. Let
us consider the central term $\gamma_0 \alpha_{|[\tau_0,\tau_1]}\gamma_1$. It
is contained in the tube around $c$ and has the same end points as $c$. Let us
decompose it according to \ref{decomposition_in_lassos} into a product
$\lambda_1\ldots \lambda_p \xi$, where the $\lambda_i$'s are lassos based at
$c(0)$ and $\xi$ is an injective path between $c(0)$ and $c(1)$. It is obvious
that $d_\infty(c,\xi) \leq d_\infty(c,\alpha)$. This tells us, by proposition
\ref{injective_paths}, that $H_\xi$ can be made arbitrarily close to $H_c$ by
taking $d_\infty(c,\alpha)$ sufficiently small.

Let us fix a positive $\epsilon$ and $\delta_1$ such that
$d_\infty(c,\alpha)<\delta_1$ implies $d_P(H_c,H_\xi)<\epsilon/2$. It is
enough now to control $d_P(H_\xi,H_{\alpha})$. 
\begin{eqnarray*}
d_P(H_\xi,H_{\alpha}) &=& d_P(H_\xi,H_{
\alpha_{|[0,\tau_0]} \gamma_0^{-1}} H_{\lambda_1} \ldots H_{\lambda_p}
H_\xi H_{\gamma_1^{-1} \alpha_{|[\tau_1,1]}}) \\
& \leq & \sum_{i=1}^p d_P(H_{\lambda_i},1) + d_P(H_{\alpha_{|[0,\tau_0]} \gamma_0^{-1}},1) + d_P(H_{\gamma_1^{-1} \alpha_{|[\tau_1,1]}},1).
\end{eqnarray*}

We are led to consider the random variables associated with loops with finite
self-intersection. According to \ref{small_lassos}, it it is sufficient to
control their lengths. We already know by \ref{golab} that we can have
$\ell(\xi)\geq \ell(c)-\varepsilon/8$ provided $\delta_1$ and so
$d_\infty(c,\xi)$ is small enough. If we impose now that
$d_1(c,\alpha)<\delta_1$, instead of $d_\infty(c,\alpha)<\delta_1$, then we
also get $\ell(\alpha)<\ell(c)+\epsilon /8$.

Then $0<\ell(\alpha) - \ell(\xi)<\varepsilon/4$.
Let us denote by $\beta_1,\ldots,\beta_p$ the buckles of the lassos
$\lambda_1,\ldots,\lambda_p$. By \ref{decomposition_in_lassos},
$$\ell(\xi) + \sum_i \ell(\beta_i) \leq \ell(\gamma_0) +
\ell(\alpha_{|[\tau_0,\tau_1]})+\ell(\gamma_1),$$
and so 
$$\sum_i \ell(\beta_i) + \ell(\alpha_{|[0,\tau_0]} \gamma_0^{-1}) +
\ell(\gamma_1^{-1} \alpha_{|[\tau_1,1]}) \leq {\epsilon \over 4} +
2\ell(\gamma_0) + 2\ell(\gamma_1).$$
Since $\ell(\gamma_i) \leq
2d_\infty(c,\alpha) + \ell(c([0,1])\cap B_i)$, the lengths appearing in the
right hand side can be made small by taking $d_\infty(c,\alpha)$ small enough.
This is exactly what was needed to control $d_P(H_\xi,H_\alpha)$.  This gives
us a $\delta_2$ such that $d_1(c,\alpha)<\delta_2$ implies
$d_P(H_c,H_{\alpha})<\epsilon$. \qed

\begin{corollary} \label{seq-cont}
Let $(\beta_n)_{n\geq 0}$ be any sequence of
piecewise geodesic paths with the same end points as $c$ that converges
to $c$. Then the sequence $(H_{\beta_n})$ converges to $H_c$.
\end{corollary}

This proves that the variable $H_c$ does not depend on the particular choice
of the sequence of paths approximating $c$.

\subsection{Continuity of the double layer potential (1)}
\label{dlp_reg}
 
Following step by step the proofs of propositions \ref{injective_paths} and
\ref{pg_continuity} and replacing statements about random variables by
statements about the double layer potential, according to the remark made at
the end of paragraph \ref{correspondance}, we get the following results:
 
\begin{lemma} For all $\epsilon>0$, there
exists $\delta>0$ such that if $\alpha$ is an injective piecewise
geodesic path with the same end points as $c$, such that
$\alpha((0,1)) \subset \psi((0,1) \times (-r,r))$ and such that
$d_\infty(c,\alpha)<\delta$, then $\parallel u_\alpha -
u_c\parallel_{L^2} <\epsilon$. 
\end{lemma}

\begin{proposition} \label{pg_continuity_dlp}
For all $\epsilon>0$, there exists $\delta>0$ such that if $\alpha$
is a piecewise geodesic path with the same end points as $c$ and such
that $d_1(c,\alpha)<\delta$, then \linebreak $\parallel u_\alpha -
u_c\parallel_{L^2} <\epsilon$.
\end{proposition}

\section{Random holonomy along arbitrary paths}
\label{arbitrary}
\subsection{Construction of the random holonomy}

Let $c$ be an element of $PM$. By definition, it can be written $c=c_1\ldots
c_p$, where the $c_i$'s are embedded paths, but this decomposition is far to
be unique. Nevertheless, we prove that the random variable $H_{c_p}\ldots
H_{c_1}$ depends only on $c$.

\begin{lemma} Let $c$ be a path. There exists a sequence
of piecewise geodesic paths with the same end points as $c$ that
converges to $c$.
\end{lemma}

\pf. Given a decomposition $c=c_1\ldots c_p$ of $c$ into a product of embedded
paths, we concatenate sequences of paths that converge to each $c_i$ with
fixed end points and get the required sequence. \qed

\begin{proposition}  \label{cutting} 
  Let $(\alpha_n)$ be a sequence of piecewise geodesic paths that converges
  with fixed end points to $c$. The sequence $(H_{\alpha_n})$ converges to the
  product $H_{c_p} \ldots H_{c_1}$.
\end{proposition}

The following corollary is in fact the main result of this paragraph. 

\begin{corollary} \label{cutting_2}
  The product $H_{c_p} \ldots H_{c_1}$ is independent of the choice of the
  decomposition of $c$ and it is equal to the common limit of all sequences
  $(H_{\alpha_n})$ associated with sequences $(\alpha_n)$ of piecewise
  geodesic paths converging to $c$ with fixed end points. We shall denote it
  by $H_c$.
\end{corollary}

\pf \textsl{of proposition \ref{cutting}.} We cut $\alpha_n$ in a way that
corresonds to the decomposition of $c$. Let us fix a parametrization of $c$
such that $c_i=c_{|[{{i-1}\over p},{i\over p}]}$. Let us also fix a
parametrization of each $\alpha_n$ such that the uniform convergence
$d(c,\alpha_n)\to 0$ holds with these parametrizations. Set
$\alpha_{i,n}={\alpha_n}_{|[{{i-1}\over p},{i\over p}]}$. Let us show that for
each $i$, $\alpha_{i,n} \build{\lra}_{n\to\infty}^{} c_i$.

The first point is that $d_\infty(c_i,\alpha_{i,n}) \leq d_\infty(c,\alpha_n)
\build{\lra}_{n\to\infty}^{} 0$, parametrizations being fixed. Now let us
choose $\epsilon>0$ and $n$ large enough such that for all $i=1,\ldots,p$,
$\ell(\alpha_{i,n}) \geq \ell(c_i) - {\epsilon \over p}$ (using lemma
\ref{golab}) and $\ell(\alpha_n) \leq \ell(c) +\epsilon$. Then
\begin{eqnarray*}
-{\epsilon \over p} \leq \ell(\alpha_{i,n}) - \ell(c_i) &=& \ell(\alpha_n)
 -\ell(c) - \sum_{j\neq i} \left( \ell(\alpha_{j;n})-\ell(c_j) \right) \\
&\leq& 2 \epsilon.
\end{eqnarray*}
So we also have $\ell(\alpha_{i,n})\build{\lra}_{n\to\infty}^{} \ell(c_i)$ for
all $i$.  Consider now for each $i$ and each $n$ the path
$\widetilde\alpha_{i,n}$ which is $\alpha_{i,n}$ concatenated at both end
points with a minimizing piece of geodesic in order to have the same end
points as $c_i$.  If $n$ is large enough, $\alpha_n$ is close enough to $c$
for these pieces of minimizing geodesic to be uniquely defined. So, these
geodesic segments cancel out in the product $\widetilde\alpha_{1,n} \ldots
\widetilde\alpha_{p,n}$ which is equivalent to $\alpha_{1,n} \ldots
\alpha_{p,n}$.  On the other hand, we have $\widetilde\alpha_{i,n}
\build{\lra}_{}^{d_1} c_i$.  Indeed, the small geodesic pieces stay close to
each $c_i$ and their length converges to zero.  Thus, the corollary
\ref{seq-cont} implies $H_{\widetilde\alpha_{i,n}}
\build{\lra}_{n\to\infty}^{} H_{c_i}$ which gives the result:
$$H_{\alpha_n}=H_{\widetilde\alpha_{p;n}} \ldots H_{\widetilde\alpha_{1;n}}
\build{\lra}_{n\to\infty}^{} H_{c_p} \ldots H_{c_1}.$$
\vskip -0.6cm \qed

\subsection{Continuity of the random holonomy}
At this point, we constructed a random holonomy along each path $c$ of $PM$.
This random holonomy is a $G$-valued random variable on the probability space
$(\Omega,\a,P(x_1,\ldots,x_q))$. Let us state some of its basic properties.

\begin{proposition} \label{multiplicative}
  Let $c,c_1,c_2$ be elements of $PM$. \\
  \indent  1. $H_{c^{-1}}=H_c^{-1}$ a.s. \\
  \indent 2. The random variable $H_c$ depends only on the equivalence class
  of $c$
  for the relation $\simeq$. \\
  \indent 3. If $c_1$ and $c_2$ satisfy $c_1(1)=c_2(0)$ then $H_{c_1 c_2}
  =H_{c_2} H_{c_1} $ a.s.
\end{proposition}

\pf : Property (3) is obvious by putting together two decompositions of $c_1$
and $c_2$ and the property (2) is a direct consequence of (1). To prove (1),
just note that this is true for piecewise geodesic paths by construction, and
that if $\alpha_n \build{\lra}_{n\to\infty}^{d_1} c$, then $\alpha_n^{-1}
\build{\lra}_{n\to\infty}^{d_1} c^{-1}$. \qed

We still have to prove that the law of this random holonomy does not depend on
choice of the Riemannian metric used in the construction. For this, we need a
regularity property which is the object of the next proposition.

\begin{proposition} \label{continuity_holonomy}
  Let $c$ be a path of $PM$. For any $\epsilon>0$, there exists $\delta>0$
  such that if $c'$ is another path of $PM$ with the same end points as $c$
  and if $d_1(c,c')<\delta$, then $d_P(H_c,H_{c'})<\epsilon$.
\end{proposition}

\pf. Let $\delta_0>0$ be given by the proposition \ref{pg_continuity} such
that for any piecewise geodesic path $\alpha$ with the same end points as $c$,
$d_1(c,\alpha)<\delta_0$, implies $d_P(H_c,H_\alpha)<\epsilon$.  Let
$\delta={{\delta_0}\over 2}$. Suppose that $c'$ is a path of $PM$ with the
same end points as $c$ such that $d_1(c,c')<\delta$.  Let $\alpha$ be a
piecewise geodesic path such that, simultaneously, $d_1(\alpha,c')<\delta$ and
$d_P(H_{c'},H_\alpha)<{\epsilon \over 2}$. Then $d_1(c,\alpha)<\delta_0$, so
that
$$d_P(H_c,H_{c'}) \leq d_P(H_c,H_\alpha) + d_P(H_\alpha,H_{c'}) < \epsilon.$$
\vskip -.6cm \qed

Let us state a result that summarizes the results of the procedure of
piecewise geodesic approximation. It is in fact the center of the continuum
limit procedure. We put together the propositions \ref{cutting},
\ref{cutting_2} and \ref{continuity_holonomy}.

\begin{proposition} \label{summarize_approx} Let $c$ be a path of $PM$. For
  any sequence $(\alpha_n)_{n\geq 0}$ of piecewise geodesic paths converging
  to $c$ with fixed end points, the sequence $(H_{\alpha_n})_{n\geq0}$
  converges to a random variable that depends only on $c$ and that we denote
  by $H_c$. Moreover, for any $\epsilon>0$, there exists $\delta>0$ such that
  if $c'$ is another path of $PM$ with the same end points as $c$ and if
  $d_1(c,c')<\delta$ then $d_P(H_c,H_{c'})<\epsilon$.
\end{proposition}

\subsection{Continuity of the double layer potential (2)}

One more time, we transpose directly the preceding arguments to the double
layer potential and get the following result:

\begin{proposition} \label{continuity_paths_dlp}
  Let $c$ be a path of $PM$. For any $\epsilon>0$, there exists $\delta>0$
  such that if $c'$ is another path of $PM$ with the same end points as $c$
  and if $d_1(c,c')<\delta$, then $\parallel u_c-u_{c'}
  \parallel_{L^2}<\epsilon$.
\end{proposition}

\begin{corollary} \label{continuity_loops_dlp}
  Let $l$ be a loop of $LM$. For any $\epsilon>0$, there exists $\delta>0$
  such that if $l'$ is another loop of $LM$ and if $d_1(l,l')<\delta$, then
  $\parallel u_l - u_{l'} \parallel_{L^2}<\epsilon$.
\end{corollary}

\pf. Let $\delta_0$ be given by the preceding proposition and set
$\delta={{\delta_0}\over 3}$. Let $l'\in LM$ be such that $d_1(l,l')<\delta$.
We have in particular $d(l(0),l'(0))<\delta$. Let $\sigma$ be a minimizing
geodesic from $l(0)$ to $l'(0)$. Then $\widetilde l'=\sigma l' \sigma^{-1}$
satisfies $u_{\widetilde l'}=u_{l'}$ a.e. and $d_1(\widetilde l',l)<\delta_0$.
Moreover, $\widetilde l'$ has the same end points as $l$. Thus,
$$\parallel u_l - u_{l'} \parallel_{L^2} = \parallel u_l - u_{\widetilde l'}
\parallel_{L^2} < \epsilon.$$
\vskip -.6cm \qed

\section{Law of the random holonomy}

For the moment, we are only able to write down the law of the holonomy along
piecewise geodesic paths. We want to show that the law of the holonomy along
arbitrary families of paths with finite self-intersection is what we expect it
to be, namely that given by the discrete theory. The goal of this section is
to prove the following proposition:

\begin{proposition} \label{law_graph}
  Let $\Gamma=\{a_1,\ldots,a_r \}$ be a graph on $M$ such that $L_1,\ldots L_q
  \in \Gamma^*$. For any function $f$ continuous on $G^\Gamma$, we have:
  $$Ef(H_{a_1},\ldots,H_{a_r})=\int_{G^\Gamma} f\; dP^\Gamma_\beta.$$
\end{proposition}

A very important consequence of this result is the independence of the
construction with respect to the Riemannian metric:

\begin{corollary} The law of the family $(H_c)_{c\in PM}$ does not
  depend on the choice of the Riemannian metric that was used throughout the
  construction.
\end{corollary}

\pf. Consider two families of variables obtained with two different choices of
metric. By the preceding proposition, these families have the same law on the
set of paths that are piecewise geodesic for, say, the first metric. By
proposition \ref{continuity_holonomy}, both families are continuous in a sense
that is strong enough to guarantee that their laws coincide on the whole set
$PM$.  \qed

In order to prove the proposition \ref{law_graph}, we need a technical result
about the approximation of graphs by piecewise geodesic graphs. Before to
state this result, let us make some remarks about the edges and faces in a
graph in $M$.

Recall that a path and hence an edge must by definition have non-zero
derivatives at its end points. This avoids pathological behaviours. For
example, consider all edges that share a given vertex of a graph and a small
geodesic circle centered at this vertex.  If the radius of this circle is
small enough, each edge cuts it only once, and the order of the intersection
points, which does not depend on the radius of the circle, defines a cyclic
order on the set of these edges.

\begin{figure}[hbtp]
\begin{center}
\scalebox{1}{\input{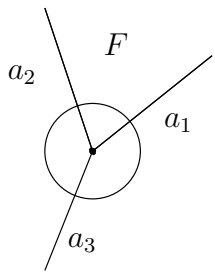}}
\end{center}
\caption{The face determined by two adjacent edges.
\label{face}}
\end{figure}

Now, consider two edges that are adjacent for this order. They bound at least
one common face. Thus, if $M$ is oriented, a couple of adjacent edges
determines a face of the graph (see fig. \ref{face}). Conversely, given a
face, any two consecutive edges of the boundary of this face are adjacent at
the vertex that they share, or eventually at both vertices if they share two.

\begin{proposition} \label{approx_graphs}
  Let $\Gamma=\{a_1,\ldots,a_r\}$ be a graph such that $L_1,\ldots,L_q \in
  \Gamma^*$. For any $\epsilon>0$, there exists a graph
  $\Gamma_\epsilon=\{\alpha_1,\ldots,\alpha_r\}$ with piecewise
  geodesic edges such that:\\
  \indent  1. $\Gamma_\epsilon$ and $\Gamma$ have the same vertices,\\
  \indent  2. For each $i=1,\ldots, r$, $d_1(\alpha_i,a_i)<\epsilon$, \\
  \indent  3. $L_1,\ldots,L_q \in \Gamma_\epsilon^*$, \\
  \indent  4. For each $i=1,\ldots,r$, $a_i$ and $\alpha_i$ are in the same
  connected component of the complementary of the unions of the images
  of the $L_i$'s.\\
  Let us denote by $\alpha:\Gamma^* \lra \Gamma_\epsilon^*$ the multiplicative
  map that sends $a_i$ to $\alpha_i$. It is possible to construct
  $\Gamma_\epsilon$ in such a way that this map induces a one-to-one
  correspondence still denoted by $\alpha : \f(\Gamma) \lra
  \f(\Gamma_\epsilon)$ such that $\partial \alpha (F) = \alpha(\partial F)$
  and $\sigma(F - \partial F)<\epsilon$, where $-$ denotes the symmetric
  difference.
\end{proposition}

\pf. The property 4 is a consequence of 2 and 3. Indeed, if $a_i$ is in
a given connected component, $\alpha_i$ meets this component if $\epsilon$ is
small enough, by 2. But $\alpha_i$ could only exit this component by
crossing $L_i$ at a point which is not an end point of $\alpha_i$, which is
impossible by 3 and by the definition of graphs.

Let $\v(\Gamma)=\{s_1,\ldots,s_p\}$ denote the set of vertices of $\Gamma$.
Let $r$ be a positive real number that "localizes the vertices
of $\Gamma$", i.e. small enough to satisfy the following properties: \\
\indent 1. The balls $B(s_i,r)$ are pairwise disjoint.\\
\indent 2. For every pair $(a_i,s_j)$ with $a_i\in \Gamma$ and $s_j$ an end
point of $a_i$, $a_i$ meets only once and transversally any circle centered at
$s_j$ and of radius smaller than $r$. Moreover, the length of the portion of
$a_i$ in the corresponding ball is smaller than $\epsilon/16$. \\
\indent 3. For any pair $(a_i,s_j)$ where $s_j$ is not an end point of $a_i$,
$a_i$ does not meet the ball $B(s_j,r)$.\\
\indent 4. The sum of the surfaces of the ball $B(s_i,r)$ is smaller than
$\epsilon/2$.\\
\indent 5. $r< \epsilon/16$ and $r<R_M$, where $R_M$ is the convexity radius of $M$.\\
All properties remain true for $r'<r$ once they are true for $r$, so that it
is not a problem to get them simultaneously.

Let $t$ be a positive real number such that
$\sigma(\{d(\cdot,\Gamma)<t\})<{\epsilon \over 2}$.  Let $\widehat a_i$ denote
the portion of $a_i$ outside the disks of radius $r$ around its end points.
Let $\delta$ be the smallest distance between the images of two distinct
$\widehat a_i$.  For each $i$, let $\gamma_i$ be an injective piecewise
geodesic path with the same end points as $\widehat a_i$, such that
$d_1(\gamma_i,\widehat a_i)<\inf(\epsilon/4, \delta/2,t)$ and that never meets
the balls $B(s_j,r)$, except at its ends. This last condition can be obtained
because $a_i$ cuts $B(s_j,r)$ transversally: in a neighbourhood of each end
point of $\widehat a_i$, there is a half-tube around $\widehat a_i$ that does
not meet $B(s_j,r)$. It is possible around each end point of $\widehat a_i$ to
construct $\gamma_i$ inside this half-tube.  By definition of $\delta$, the
$\gamma_i$'s are disjoint.

\begin{figure}[hbtp]
\begin{center}
\scalebox{1}{\input{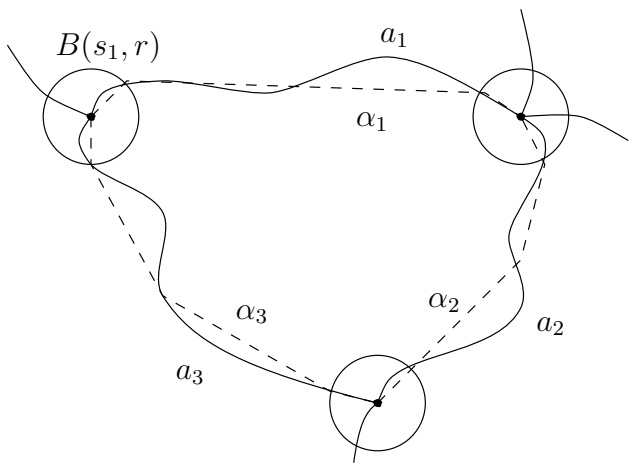}}
\end{center}
\caption{Definition of the edges of the graph $\Gamma_\epsilon$.
\label{graph}}
\end{figure}

Now define $\alpha_i$ for each $i$ such that $a_i$ is not piecewise geodesic
as the concatenation of the minimizing geodesic from $a_i(0)$ to $\widehat
a_i(0)$, of $\gamma_i$ and of the minimizing geodesic from $\widehat a_i(1)$
to $a_i(1)$ (see fig. \ref{graph}). Assumption 5 ensures that these minimizing
geodesics are well defined. Assumptions 2 and 5 imply that
$d_1(\alpha_i,a_i)<\epsilon$. The edges of the decompositions of the $L_i$'s
are already piecewise geodesic.  Hence we only rename them, setting
$\alpha_i=a_i$.

The $\alpha_i$'s are edges. Moreover, they were constructed in such a way that
they meet only at their ends: we already noticed that they do not meet outside
the balls around the vertices of $\Gamma$, and they cannot meet more than once
inside these balls according to the local properties of geodesics. Thus, the
graph $\Gamma_\epsilon=\{\alpha_1, \ldots, \alpha_r \}$ exists and has the
same vertices as $\Gamma$.

We just proved that properties 1 and 2.
Property 3 is true because we kept the edges corresponding to the $L_i$'s
and property 4 follows, according to the remark made at the beginning of the
proof.  It remains to prove the last part of the statement.

Consider edges of $\Gamma$ that share a given vertex. They are given a cyclic
order. By definition, the corresponding $\alpha_i$'s cut the circle of radius
$r$ around this vertex in the same order, so that the multiplicative
application $\alpha : \Gamma^* \to \Gamma_g^*$ defined by
$\alpha(a_i)=\alpha_i$ preserves the cyclic order at each vertex.

Given a pair of edges of $\Gamma$ that determine the face $F$, the pair of
corresponding edges of $\Gamma_\epsilon$ is a pair of adjacent edges that
determine a face of $\Gamma_\epsilon$. This face does not depend on the
particular choice of the edges that represent $F$ and we denote it by
$\alpha(F)$. By construction, we have the relation
$\partial(\alpha(F))=\alpha(\partial F)$.

The symmetric difference of $F$ and $\alpha(F)$ is contained in the reunion of
the balls $B(s_j,r)$ and the sets $\{d(\cdot,\widehat a_i)<t\}$. By assumption
4) and by definition of $t$, we know that the total volume of these sets is
smaller than $\epsilon$. Thus, $\vol (F - \partial F) < \epsilon$. Moreover,
this inequality characterizes $\alpha(F)$ among the faces of $\Gamma_\epsilon$
that have $\alpha(\partial F)$ as boundary, if there is more than one,
provided $\sigma(M)$ is greater than $2\epsilon$. \qed

\pf \textsc{ of \ref{law_graph}}. For each integer $n$, the
preceding proposition gives a graph $\Gamma_{1\over
  n}=\{\alpha_{1,n},\ldots,\alpha_{r,n}\}$. For each $i=1,\ldots,r$, the
sequence $(\alpha_{i,n})$ converges to $a_i$ with fixed end points, so that
$H_{\alpha_{i,n}} \lra H_{a_i}$. In particular, we have the convergence in
law:
$$(H_{\alpha_{1,n}},\ldots,H_{\alpha_{r,n}}) \build{\lra}_{n\to\infty}^{\law}
(H_{a_1},\ldots,H_{a_r}).$$
Thus, for any function $f$ continuous on
$G^\Gamma$,
\begin{eqnarray*}
Ef(H_{a_1},\ldots,H_{a_r}) &=& \lim_{n\to\infty} {1\over
{Z^{\Gamma_{1\over n}}}} \int_{G^\Gamma} f(g_1,\ldots,g_r) \prod_{F\in
\f(\Gamma_{{1\over n}})} p_{\sigma(F)}(h_{\partial F}) \\
&& \hskip 6.2cm d\nu_{x_1} \ldots d\nu_{x_q} \; dg' \\
&=& \left(\lim_{n\to\infty} {1\over{Z^{\Gamma_{1\over n}}}}\right) \int_{G^\Gamma}
f(g_1,\ldots,g_r) \prod_{F\in \f(\Gamma)}
p_{\sigma(F)}(h_{\partial F}) \\
&& \hskip 6.2cm d\nu_{x_1} \ldots d\nu_{x_q} \; dg',
\end{eqnarray*}
using the fact that $\sigma(\alpha(F))$ tends to $\sigma(F)$ when $n$
tends to infinifty. 

Recall from proposition \ref{Z_constant} that the conditional partition
functions computed in two graphs, one being finer than the other, are equal.
But given two piecewise geodesic graphs, there exists a third one which is
finer than both others, as was proved in \ref{finer}. Thus the partition
function is the same for all piecewise geodesic graphs, and the sequence
$\left(Z^{\Gamma_{1\over n}}\right)$ is constant. Its value can be computed by
setting $f$ identically equal to 1: we find that it is equal to $Z^\Gamma$.
This proves the result. \qed

By the way, we proved the following important result:

\begin{proposition} \label{inv_part_gr}
  Let $\Gamma$ be a graph such that $L_1,\ldots,L_q\in \Gamma^*$. Then the
  value of the conditional partition function $Z^\Gamma(x_1,\ldots,x_q)$ does
  not depend on $\Gamma$.
\end{proposition}

\section{Surfaces with boundary}
\label{boundary_case}

At the beginning of this chapter, we restricted ourselves for technical
reasons to surfaces without boundary. In this section, we will extend the
construction of the random holonomy to the case of surfaces with boundary.

\subsection{Natural law of the holonomy along the boundary} 

Let $(M,\sigma)$ be a surface with a boundary $\partial M=N_1\cup \ldots\cup
N_p$. In order to construct the holonomy along the paths of $M$, we shall
embed $M$ in a minimal closure and use the construction described in the
preceding sections. But if we want this procedure to give a result independent
of the closure of $M$, and we do, it is necessary to condition the holonomy
along every component of $\partial M$. If this was not our first intention,
say if we expected only to impose the holonomy along $N_1$, to be equal to $x$
for example, we need to know the natural law of the holonomy along the whole
boundary under $P(x)$. Then, we will artificially impose this natural law when
working on the closure of $M$. We begin by defining this natural law.

Let $L_1,\ldots,L_q$ be disjoint simple loops on $M$ whose image is included
in the interior of $M$. Let $N_1,\ldots,N_k$ be the components of $\partial M$
along which we want to impose the holonomy. Let
$x_1,\ldots,x_k,y_1,\ldots,y_q$ be elements of $G$. The $k$ first elements
correspond to the components of $\partial M$, the $q$ others to the interior
loops.

\begin{proposition} \label{natural_law}
  Let $\Gamma$ be a graph on $M$ such that $L_1,\ldots,L_q\in\Gamma^*$. The
  law of the random variable $(h_{N_1},\ldots,h_{N_p},h_{L_1},\ldots,h_{L_q})$
  defined on the probability space
  $(G^\Gamma,P(x_1,\ldots,x_k,y_1,\ldots,y_q))$ does not depend on $\Gamma$.
  We will denote it by $\beta(x_1,\ldots,x_k,y_1,\ldots,y_q)$.
\end{proposition}

\begin{lemma} \label{valid_boundary}
  The propositions \ref{approx_graphs} and \ref{inv_part_gr} hold on surfaces
  with boundary.
\end{lemma}

\pf. Note that \ref{approx_graphs} implies \ref{inv_part_gr}, using the
computation done at the end of its proof. Thus, it is sufficient to prove that
\ref{approx_graphs} holds. For this, embed $M$ in a minimal closure $M_1$. A
graph $\Gamma$ on $M$ induces a graph on $M_1$, which can be approximated by
piecewise geodesic graphs. If each component of $\partial M$ is the image of
one of the $L_i$'s, then the property (4) says exactly that the approximating
graphs stay inside $M$. \qed
 
\pf \textsl{of proposition \ref{natural_law}}. Let us endow $M$ with a
Riemannian metric for which $N_1,\ldots,N_p$, $L_1,\ldots,L_q$ are geodesics.
The law of $(h_{N_1},\ldots,h_{N_p},h_{L_1},\ldots,h_{L_q})$ does not depend
on $\Gamma$ provided it is piecewise geodesic, by invariance by subdivision.

Now consider an arbitrary graph $\Gamma$. According to the preceding lemma, we
can approximate it by piecewise geodesic graphs, for which the law we are
interested in is always the same. The convergence of the joint law of the
holonomy along all edges proves the result. \qed

The lemma \ref{partition_law} gives us an expression of
$\beta(x_1,\ldots,x_k,y_1,\ldots,y_q)$. We state it here again.

\begin{lemma} The following equality between measures on $G^{p+q}$
holds:
\begin{eqnarray*}
&&\hskip -.5cm \beta(x_1,\ldots,x_k,y_1,\ldots,y_q)= \\
&& \hskip .3cm \delta_{(x_1,\ldots,x_k)} \otimes
{{Z(x_1,\ldots,x_k,x'_{k+1},\ldots,x'_p,y_1,\ldots,y_q)}\over
{Z(x_1,\ldots,x_k,y_1,\ldots,y_q)}} \; dx'_{k+1}\ldots dx'_p \otimes
\delta_{(y_1,\ldots,y_q)}.
\end{eqnarray*}

\end{lemma}

In the particular case where we do not want to condition the measure
at all, the last expression is still true, with the convention that a
conditional partition function without parameters is equal to $1$.

\subsection{Definition of the random holonomy}

Let $M_1$ be a minimal closure of $M$ endowed with a surface measure
$\sigma_1$ that extends $\sigma$. We see $N_1,\ldots,N_p$,
$L_1,\ldots,L_q$ as loops on $M_1$. There is a measurable space
$(\Omega_1,\a_1)$ on which we constructed
a family of measurable functions $(H_c)_{c\in PM_1}$. On this
measurable space, we put the following probability:
$$P_1=\! \int_{G^{p+q}} \!\!\! P(x'_1,\ldots,x'_p,y'_1,\ldots,y'_q) \;
d(\beta(x_1,\ldots,x_k,y_1,\ldots,y_q))(x'_1,\ldots,x'_p,y'_1,\ldots,y'_q).$$
In other words, we insist on the law of
$(H_{N_1},\ldots,H_{N_p},H_{L_1},\ldots,H_{L_q})$ being the natural
one under $P(x_1,\ldots,x_k,y_1,\ldots,y_q)$.

We consider the restriction of the family $(H_c)_{c\in CM_1}$ to $M$, i.e. we
restrict the index set to $PM$.

\begin{proposition}
The law of the restriction $(H_c)_{c\in PM}$ does not depend on
$M_1$. If $\Gamma=\{a_1,\ldots,a_r\}$ is a graph on $M$ such that
$L_1,\ldots,L_q \in \Gamma^*$, then the law of
$(H_{a_1},\ldots,H_{a_r})$ under $P_1$ is the discrete Yang-Mills measure
$P_M(x_1,\ldots,x_k,y_1,\ldots,y_q)$.
\end{proposition}

\pf. The regularity property \ref{continuity_holonomy} of the random
holonomy on $M_1$ is still true for its restriction to $M$. Thus, the
second assertion implies the first one, using the fact that any family
of paths on $M_1$ can be approximated by piecewise geodesic families. 

Let $\Gamma$ be a graph as in the statement and $f$ be a continuous function
on $G^\Gamma$.
\begin{eqnarray*}
&&E_{P_1}f(H_{a_1},\ldots,H_{a_r}) = \int_{G^{p+q}}
{{d(\beta(x_1,\ldots,x_k,y_1,\ldots,y_q))}\over
{Z_{M_1}(x'_1,\ldots,x'_p,y'_1,\ldots,y'_q)}} \int_{G^\Gamma}
f(g_1,\ldots,g_r) \\
&& \hskip 3.8cm \prod_{F\in\f(\Gamma), F\subset M_1}
p_{\sigma_1(F)}(h_{\partial F})  \; d\nu_{x'_1} \ldots d\nu_{x'_p}
d\nu_{y'_1} \ldots d\nu_{y'_q} \; dg'.
\end{eqnarray*}
In $M_1$, the loops $N_1,\ldots,N_p$ bound $p$ disks
$D_1$,\ldots,$D_p$ which are the only faces of $\Gamma$ that are not
inside $M$. Thus,
\begin{eqnarray*}
Z_{M_1}(x'_1,\ldots,x'_p,y'_1,\ldots,y'_q)&=&\int_{G^\Gamma} \prod_{F
\in\f(\Gamma) \; F\subset M} p_{\sigma(F)}(h_{\partial F})
\prod_{i=1}^p p_{\sigma_1(D_i)}(h_{N_i})\\
&&\hskip 3.3cm d\nu_{x'_1} \ldots d\nu_{x'_p}
d\nu_{y'_1} \ldots d\nu_{y'_q} \; dg' \\
&=& \prod_{i=1}^p p_{\sigma_1(D_i)}(x'_i)
Z_M(x'_1,\ldots,x'_p,y'_1,\ldots,y'_q).
\end{eqnarray*}

Using this last relation, we get:
\begin{eqnarray*}
E_{P_1}f(H_{a_1},\ldots,H_{a_r}) &=& \int_{G^{p+q}}
{{d(\beta(x_1,\ldots,x_k,y_1,\ldots,y_q))}\over
{Z_M(x'_1,\ldots,x'_p,y'_1,\ldots,y'_q)}} \int_{G^\Gamma}
f(g_1,\ldots,g_r) D^\Gamma_M \\
&&\hskip 4.3cm d\nu_{x'_1} \ldots d\nu_{x'_p}
d\nu_{y'_1} \ldots d\nu_{y'_q} \; dg'\\
&& \hskip -4cm = \int_{G^{p-k}} \!\!\!
{{Z(x_1,\ldots,x_k,x'_{k+1},\ldots,x'_p,y_1,\ldots,y_q)}\over 
{Z(x_1,\ldots,x_k,y_1,\ldots,y_q)}} {{dx'_{k+1}\ldots dx'_p}\over
{Z(x_1,\ldots,x_k,x'_{k+1},\ldots,x'_p,y_1,\ldots,y_q)}} \\
&&\hskip -1.4cm \int_{G^\Gamma}
f(g_1,\ldots,g_r) D^\Gamma_M \; d\nu_{x_1} \ldots
d\nu_{x_k}d\nu_{x'_{k+1}}\ldots d\nu_{x'_p}
d\nu_{y_1} \ldots d\nu_{y_q} \; dg'\\
&& \hskip -3.5cm ={1 \over {Z(x_1,\ldots,x_k,y_1,\ldots,y_q)}} \int_{G^\Gamma} f(g_1,\ldots,g_r) D^\Gamma_M \; d\nu_{x_1} \ldots d\nu_{x_k}
d\nu_{y_1} \ldots d\nu_{y_q} \; dg'\\
&& \hskip -3.5cm =P(x_1,\ldots,x_k,y_1,\ldots,y_q)(f).
\end{eqnarray*}
\vskip -.6cm \qed

\section{Summary of the properties of the random holonomy}

\subsection{Existence, unicity in law and main properties}

Let us summarize what has been done in this chapter. We started with a surface
$(M,\sigma)$, with or without boundary. We choosed on $M$ disjoint simple
loops $L_1,\ldots,L_q$, whose image is either a boundary component of $M$ or
contained in the interior of $M$. We picked $q$ elements $x_1,\ldots,x_q$ in
$G$. We almost proved the following theorem:

\begin{theorem} \label{main} There exists a probabilty space
  $(\Omega,\a,P(x_1,\ldots,x_q))$ and a family of $G$-valued random
  variables $(H_c)_{c\in PM}$ on this space, such that:\\
  \indent 1. For any graph $\Gamma=\{a_1,\ldots,a_r\}$ on $M$ such that
  $L_1,\ldots,L_q \in\Gamma^*$, the law of $(H_{a_1},\ldots,H_{a_r})$ is
  the discrete Yang-Mills measure $P^\Gamma(x_1,\ldots,x_q)$ on $G^\Gamma$.\\
  \indent 2. For any path $c$ of $PM$ and any sequence $(c_n)_{n\geq 0}$ of
  paths of $PM$ such that $c_n \build{\lra}_{n\to\infty}^{d_1} c$ with fixed
  end points, we have $H_{c_n} \build{\lra}_{n\to\infty}^{d_P}
  H_c$.\\
  The law of this family of random variables is uniquely defined by
  these two properties. Moreover, it has the following properties:\\
  \indent 3. If $c_1$ and $c_2$ are paths that can be concatenated to form
  $c_1
  c_2$, then $H_{c_1 c_2}=H_{c_2} H_{c_1}$ a.s. \\
  \indent 4. If $\varphi:M\lra M$ is a diffeomorphism such that $\varphi_*
  \sigma=\sigma$, then $\varphi$ induces a permutation of the set of paths
  $PM$ and the families $(H_c)_{c\in PM}$ and $(H_{\varphi(c)})_{c\in PM}$
  have the same law.
\end{theorem}

\pf. We already proved the existence of the family. When $M$ has a boundary,
the probability space is that associated with a minimal closure of $M$.  Let
us prove the uniqueness in law. This law is a probability measure on the set
$\f(PM,G)$ endowed with the $\sigma$-algebra generated by cylinder sets. So it
is characterized by its finite-dimensional marginals. Since any family of
paths can be approximated by families of paths in graphs, for example
piecewise geodesic paths for some metric, the law of the random holonomy along
an arbitrary finite family of paths is determined by properties (1) and (2).

Property (3) was already proved in proposition \ref{multiplicative} for closed
surfaces. For surfaces with boundary, the construction by restriction of the
random holonomy on a minimal closure obviously preserves the multiplicativity.

Property (4) was proved at the discrete level in proposition \ref{inv_by_apd}.
Since the law of the whole family is determined by discrete laws, it is also
true is the continuous setting. \qed

Given $(M,\sigma)$, $L_1,\ldots,L_q$, $x_1,\ldots,x_q$, the law whose
existence and uniqueness is stated by this theorem is a measure on
$(\f(LM,G),\c)$, where $\c$ is the $\sigma$-algebra generated by the cylinder
sets. We shall denote this measure by $\mu_0(x_1,\ldots,x_q)$, or just $\mu_0$
if $q=0$. We keep the notation $(H_c)_{c\in PM}$ for the canonical process on
the space $(\f(LM,G),\c)$.

\subsection{Disintegration formula}

Consider a surface $(M,\sigma)$. Recall from proposition
\ref{discrete_disintegration} that the conditional discrete measures
constitute a disintegration of the free discrete Yang-Mills measure. We want
to extend this result to the continuous setting. As usual, $L_1,\ldots,L_q$
are loops on $M$.

\begin{proposition} \label{disint_mu0}
  The map $(x_1,\ldots,x_q)\mapsto \mu_0(x_1,\ldots,x_q)$ provides a
  disintegration of the measure $\mu_0$ on $(\f(PM,G),\c)$ with respect to the
  random variable $(H_{L_1},\ldots,H_{L_q})$.
\end{proposition}

\pf. By construction, $(H_{L_1},\ldots,H_{L_q})=(x_1,\ldots,x_q)$
$\mu_0(x_1,\ldots,x_q)$-a.s. Let $(c_1,\ldots,c_n)$ be a family of paths of
$PM$. We need to prove that, for any function $f$ continuous on $G^q$,
$$E_{\mu_0} f(H_{c_1},\ldots,H_{c_q}) = \int_{G^q} E_{\mu_0(x_1,\ldots,x_q)}
f(H_{c_1},\ldots,H_{c_q}) \; d\eta(x_1,\ldots,x_q),$$
where $\eta$ is the law
of $(H_{L_1},\ldots,H_{L_q})$ under $\mu_0$.  We already know that this result
is true if $c_1,\ldots,c_n$ are paths in a graph. If they are not, we can
approximate them in the $d_1$-topology by paths in graphs so that both
expectations appearing in the formula converge. Since $G$ is compact, $f$ is
bounded and the dominated convergence theorem applies. \qed

\section{Yang-Mills measure}
\label{YM_measure}

\subsection{Definition of the Yang-Mills measure}

In this paragraph, we will explain why and how the measure $\mu_0$ defined in
the preceding section
still has to be transformed in order to become something that
might be called Yang-Mills measure. 

According to the formal description, Yang-Mills measure should be a measure on
the quotient space $\a/\j$ of connections modulo gauge transformations. But an
element of this space does not determine a holonomy along each path on $M$
that could be intrinsically represented by an element of $G$. Indeed, the
holonomy along an open path $c$, i.e. such that $c(0)\neq c(1)$, can be
transformed into any other $G$-equivariant diffeomorphism of the fiber over
$c(0)$ into the fiber over $c(1)$ by an appropriate gauge transformation. The
fact that the law of the random holonomy along an edge is always uniform on
$G$ could be thought of as a reflect of this geometric property. This is why
we will restrict to the set $LM$ of loops on $M$ instead of $PM$. Thus we will
consider the family $(H_l)_{l\in LM}$ whose law is a probability measure on
$(\f(LM,G),\c)$, where we keep the notation $\c$ for the $\sigma$-algebra
generated by the cylinders.

But it is still not true that an element of $\a/\j$ determines an element of
$G$ as holonomy along each loop. Gauge transformations act by conjugation on
the holonomy along loops. More precisely, they conjugate in the same way the
holonomies along loops based at the same point.  Let us denote by $\Ad$ the
diagonal action of $G$ on $G^n$ defined by:
$$\Ad(g)(g_1,\ldots,g_n)=(\Ad(g)g_1,\ldots,\Ad(g)g_n).$$
Orbits of this action
will be called joint conjugacy classes and the joint class of
$(g_1,\ldots,g_n)$ will be denoted by $[g_1,\ldots,g_n]$. We can reformulate
our observation by saying that an element of $\a/\j$ determines the the joint
conjugacy class of the holonomy along all loops based at the same point.
Sengupta proved the converse of this statement (prop.  $2.1.2$ in
\cite{Sengupta_S2}):

\begin{proposition} [\cite{Sengupta_S2}] Let $\omega_1$ and $\omega_2$ be two connections
  on $M$. Let $m_0$ be a point on $M$. Suppose that along any finite family of
  loops $l_1,\ldots,l_n$ based at $m_0$, the joint conjugacy classes of the
  holonomies defined by $\omega_1$ and $\omega_2$ are equal. Then $\omega_1$
  and $\omega_2$ belong to the same class in $\a/\j$.
\end{proposition}

Let the group $\f(M,G)$ act on $\f(LM,G)$ in the following way: if $j \in
\f(M,G)$, $f\in \f(LM,G)$ and $l\in LM$, set
$$(j\cdot f)(l)=j(l(0))^{-1} f(l) j(l(0)).$$
This action extends the action of
a discrete gauge transformation.  We can summarize our observations as
follows:

\begin{proposition} The holonomy allows to define an injective map 
$$\a/\j \lra \f(LM,G)/\f(M,G).$$
\end{proposition}

This result says that the quotient space $\f(LM,G)/\f(M,G)$ can be viewed as
an extension of the space of connections modulo gauge transformations. We want
to define the Yang-Mills measure on this space. To begin with, we must define
a convenient $\sigma$-algebra.

There is a set of natural functions on the quotient space: given $n$ loops
$l_1,\ldots,l_n$ based at the same point, the joint class
$[H_{l_1},\ldots,H_{l_n}]$ is a well-defined function that we denote by
$\h_{l_1,\ldots,l_n}$. We will consider the $\sigma$-algebra $\a$ generated by
the set of these functions. Of course, we want to be able to consider random
variables associated with families of loops that are not based at the same
point. We claim that the $\sigma$-algebra $\a$ allows to do this. Indeed, let
$l_1,\ldots,l_n$ be a family of loops that we rewrite
$(l_1,\ldots,l_{i_1}),\ldots,(l_{i_p+1},\ldots,l_n)$, putting together the
loops based at the same points. Then we can define the variable
$\h_{l_1,\ldots,l_n}$ by
$$\h_{l_1,\ldots,l_n}=(\h_{l_1,\ldots,l_{i_1}},\ldots,\h_{l_{i_p+1},\ldots,l_n})$$
and this random variable is measurable with respect to $\a$.  Remark that $\a$
may also be seen as a $\sigma$-algebra on $\f(LM,G)$, invariant by the action
of $\f(M,G)$, since the functions $\h_{l_1,\ldots,l_q}$ are also naturally
defined on this space. Another natural choice for $\a$ would have been to
consider the $\f(M,G)$-invariant sets of the cylinder $\sigma$-algebra $\c$.
We shall discuss this point at the end of this section.

\begin{proposition}\label{restriction_de_mu} Let $(M,\sigma)$ be a
  surface. Let $L_1,\ldots,L_q$ be disjoint simple loops on $M$ whose image is
  either a component of the boundary of $M$ or contained in the interior of
  $M$. Let $(x_1,\ldots,x_q)$ be an element of $G^q$. The restriction of
  $\mu_0(x_1,\ldots,x_q)$ to $\a$ depends on each $x_i$ only through its
  conjugacy class.
\end{proposition}

\pf. The point is to understand how $\mu_0(x_1,\ldots,x_q)$ is transformed
under the action of $\f(M,G)$. Similarly to what we proved in \ref{gt_of_P},
if $j$ is an element of $\f(M,G)$, then, setting $y_i=j(L_i(0))$, we have
$$j_*P(x_1,\ldots,x_q)=P(y_1^{-1} x_1 y_1,\ldots,y_q^{-1} x_q y_q).$$
Indeed,
we already know that this equality holds when we evaluate these measures
against functions of the holonomy along paths in a graph, and we extend it to
general measurable functions by the usual approximation scheme.

Thus, the $\mu_0(x_1,\ldots,x_q)$-measure of sets invariant under the action
of $\f(M,G)$ depends only on the conjugacy classes $[x_1],\ldots,[x_q]$. \qed

We denote by $t_i$ the conjugacy class of each $x_i$.

\begin{definition} \label{maindef}
  We call Yang-Mills measure on $M$ and denote by $\mu$ the image measure of
  $\mu_0$ on the quotient space $(\f(LM,G)/\f(M,G),\a)$, or equivalently the
  restriction of $\mu_0$ to $(\f(LM,G),\a)$.
  
  Similarly, we call conditional Yang-Mills measure with respect to
  $L_1,\ldots,L_q$ and we denote by $\mu(t_1,\ldots,t_q)$ the image measure of
  $\mu_0(x_1,\ldots,x_q)$ on the quotient space $(\f(LM,G)/\f(M,G),\a)$, or
  equivalently the restriction of $\mu_0(x_1,\ldots,x_q)$ to $(\f(LM,G),\a)$.
\end{definition}

The first point of view keeps track of the quotient structure of the space
$\a/\j$. Nevertheless, the second will often be technically more convenient.

\begin{proposition} \label{disint}
  The map $(t_1,\ldots,t_q) \mapsto \mu(t_1,\ldots,t_q)$ defined on
  $(G/\Ad)^q$ provides a disintegration of the measure $\mu$ with respect to
  the random variable $\h_{L_1},\ldots,\h_{L_q}$.
\end{proposition}

Note that since the $L_i$'s are not based at the same point, the variables
$\h_{L_1,\ldots,L_q}$ and $(\h_{L_1},\ldots,\h_{L_q})$ are equal. \\

\pf. Let $(t_1,\ldots,t_q)$ be an element of $(G/\Ad)^q$ and
$(x_1,\ldots,x_q)\in G^n$ be such that $[x_i]=t_i$ for each $i$. Then
$([H_{L_1}],\ldots,[H_{L_q}])=([x_1],\ldots,[x_q])$
$\mu_0(x_1,\ldots,x_q)$-a.s., so that
$$(\h_{L_1},\ldots,\h_{L_q})=(t_1,\ldots,t_q) \;\; \mu(t_1,\ldots,t_q)-{\rm
  a.s.}$$
By \ref{partition_law}, we know that the law of
$(H_{L_1},\ldots,H_{L_q})$ under $\mu_0$ is $Z^{-1} Z(x_1,\ldots,x_q) \; dx_1
\ldots dx_q$. We also proved in \ref{disint_mu0} that $\mu_0$ is disintegrated
by the $\mu_0(x_1,\ldots,x_q)$, so that
$$\mu_0={1\over Z} \int_{G^q} Z(x_1,\ldots,x_q) \mu_0(x_1,\ldots,x_q) \; dx_1
\ldots dx_q.$$
If we evaluate these measures on sets of $\a$ and use the
invariance by conjugation of the conditional partition function stated in
\ref{inv_part}, we find:
\begin{eqnarray*}
\mu_0 &=& {1\over Z} \int_{G^q} Z([x_1],\ldots,[x_q]) \mu_0([x_1],\ldots,[x_q])
\; dx_1 \ldots dx_q \\
&=& {1\over Z} \int_{(G/\Ad)^q} Z(t_1,\ldots,t_q)
\mu_0(t_1,\ldots,t_q) \; dt_1 \ldots dt_q,
\end{eqnarray*}
where $dt$ is the image measure on $G/\Ad$ of the Haar measure. The last
equality restricted to $\a$-measurable sets proves the result. \qed

\subsection{Regularity properties}

In order to study Yang-Mills measure, we need to say more about the set of
joint conjugacy classes $G^n/\Ad$. We regard it as a set of compact subsets of
$G^n$ and endow it with the Hausdorff distance, defined in general between two
compact sets by
$$d(K_1,K_2)=\sup(\sup_{k_1 \in K_1} \inf_{k_2 \in K_2} d(k_1,k_2),\sup_{k_2
  \in K_2} \inf_{k_1\in K_1} d(k_1,k_2)).$$

\begin{lemma} The canonical projection $G^n \lra G^n/\Ad$ is
  $1$-Lipchitz.
\end{lemma}

\pf. Let $(g_1,\ldots,g_n)$ and $(h_1,\ldots,h_n)$ be two elements of $G^n$.
\begin{eqnarray*}
d([g_1,\ldots,g_n],[h_1,\ldots,h_n]) &=& \sup_{g\in G} \inf_{h\in G}
d(\Ad(g)(g_1,\ldots g_n),\Ad(h)(h_1,\ldots,h_n)) \\
&=& \inf_{h\in G} d((g_1,\ldots,g_n),\Ad(h)(h_1,\ldots,h_n)) \\
&\leq & d((g_1,\ldots,g_n),(h_1,\ldots,h_n)).
\end{eqnarray*}
\vskip -.6cm \qed

As before, this distance on $G^n/\Ad$ allows to define the distance
$$d_P(X,Y)=Ed(X,Y)$$
between $G^n/\Ad$-valued random variables defined on the
same probability space. The regularity property \ref{continuity_holonomy} of
the random holonomy becomes the following regularity property for the
Yang-Mills measure:

\begin{proposition} \label{reg_YM}
  Let $((l_{1,k},\ldots,l_{n,k}))_{k\geq 0}$ be a sequence of $n$-uples of
  loops such that\\
\indent  1. for each $k\geq 0$, the loops $l_{1,k},\ldots,l_{n,k}$ are based
  at the same point,\\
\indent  2. for each $i=1,\ldots,n$, there exists a loop $l_i$ such that
  $l_{i,k} \build{\lra}_{k\to\infty}^{d_1} l_i$. \\
  Then
  $$\h_{l_{1,k},\ldots,l_{n,k}} \build{\lra}_{k\to\infty}^{d_P}
  \h_{l_1,\ldots,l_n}.$$
\end{proposition}

\pf. The loops $l_i$ are necessarily based at the same point, denoted by $m$.
Denoting by $m_k$ the base point of the $l_{i,k}$'s, we have $m_k \lra m$. For
each $k$, let $z_k$ denote an arbitrary path joining $m$ to $m_k$.  Then for
each $i$ $z_k l_{i,k} z_k^{-1} \build{\lra}_{}^{d_1} l_i$ with fixed
basepoint, so that
$$H_{z_k l_{i,k} z_k^{-1}} \build{\lra}_{}^{d_P} H_{l_i}.$$
Since the
projection on $G^n/\Ad$ reduces the distances, this implies
$$[H_{z_k l_{1,k} z_k^{-1}},\ldots,H_{z_k l_{n,k} z_k^{-1}}]
\build{\lra}_{k\to\infty}^{d_P} [H_{l_1},\ldots,H_{l_n}].$$
The left hand side
term is equal to
$$[H_{z_k}^{-1} H_{l_{1,k}} H_{z_k},\ldots,H_{z_k}^{-1} H_{l_{n,k}} H_{z_k}]=
[H_{l_{1,k}},\ldots,H_{l_{n,k}}],$$
so that the result is proved. \qed

\subsection{Remarkable subfamilies of random variables}
\label{remarkable}

We study two special subfamilies of random variables defined on
$(\f(LM,G)/\f(M,G),\a,\mu)$, using the results proved in the preceding
paragraphs.

We begin by the family $(\h_l)_{l\in LM}$. Each variable is
$G/\Ad$-valued and this family satisfies a very nice regularity
property:

\begin{proposition} \label{cont_5}
  Let $l$ be in $LM$ and $(l_n)_{n\geq 0}$ be a sequence of loops that
  converges to $l$. Then $\h_{l_n} \build{\lra}_{}^{d_P} \h_l$.
\end{proposition}

This is the only situation where we can forget about end points.
Unfortunately, this family does not generate $\a$, since it does not contain
any information about joint conjugacy classes.

Now fix a point $m\in M$ and consider the set $L_mM$ of loops based at $m$. We
are interested in the family $(\h_{l_1,\ldots,l_n})_{l_i \in L_mM}$. It has
the same property as that stated in \ref{reg_YM}, the condition on end points
being always satisfied. What is interesting here is the following fact:

\begin{proposition} \label{based_loops_gen_t}
  The family $(\h_{l_1,\ldots,l_n})_{l_i \in L_mM}$ generates the
  $\sigma$-algebra $\a$.
\end{proposition}

\pf. Let $l_1,\ldots,l_n$ be $n$ loops on $M$ based at a point
$m_1$. Let $c$ be a path joining $m$ to $m_1$. Then the equality
$$\h_{cl_1c^{-1},\ldots,cl_nc^{-1}}=\h_{l_1,\ldots,l_n}$$
proves that it is
always possible to get back to loops based at $m$. \qed

This subfamily satisfies also a multiplicativity property. Indeed, the joint
conjugacy class of some elements of $G$ determines the joint class of all
products of these elements. For example, there is a well defined map from
$G^n/\Ad$ to $G^{n-1}/\Ad$ that sends $[g_1,\ldots,g_{n-1},g_n]$ to
$[g_1,\ldots,g_{n-2},g_ng_{n-1}]$. The multiplicativity can be expressed by
saying that for any $l_1,\ldots,l_n \in L_mM$, this maps sends almost surely
$\h_{l_1,\ldots,l_n}$ to $\h_{l_1,\ldots,l_{n-1}l_n}$.

\subsection*{Remark.} Let us discuss the definition of the $\sigma$-algebra
$\a$. For this, we consider the Yang-Mills measure an invariant measure on the
space $\f(L_mM,G)$, because in this setting, the action of the gauge group is
that of the finite dimensional group $G$. We use this fact below in order to
integrate functions over the orbits of this action.

We could have made another natural choice of an invariant $\sigma$-algebra on
the space $\f(L_mM,G)$, namely that of invariant sets of the cylinder
$\sigma$-algebra $\c$ on $\f(L_mM,G)$. Let us denote by $\c_\i$ this
$\sigma$-algebra. It is clear that $\a \subset \c_\i$ and it is very likely
that this inclusion is in fact an equality. We prove that the completions
$\widetilde\a$ and $\widetilde{\c_\i}$ with respect to $\mu_M(x)$ are equal.

We use the separability of $L_mM$ proved in the next lemma.

\begin{lemma} \label{separable}
  Let $M$ be a surface and $m$ a point of $M$. The loop space $L_mM$ endowed
  with the $d_1$-topology is separable.
\end{lemma}

\pf. We construct a countable dense subset of $L_mM$. The first point is that
$M$ itself is separable. Choose a countable dense subset $\Pi\subset M$
containing $m$. Endow $M$ with a Riemannian metric such that $\partial M$ is
geodesic if it is non empty. Let $R_M$ be the convexity radius of $M$: two
points at distance smaller than $R_M$ are joined by a unique minimizing
geodesic. Define $\Lambda$ to be the set of loops obtained by concatenation of
a finite number of geodesic segments joining two points of $\Pi$ at distance
smaller than $R_M$. The set $\Lambda$ is countable because it is equipotent to
a subset of finite sequences of $\Pi$. We claim that it is dense in $L_mM$.
Indeed, any geodesic segment of length smaller than $R_M$ can be approximated
by segments joining points of $\Pi$, since a small piece of geodesic depends
continuously on its end points.  Thus, the $d_1$-closure of $\Lambda$ contains
the set of piecewise geodesic loops and we already know that this set is dense
in $L_mM$.\qed

Let $\Lambda \subset L_mM$ be a countable dense subset. Let $\a^\Lambda$
denote $\sigma(\h_{\lambda_1,\ldots,\lambda_n} \in \Lambda)$. Let $\c^\Lambda$
denote $\sigma(H_\lambda, \lambda \in \Lambda)$ and $(\c^\Lambda)_\i$ denote
the invariant sets of this $\sigma$-algebra. It is clear that $\a^\Lambda
\subset (\c^\Lambda)_\i$. It is also clear that $\a^\Lambda$ is a separable
$\sigma$-algebra. Finally, one easily checks that any atom of $\a^\Lambda$ is
contained in an atom of $(\c^\Lambda)_\i$. Thus, Blackwell's theorem implies
(see \cite{Dellacherie_Meyer}) that $\a^\Lambda = (\c^\Lambda)_\i$.

We use this equality to prove the inclusion $\c_\i \subset \widetilde\a$,
which implies the result. The point is that the continuity in probability of
the map $l\mapsto \h_l$ and the density of $\Lambda$ in $L_mM$ imply that
$\widetilde{\c^\Lambda}$ contains $\c$ and that $\widetilde{\a^\Lambda}$
contains $\a$. This last inclusion implies the equality
$\widetilde{\a^\Lambda}=\widetilde\a$. Thus, it is sufficient to prove that
$\c_\i \subset \widetilde{\a^\Lambda}$. Let $f$ be a $\c_\i$-measurable
function. As a $\c$-measurable function, it is $\mu_M(x)$-almost surely the
limit of a sequence $(f_n)$ of $\c^\Lambda$-measurable functions. Let us
integrate these functions $f_n$ over the orbits of the action of $G$ on
$\f(L_mM,G)$, using the Haar measure on $G$. We get a sequence of
$(\c^\Lambda)_\i$-measurable functions still converging to $f$. Thus, $f$ is
measurable with respect to the completion
$\widetilde{(\c^\Lambda)_\i} =\widetilde{\a^\Lambda}$ and we get the result.


\pagevide
\chapter{Abelian theory} 
\label{abelian_theory}

In this chapter, we continue the investigation of the case $G=U(1)$ started in
section \ref{Abelian_theory_1}. Recall that we had reconstructed the random
holonomy along loops homologous to zero in a graph, using a white noise on $M$
(see proposition \ref{wn_rep_loops}). We extend now this reconstruction to all
cycles of $M$, using the unicity properties of the Yang-Mills measure proved
in chapter \ref{Continuous_theory}.

Then we show that it is possible to proceed backwards, namely to extract a
white noise on $M$ from the Yang-Mills measure on $M$, more precisely, using
the random holonomy along very small loops. This makes clear the relationship
between the random holonomy and the white noise in this Abelian case.

\section{The random holonomy as a white noise functional}
\label{rh_as_wn_functional}

As usual, $M$ may have a boundary $\partial M=N_1 \cup \ldots \cup N_p$.
Choose elements $x_1,\ldots,x_p$ in $U(1)$ and set $x=x_1 \ldots x_p$ or $x=1$
if $M$ has no boundary. We denote by $CM$ the set of cycles on $M$, i.e. the
set of linear combination of loops with integer coefficients and by
$C_0\Gamma$ the set of cycles homologous to zero. The family of random
variables $(H_l)_{l\in LM}$, that we used in \ref{YM_measure} to construct the
Yang-Mills measure $\mu_M(x_1,\ldots,x_p)$, extends by multiplicativity to the
cycles of $CM$ and gives rise to a measure on $\f(CM,U(1))$. In this Abelian
setting, the action of $\f(M,U(1))$ is trivial.

We seek a result similar to \ref{wn_representation}, valid for all cycles on
$M$. We begin by defining a family of random variables using a white noise on
$M$ and prove later that it has the law of a Yang-Mills random holonomy.

Recall that we proved earlier that the holonomies along a system of loops
representing a basis of $H_1(M;\Z)$ are independent uniform variables on
$U(1)$, independent of the holonomies along loops homologous to zero (see
proposition \ref{law_1}).

\subsection{Definition of the white noise functional}

Recall that $g$ denotes the genus of $M$. In order to define the double layer
potential (see \ref{def_dlp}), we need a Riemannian metric on $M$, that we
choose such that the boundary of $M$ is geodesic. Let $\ell_1,\ldots,\ell_{2g}$ be
piecewise geodesic loops such that
$$\b=([\ell_1],\ldots,[\ell_{2g}],[N_1],\ldots,[N_{p-1}])$$
is a basis of $H_1(M;\Z)$. 

Let $c$ be a cycle of $CM$. We can decompose it in $H_1(M;\Z)$:
$$[c]=\lambda_1 [\ell_1] + \ldots + \lambda_{2g} [\ell_{2g}] + \nu_1 [N_1] +
\ldots + \nu_{p-1} [N_{p-1}].$$
Let $U_1,\ldots,U_{2g}$ be $2g$ independent uniform random variables
on $U(1)$. Set
$$\Theta_c=U_1^{\lambda_1} \ldots U_{2g}^{\lambda_{2g}} x_1^{\nu_1}
\ldots x_{p-1}^{\nu_{p-1}}.$$
There is a cycle $c^\perp$ of $C_0M$, i.e. a cycle homologous to zero,
associated with $c$, defined by
$$c^\perp=c-(\lambda_1 \ell_1 + \ldots + \lambda_{2g} \ell_{2g} + \nu_1 N_1 +
\ldots + \nu_{p-1} N_{p-1}).$$
Recall that $u_{c^\perp}$ denotes the doule
layer potential of $c^\perp$.  The cycle $c^\perp$ is the boundary of a
$2$-chain $\alpha$. We defined an element $\sigmaint(c^\perp)$ of $\R$ (resp.
\R/\Z) when $M$ has a boundary (resp. no boundary) by
$\sigmaint(c^\perp)={{|\langle \sigma,\alpha \rangle|}\over {\sigma(M)}}$.

Let $W$ be a white noise on $M$, independent of the $U_i$'s. Let $T$ be a
variable independent of $W$ and the $U_i$'s, whose law is that described in
proposition \ref{gaussian_representation}. Finally, denote by $W_0$ the
projection of $W$ on the hyperplane of zero-mean functions: for any function
$u\in L^2(M,\sigma)$,
$$W_0(u)=W\left( u- {1\over {\sigma(M)}} \int_M u\; d\sigma \right).$$

We are able to define what will be proved to be a second realization
of the random holonomy along cycles. Denote by $(\Omega,P)$ a probability
space that supports $W$, $T$ and the $U_i$'s. 

\begin{definition} \label{def_u1}
For each cycle $c\in CM$, define the following random variable on
$(\Omega,P)$: 
$$\wh_c=\exp i(W_0(u_{c^\perp})+\sigmaint(c^\perp) T) \; \Theta_c.$$
\end{definition}

\subsection{Regularity of the new random holonomy}

In order to prove that this family has the law of the random holonomy,
we will check that this is true for a restricted class of paths,
namely piecewise geodesic paths, and extend this partial result by
continuity. This is why we are interested in the regularity of this
new family. 

We begin by extending the distance $d_1$ to the space $CM$ of cycles. 

\begin{definition} Let $c=n_1 l_1 +\ldots + n_k l_k$ and $c'=n'_1 l'_1
+\ldots + n'_{k'} l'_{k'}$ be two cycles on $M$, written as
combinations of loops. If $k\neq k'$, set $d_1(c,c')=1$. If $k=k'$,
let $\pi(c,c')$ be the set of permutations $\tau \in S_k$ such that
$n'_{\tau(i)}=n_i$ for all $i=1,\ldots,k$. If $\pi(c,c')=\emptyset$,
set $d_1(c,c')=1$. Otherwise, set
$$d_1(c,c')=\inf_{\tau \in \pi(c,c')} \sum_{i=1}^k n_i
d_1(l_i,l'_{\tau(i)}).$$
\end{definition}

We will use proposition \ref{continuity_loops_dlp} about the continuity of
the double layer potential of loops to prove the next
proposition. Recall that \ref{continuity_loops_dlp} was proved only on
surfaces without boundary.

\begin{proposition} \label{continuity_cycles_dlp}
Suppose that $M$ has no boundary. Let $c$ be
a cycle and $(c_n)_{n\geq 0}$ be a sequence of cycles such that $c_n
\build{\lra}_{}^{d_1} c$. Then $u_{c_n}
\build{\lra}_{}^{L^2} u_c$.
\end{proposition}

\pf: Decompose $c$ as $n_1 l_1 +\ldots +n_k l_k$. Let $N=n_1+\ldots +
n_k$. Fix $\epsilon>0$. For each $l_i$, there exists $\delta_i<1$ such that
$d_1(l'_i,l_i)<\delta_i$ implies $\parallel u_{l'_i} - u_{l_i}
\parallel_{L^2} < {\epsilon \over {kN}}$ for any loop $l'_i$. Suppose that
$n$ is such that $d_1(c_n,c)<\inf \delta_i$. Then $c_n$ can be written
$n_1 l'_1+ \ldots + n_k l'_k$, with $\sum n_i d_1(l_i,l'_i)<\inf
\delta_i$. For each $i$, we have in particular
$d_1(l_i,l'_i)<\delta_i$. Thus, 
$$\hskip 3cm \parallel u_{c_n} - u_c \parallel_{L^2} \leq \sum_{i=1}^k n_i
\parallel u_{l'_i} - u_{l_i} \parallel_{L^2} < {{kN\epsilon}\over {kN}} =
\epsilon. \hskip 2cm \qed$$

The homology class of a cycle $c$ depends continuously on $c$, even for the
distance $d_\infty$. Thus, $c^\perp$ and hence $u_{c^\perp}$ depends
continuously on $c$, and, by continuity of the white noise which is an
isometry, the map $c\mapsto W_0(u_{c^\perp})$ is continuous for the $L^2$ norm
when $M$ is closed. We want to extend this result to surfaces with boundary.
For this, we study $u_c$ when $c$ is homologous to zero.

Let $M$ be a surface, with or without boundary and fix $c\in C_0M$. Let $x$
and $y$ be two points of $M$ outside the image of $c$. If necessary, we modify
locally the $\ell_i$'s in a neighbourhood of $x$ in order to make sure that
$x$ meets none of them. Let $\ell_x$ be the boundary of a small disk $D_x$
around $x$, small enough not to meet $c$ and not to contain $y$. The module
$H_1(M-\{x,y\})$ is generated by
$[\ell_1],\ldots,[\ell_{2g}],[N_1],\ldots,[N_{p-1}],[\ell_x],[\ell_y]$. In
$M-\{x,y\}$, we have the equality
$$[c]=\sum_{i=1}^{2g} \lambda_i [\ell_i] + \sum_{j=1}^{p-1} \nu_j [N_j] + p
[\ell_x]+q[\ell_y]$$
for some $p\in \Z$. This equality also holds in $M$,
where $[c]=[\ell_x]=[\ell_y]=0$, and this proves that $\lambda_i=\nu_j=0$.
Thus, $[c]=p[\ell_x]+q[\ell_y] \;\; {\rm in} \;\; M-\{x,y\}.$

\begin{lemma} \label{dlp_h1}
  With the preceding notations, $u_c(x)-u_c(y)=p-q.$
\end{lemma}

\pf: Recall the definition of the Green function $G(\cdot,\cdot)$ from
equation (\ref{def_G}), and the value of the double layer potential of small
loops from proposition \ref{compute_dlp}. The $1$-form $*dG_x-*dG_y$ is closed
on $M-\{x,y\}$. So,
\begin{eqnarray*}
u_c(x)-u_c(y) &=& \int_c *dG_x-*dG_y \\
&=& p \int_{\ell_x} *dG_x-*dG_y + q \int_{\ell_y} *dG_x-*dG_y\\
&=& p(u_{\ell_x}(x) - u_{\ell_x}(y)) + q(u_{\ell_y}(x) - u_{\ell_y}(y)) \\
&=& p\left(\1_{D_x}(x) - {{\sigma(D_x)}\over {\sigma(M)}} -\1_{D_x}(y) +
{{\sigma(D_x)}\over {\sigma(M)}}\right) + q \left(\1_{D_y}(x)  -\1_{D_y}(y)
\right) \\ 
&=& p-q.
\end{eqnarray*}
\vskip -.6cm \qed

\subsection*{Remark.} If $M$ is closed, then $[\ell_x]+[\ell_y]=0$ in
$H_1(M-\{x,y\})$. In this case, $p$ and $q$ are defined up to an additive
constant but the difference $p-q$ is well defined.

\begin{corollary} \label{dlp_bnd}
\indent  1. When $M$ has no boundary, the double layer potential of cycles of
  $C_0M$ does not depend on the choice of the metric. \\
\indent  2. The potential $u_c$ is constant on each connected component of the
  complementary of the image of $c$. \\
\indent  3. If $M_1$ is a minimal closure of $M$ and if we identify $M$ with a
  submanifold of $M_1$, then for any $c\in C_0M$, the potentials $u_c^M$ and
  ${u_c^{M_1}}_{|M}$ computed respectively in $M$ and $M_1$ differ only by an
  additive constant.
\end{corollary}

\pf. $1.$ The lemma determines the potential of any cycle of $C_0M$ up to a
constant. When $M$ has no boundary this constant is determined by the
condition $\int_M u_c \; d\sigma=0$.

$2.$ Given a fixed point $y_0$, $u_c(x)-u_c(y_0)$ depends only on the homology
class of $c$ in $M-\{x,y\}$, which does not change if $x$ stays in a given
connected component.

$3.$ Both functions $u_c^M$ and ${u_c^{M_1}}_{|M}$ satisfy the property shown
in the lemma, with the same values of $p$ and $q$. Indeed, if
$[c]=p[\ell_x]+q[\ell_y]$ in $H_1(M)$, then the same equality holds in
$H_1(M-\{x,y\})$. Thus, they cannot differ by more than an additive
constant. \qed

\begin{proposition} Even if $M$ has a boundary, the map $c \mapsto
    W_0(u_{c^\perp})$ is continuous.
\end{proposition}

\pf: By property (3) of the preceding corollary, 
$$u_{c^\perp}^M - {1\over {\sigma(M)}} \int_M u_{c^\perp}^M \; d\sigma =
{u_{c^\perp}^{M_1}}_{|M} - {1\over {\sigma(M)}} \int_M u_{c^\perp}^{M_1} \; d\sigma.$$
Together with \ref{continuity_cycles_dlp}, this shows that $u_{c^\perp}$
depends continuously on $c_0$. This was the only missing point. \qed

Now we study the term $\sigmaint(c)$, when $c\in C_0M$. We will show that
it can be extracted from the double layer potential of $c$.

\begin{lemma} Let $c$ be a cycle of $C_0M$.\\
\indent 1. If $M$ has no boundary, $\sigmaint(c)$ is the element $t \in \R/\Z$
    such that $u_c$ takes its values in $\Z-t$. \\
\indent 2. If $M$ has a boundary, consider $M_1$ a minimal closure of
    $M$. Then $\sigmaint(c)$ is equal to $-{{\sigma(M_1)}\over
    {\sigma(M)}}$ times the value of $u_c^{M_1}$ at any point of
    $M_1-M$.
\end{lemma}

\pf: $1.$ Let $x$ be a point of $M$ outside the image of $c$. Denote by
    $\alpha$ a $2$-chain such that $c=\partial \alpha$. Then
\begin{eqnarray*}
u_c(x) = \int_c *dG_x &=&  \int_\alpha \delta_x -{{d\sigma}\over
{\sigma(M)}} \\
&\equiv& -{{|\langle\sigma,\alpha\rangle|}\over {\sigma(M)}} \;\; (\mod 1).
\end{eqnarray*}

$2.$ Let $x$ be a point of $M_1 - M$. We have
$$u_c^{M_1}(x)= \int_\alpha \delta_x - {{d\sigma}\over {\sigma(M_1)}} =
-{{|\langle\sigma,\alpha\rangle|}\over {\sigma(M_1)}} = - {{\sigma(M)}\over
    {\sigma(M_1)}} \sigmaint(c).$$
\vskip -.6cm \qed

Let $(c_n)_{n\geq 0}$ be a sequence of cycles of $C_0M$ such that $c_n
\build{\lra}_{}^{} c$. Since $u_{c_n} \build{\lra}_{}^{L^2}
u_c$ and since all these functions are locally constant, there is
pointwise convergence outside the image of $c$. This implies:

\begin{proposition} The map $c \mapsto \sigmaint(c)$ defined on $C_0M$ is
continuous. 
\end{proposition}

Finally, the map $c\mapsto \Theta_c$ is locally constant on $C_0M$. We
proved:

\begin{proposition} \label{cont_wh}
The map $c\mapsto \wh_c$ is continuous from $(CM,d_1)$
into the space of square integrable random variables on $(\Omega,P)$.
\end{proposition}

The last property that we need is the multiplicativity:

\begin{proposition} For any cycles $c_1$ and $c_2$ in $CM$, 
$$\wh_{c_1+c_2}=\wh_{c_1} \wh_{c_2} \;\; P-{\rm a.s.}$$
\end{proposition}

\pf: This follows immediately from the following facts: $c^\perp$ depends
linearly on $c$, the double layer potential is additive, the map $c\mapsto
\sigma(c)$ is also additive and the map $c\mapsto \Theta_c$ is multiplicative.
\qed

\subsection{Identification of the random holonomies}

We are now able to prove the main theorem of this section:

\begin{theorem} \label{abelian_main}
The family of random variables $(\wh_c)_{c\in CM}$ has the same law as
the family $(\h_c)_{c\in CM}$ under $\mu_M(x_1,\ldots,x_p)$.
\end{theorem}

\pf: According to the unicity statement of the theorem \ref{main} and to the
regularity property \ref{cont_wh}, it is enough to prove the equality of the
laws for piecewise geodesic cycles. Let $\Gamma$ be a piecewise geodesic graph
such that $\ell_1,\ldots,\ell_{2g} \in \Gamma^*$.  Recall that these loops
generate the first homology group of a minimal closure of $M$. Denote by
$F_1,\ldots,F_n$ the faces of $\Gamma$. The arguments developped in
\ref{decomposition_of_cycles} explain why it is enough to prove the equality
of the laws for the fundamental system
$(\ell_1,\ldots,\ell_{2g},N_1,\ldots,N_p,\partial F_1,\ldots,\partial F_n)$.

On one hand, $\ell_1^\perp=\ldots=\ell_{2g}^\perp=0$, so that
$\wh_{\ell_i}=U_i$ for all $i=1,\ldots,2g$. On the other hand,
$N_1^\perp=\ldots=N_{p-1}^\perp=0$ and $N_p^\perp=N_1+\ldots+N_p$. Thus,
$\sigmaint(N_p^\perp)=1$ and $u_{N_p^\perp}$ is equal to zero so that
$W_0(u_{N_p^\perp})=0$. This implies $\wh_{N_j}=x_j$ for all $j=1,\ldots,p-1$
and also for $j=p$. Finally, we already proved in \ref{wn_representation} that
$(\wh_{\partial F_1},\ldots, \wh_{\partial F_n})$ and $(\h_{\partial
  F_1},\ldots,\h_{\partial F_n})$ have the same law.  This terminates the
proof, since we know that $(\h_{\partial F_1},\ldots,\h_{\partial F_n})$ and
$(\h_{\ell_1},\ldots,\h_{\ell_{2g}})$ are independent under the Yang-Mills
measure. \qed

\section{Small scale structure of the Yang-Mills field}
\label{small_scale}
\subsection{Extraction of a white noise}

In the first part of this chapter, we explained how the data of a white noise
on $M$ and a bit more alea allows to reconstruct the Yang-Mills measure. We
proceed now backwards: we try to extract a white noise from the Yang-Mills
measure on a surface. In some sense, this amounts to compute the curvature of
a Yang-Mills random connection.

As usual, $(M,\sigma)$ is given, as well as elements $x_1,\ldots,x_p$ of $G$,
associated with the components of $\partial M$ and $x=x_1\ldots x_p$ or $x=1$
if $M$ has no boundary. We denote by $(\Omega_M,\mu_M(x_1,\ldots,x_p))$ the
space of the Yang-Mills measure on $M$.

In order to study the measure at small scale, we construct on $M$ a sequence
of partitions in the following way. Let $(\Gamma_n)_{n\geq 0}$ be a sequence
of graphs on $M$ such that $\Gamma_n$ has exactly $n$ faces denoted by
$F_{j,n}$, $j=1,\ldots,n$. We assume that $\sigma(F_{j,n})={{\sigma(M)}\over
  n}$ and also that the diameter of the faces decreases uniformly to $0$, i.e.
that for any metric on $M$, $\sup_j \diam(F_{j,n}) \lra 0$.  We fix an
orientation of $M$ and assume that the boundaries of the $F_{j,n}$'s are
oriented with the usual convention. For each couple $(j,n)$ with $n\geq 0$,
$1\leq j \leq n$, we denote the random variable $\h_{\partial F_{j,n}}$
defined on $(\Omega_M,\mu_M(x_1,\ldots,x_p))$ by $\h_{j,n}$ and see it as a
$\C$-valued random variable, identifying $U(1)$ with $\{ z\in \C, |z|=1 \}$.

For each $n\geq 0$, let $E_n$ denote the space of functions on $M$ constant on
each face of $\Gamma_n$. Set $E_\infty=\cup_n E_n$. The assumption on the
diameter of the faces $F_{j,n}$ imply that any continuous function on $M$ can
be uniformly approximated by functions of $E_\infty$. Thus, $E_\infty$ is
dense in $L^1(M,\sigma)$ and $L^2(M,\sigma)$ with their respective usual
norms.

In order to define a kind of white noise, we will proceed as for the
construction of the standard Wiener integral. We define a linear
form $I_n$ on each $E_n$. Let $f_n$ be a function of $E_n$ and let
$f_{j,n}$ be its value on $F_{j,n}$. We set
$$I_n(f_n)={1\over i} \sum_{j=1}^{n} f_{j,n} (\h_{j,n} -1).$$

\begin{theorem} \label{construct_wn}
  Let $f$ be a square-integrable function on $M$ and $(f_n)_{n\geq 0}$ a
  sequence of functions converging to $f$ in $L^2$ norm and such that $f_n \in
  E_n$. Then the sequence $(I_n(f_n))_{n\geq 0}$ converges in
  $L^2(\Omega_M,\mu_M(x_1,\ldots,x_p))$ to a random variable $I(f)$ that does
  not depend on the choice of the sequence $(f_n)$. The law of this random
  variable can be described in the following way. Let $W^0_f$ be a centered
  gaussian random variable with variance $\parallel f \parallel^2_{L^2_0} =
  \parallel f - {1 \over \sigma(M)} \int_M f d\sigma \parallel_{L^2}^2$. Let
  $T$ be a $\n(0,\sigma(M))$ random variable conditioned to take its values
  in $\exp^{-1}(x)$, independent of $W^0_f$. Then, the following identity
  holds in distribution:
\begin{equation} \label{complete-law}  
I(f) \build{=}_{}^{\rm law} W^0_f + \left(T+ {i\over 2}\right) {1\over
    {\sigma(M)}}\int_M f\; d\sigma.
\end{equation}
\end{theorem}
This proves in particular that the law of $I(f)$ does not depend on the choice
of the orientation of $M$. \\

\pf: To prove this theorem, it is convenient to use the white noise
realization of the Yang-Mills measure. Let $(\Omega,P)$ be a probability space
on which a pair $(W,T)$ is defined, consisting in a white noise $W$ and a
random variable $T$ independent of $W$, whose law is that described in the
theorem. We do not need the variables $U_i$, because we are only computing the
holonomy along loops that are homologous to zero. Set
$$Y_{j,n}=W(\1_{F_{j,n}})\; , \;\; S_n=\sum_{j=1}^{n} Y_{j,n} \; , \;\;
X_{j,n}=Y_{j,n}-{1\over {n}} S_n.$$
We know by the theorem \ref{abelian_main} that the law
of the sequence $(I_n(f_n))$ can be represented on $(\Omega,P)$ by
$$\left( {1\over i} \sum_{j=1}^{n} f_{j,n} \left( e^{i (X_{j,n} + {T\over
        {n}})} - 1 \right) \right)_{n\geq 0}.$$
We will prove the theorem for
this sequence. For this, we study the following Lagrange inequality:
\begin{eqnarray}
&&\hskip -.5cm \left| \sum_{j=1}^{n} f_{j,n} \left( e^{i (X_{j,n} + {T\over
{n}})} - 1 \right) - i\sum_{j=1}^{n} f_{j,n} \left(X_{j,n} + {T\over
{n}}\right) +{1\over 2} \sum_{j=1}^{n} f_{j,n} \left(X_{j,n} +{T\over
{n}}\right)^2 \right | \leq \nonumber \\
&& \hskip 8 cm  \sum_{j=1}^{n} |f_{j,n}| \left|X_{j,n} + {T \over
{n}} \right|^3 \label{TAF}
\end{eqnarray}

We will often use of the following lemma:

\begin{lemma} For each positive integer $p$, there exists a constant
$C_p$ such that
$$E |X_{j,n}|^p \leq {{C_p}\over {n^{p\over 2}}}.$$
\end{lemma}

\pf: This is just a consequence of the fact that a centered gaussian random
variable $X$ of variance $t$ satisfies $E |X|^p = C_p t^{p\over 2}$ for come
constant $C_p$ independent of $t$ and that $X_{j,n}$ has variance
${{\sigma(M)}\over {n}} - {\sigma(M)\over {n^2}}$. \qed

We begin by showing that the right hand side term converges to zero in
$L^2(\Omega,P)$.y
$$ E \left| X_{j,n} + {T\over {n}} \right|^6 \leq  \sum_{k=0}^6
\left(\!\matrix{6\cr k}\!\right) E |X_{j,n}|^{6-k} E {{|T|^k}\over {n^k}} \leq  {K\over {n^3}},$$
so that 
\begin{eqnarray*}
E\left| \sum_{j=1}^{n} |f_{j,n}|\left|X_{j,n} +{T\over {n}}\right|^3
\right|^2 &=& \sum_{j,k=1}^{n} |f_{j,n}| |f_{k,n}| E \left|X_{j,n} +{T\over
{n}}\right|^3 \left|X_{k,n} +{T\over {n}}\right|^3 \\
&\leq & \left( \sup_{j=1}^{n} n \parallel (X_{j,n} +{T\over
        {n}})^3 
\parallel_{L^2} \right)^2 \left(\sum_{j=1}^{n} {{|f_{j,n}|}\over {n}}
\right)^2 \\ 
&\leq & {1\over {n}} \parallel f_n \parallel_{L^1}^2
\build{\lra}_{n\to\infty}^{} 0.
\end{eqnarray*}

Now look at the second order term of the left hand side of (\ref{TAF}). Let
$m_n={1\over {\sigma(M)}} \int_M f_n \; d\sigma$ denote the mean of $f_n$ and
$f_n^0=f_n-m_n$ denote its zero-mean part. We will use several times the fact
that $\sum_j f_{j,n}^0 =0$.  We have
\begin{equation}
\label{2d_order}
\sum_{j=1}^{n} f_{j,n} (X_{j,n} +{T\over{n}})^2 = \sum_{j=1}^{n}
f_{j,n}^0 (X_{j,n} +{T\over {n}})^2 + m_n \sum_{j=1}^{n}
(X_{j,n} +{T\over {n}})^2.
\end{equation}

Let us study the first term of this decomposition. In all estimations,
$C$ denotes a constant, i.e. a number that depends neither on $j$ nor
on $n$. It may denote different constants at different lines. 
\begin{equation}
\label{2d_order_1}
\sum_{j=1}^{n} f_{j,n}^0 (X_{j,n} +{T\over {n}})^2 = \sum_{j=1}^{n}
f_{j,n}^0 X_{j,n}^2 + 2 \sum_{j=1}^{n} {{f_{j,n}^0}\over {n}}
X_{j,n} T.
\end{equation}

The first term of the right hand side term can be written:
$$\sum_{j=1}^{n} f_{j,n}^0 X_{j,n}^2 = \sum_{j=1}^{n} f^0_{j,n}
(Y_{j,n}-{1\over {n}} S_n)^2 = \sum_{j=1}^{n} f_{j,n}^0 Y_{j,n}^2 - 2
\sum_{j=1}^{n} {{f_{j,n}^0}\over {n}} Y_{j,n} S_n. $$

On one hand, 
$$E\left| \sum_{j=1}^{n} f_{j,n}^0 Y_{j,n}^2 \right|^2 = E\left|
  \sum_{j=1}^{n} {{f_{j,n}^0}\over {n}} n Y_{j,n}^2 \right|^2 = \sum_{j=1}^{n}
{{|f_{j,n}^0|^2}\over {2^{2n}}} E|n Y^2_{j,n}|^2 \leq {{C}\over {n}} \parallel
f_n^0 \parallel_{L^2} \build{\lra}_{n\to\infty}^{} 0,$$
since $E|nY_{j,n}|^2$ depends neither on $j$ nor on $n$.

On the other hand, 
$$\parallel Y_{j,n} S_n \parallel_{L^2}^2 = E(Y_{j,n}^2 S_n^2) \leq n
(E Y_{j,n}^2)^2 + E Y_{j,n}^4 \leq {{C}\over {n}}$$
implies
$$\parallel \sum_{j=1}^{n} {{f_{j,n}^0}\over {n}} Y_{j,n} S_n \parallel_{L^2}
\leq \sum_{j=1}^{n} {{|f_{j,n}^0|}\over {n}} \parallel Y_{j,n} S_n
\parallel_{L^2} \leq {{C}\over {2^{{n\over 2}}}} \parallel f_n^0
\parallel_{L^1} \build{\lra}_{n\to\infty}^{} 0.$$

We proved that the first term of the r.h.s.  of (\ref{2d_order_1})
tends to $0$. To study the second one, note that
$$\parallel X_{j,n} T \parallel_{L^2}^2 = E X_{j,n}^2 E T^2 \leq {C
\over {n}},$$
so that
$$\parallel \sum_{j=1}^{n} {{f_{j,n}^0}\over {n}} X_{j,n} T \parallel_{L^2}
\leq \sum_{j=1}^{n} {{|f_{j,n}^0|}\over {n}} \parallel X_{j,n} T
\parallel_{L^2} \leq {C\over {\sqrt{n}}} \parallel f_{j,n}^0
\parallel_{L^1} \build{\lra}_{n\to\infty}^{} 0.$$

We proved that the zero-mean part of $f_n$ does not contribute to the
second order term. Let us study the last term of (\ref{2d_order}). 
$$m_n \sum_{j=1}^{n} \left(X_{j,n} +{T\over {n}}\right)^2 = m_n
\sum_{j=1}^{n} X_{j,n}^2 + {{m_n T}\over {2^{n-1}}} \sum_{j=1}^{n}
X_{j,n} + {{m_n}\over {n}} T^2.$$
We have $\sum_j X_{j,n}=0$ a.s. and ${{m_n}\over {n}} T^2 \lra 0$
a.s. . It remains 
$$ m_n \sum_{j=1}^{n} X_{j,n}^2 = m_n \sum_{j=1}^{n} Y_{j,n}^2 +
{{m_n}\over {n}} S_n^2 - {{m_n}\over {2^{n-1}}} S_n^2.$$
Since the law of $S_n^2$ does not depend on $n$, the two last terms
tend to zero. In order to determine the limit of the first one, we
compute
\begin{eqnarray*}
&& E\left| \sum_{j=1}^{n} Y_{j,n}^2 -{1\over {n}} \right|^2 = \sum_{j=1}^{n}
E\left|Y_{j,n}^2 -{1\over {n}}\right|^2 = \sum_{j=1}^{n} E\left(Y_{j,n}^4 +
  {1\over {n^2}} - {{2Y_{j,n}^2}\over {n}}\right)\\
&& \hskip 5.8cm  = \sum_{j=1}^{n} {C \over
  {n^2}} +{1\over {n^2}} - {2\over {n^2}} \leq {C\over {n^2}}.
\end{eqnarray*}
Thus,
$$m_n \sum_{j=1}^{n} Y_{j,n}^2 = \sum_{j=1}^{n} (Y_{j,n}^2-{1\over {n}}) +
m_n \build{\lra}_{n\to\infty}^{L^2} \lim_{n\to\infty} m_n= {1\over
  {\sigma(M)}} \int_M f\; d\sigma.$$
We are done with the second order term.
We finish the proof by studying the first order one.
\begin{eqnarray*}
\sum_{j=1}^{n} f_{j,n} \left(X_{j,n}+{T\over {n}}\right) &=& \sum_{j=1}^{n}
f^0{j,n} X_{j,n} + m_n T \\
&=& m_n T + \sum_{j=1}^{n} f_{j,n}^0 Y_{j,n} \\
&=& m_n T + W(\sum_{j=1}^{n} f_{j,n}^0 \1_{F_{j,n}}) \\
&=& W(f_n^0) + m_n T \build{\lra}_{n\to\infty}^{L^2} W(f^0) + {T\over
{\sigma(M)}} \int_M f\; d\sigma.
\end{eqnarray*}

We have proved that
$${1\over i} \sum_{j=1}^{n} f_{j,n} (e^{i(X_{j,n} + {T\over
{n}})}-1) \build{\lra}_{n\to\infty}^{L^2} W(f^0) + \left(T+{i\over 2}\right)
{1\over {\sigma(M)}} \int_M f\; d\sigma.$$

This limit does not depend on the choice of the sequence $(f_n)$. Thus
the sequence $(I_n(f_n))$ converges also to a limit $I(f)$ that does
not depend on the choice of $(f_n)$ and whose law is the law announced
in the theorem. \qed

\subsection{Meaning of the variable $T$}
\label{meaning_of_T}
As a conclusion for this chapter, we will spend a few lines to suggest a
geometric interpretation for the variable $T$, whose meaning could seem to be
quite mysterious.

In a deterministic setting with a smooth connection $\omega$, a construction
similar to that of the map $f\mapsto I(f)$ would have given the map:
$$f\mapsto \int_M f F(\omega),$$
where $F(\omega)$ is the curvature $2$-form
of $\omega$. As long as we consider zero-mean functions, the comparison
between this formula and (\ref{TAF}) is in agreement with the heuristic
principle saying that the curvature of a Yang-Mills random connection is a
white noise, as explained in the introduction.

If we take the function $f$ identically equal to $1$ in the deterministic
setting, we get the total curvature $\int_M F(\omega)$ of the fiber bundle $P$
on which $\omega$ lives. This quantity is well known to be independent of
$\omega$ and to be a topological invariant of $P$, namely its first Chern
class. The probabilistic counterpart of this total curvature seems then to be
$I(1)=T$, droping out the imaginary part.  This discussion becomes really
meaningful when $M$ is closed, because $P$ is not necessarily trivial. We
mentioned at the beginning of the discretization procedure in section
\ref{discrete_holonomy} that we had lost any topological information about the
structure of $P$.  If we compute the ``random Chern class'' of $P$ at the end
of the construction, we find a weighted sum of all possible Chern classes,
with the smallest weights for the most complicated types of bundles. This was
already suggested by Witten \cite{Witten}.

On the other hand, we can change our point of view in the following way: we
have an expression of the random holonomy which depends explicitely on the
Chern class of $P$. So if we replace $T$ by a deterministic multiple of $2\pi$
in the definition \ref{def_u1}, we are able to construct a random holonomy
consistent with any prescribed type of bundle $P$.

\section{Square-integrability of the double-layer potential}

In this section, we prove the theorem \ref{dlp_l2}. We claim that it is
enough to prove the theorem on closed surfaces. Indeed, we did not use the
square-integrability of the double-layer potential to prove the results
\ref{dlp_h1} and \ref{dlp_bnd}, which show that the result on a surface with
boundary can be deduced from the result on a minimal closure of this surface.
Thus, we assume that $M$ is closed.

\begin{proposition} There exists $R_F>0$ such that for all embedded path $c\in
  PM$ such that $\ell(c)<{1\over 4}R_F$, the double layer potential $u_c$ of
  $c$ is in $L^\infty(M)$.
\end{proposition}

This proposition implies obviously the theorem. It implies even more, namely
that the double layer potential of any path is in $L^\infty$.\\

\pf. The proof relies on three facts. The first one is that we know the
divergence near the diagonal of the Green function in an open subset of
$\R^2$. The second one is that, according to a classical theorem due to Gauss
\cite{Chern}, any metric on $M$ is locally conformally flat. The third point
is that the Green function is conformally invariant.

Since $M$ is compact, the second remark implies that there exists a radius
$R_F$ such that any geodesic ball of $M$ of radius smaller than $R_F$ is
conformally flat. Let us choose an embedded path $c$ of length smaller than
$R_F/4$. For each $r>0$, we denote by $B_r$ the ball $B(c(0),r)$.  Since
$\ell(c)$ is smaller than $R_F/4$, $c$ is contained in $B_{{1\over 2}R_F}$.
Since the Green function $G$ is smooth outside the diagonal, $u_c$ is smooth
outside $B_{{1\over 2}R_F}$. It is enough now to prove that it is bounded on
$B_{{3\over 4}R_F}$ for example. Set $r={3\over 4}R_F$.

The values of $u_c$ inside $\overline B_r$ depend only on the restriction of
$G$ to $\overline B_r\times \overline B_r$. On this set, $G$ satisfies $\Delta
G_x=\delta_x - {1\over \sigma (M)}$. Our idea is to substract smooth functions
to $G$ until we get something easier to compute than $G$ itself. Denote by
$G^0$ the solution of:
$$\cases{\Delta G^0_x=-{1\over {\sigma(M)}} \cr G^0_x(y)=G_x(y) \;\; \forall
  y \in \partial B_r.}$$
It is a smooth function inside $B_r$. The function $G^1$
defined by $G^1=G-G^0$ satisfies
$$\cases{\Delta G^1_x=\delta_x \cr G^1_x(y)=0 \;\; \forall y \in \partial
  B_r.}$$
So, $G^1$ is the fundamental solution of $\Delta$ inside $B_r$ with
Dirichlet boundary conditions. We can decompose $u_c(x)$ for any $x$ in $B_r$
according to:
$$u_c(x)= \int_c *dG^0_x + \int_c *dG^1_x = u_c^0(x)+u_c^1(x).$$
The first
term is smooth and we are led to study the second. This is where conformally
flat coordinates are useful: we choose a local chart $\varphi: U
\lra B_r$, where $U$ is an open subset of $\R^2$, such that the pull-back of
the metric $g$ of $M$ by $\varphi$ is conformally equivalent to $dx^2 + dy^2$
on $U$. The point is that $\varphi^* G^1$ is not only the fundamental solution
of $\Delta_{\varphi^* g}$ with respect to the measure induced by $\varphi^*
g$, but also the fundamental solution of $\Delta_0=\partial_x^2 +\partial_y^2$
with respect to the flat metric on $U$, by conformal invariance of the Green
function. This tells us that $\varphi^* G^1$ diverges like ${1 \over {4\pi}}
\log d(\cdot,\cdot)$ on the diagonal. In other words, there exists a smooth
function $G^2$ on $B_r \times B_r$ such that $\varphi^* G^1={1 \over {4\pi}}
\log d(\cdot,\cdot) + \varphi^*G^2$.

What we want is to prove that $u_c^1(x)$ is bounded inside $B_r$. It is
equivalent to prove that $u_c^1 \circ \varphi$ is bounded on $U$. But for any
$y\in U$,
$$u_c^1 \circ \varphi (y)= \int_c *dG^1_{\varphi(y)} = \int_{\varphi^{-1}(c)}
*d(\varphi^*G^1_y).$$
Note that in this last term, the Hodge operator $*$ is
that of the metric $\varphi^* g$. It is not the same operator as that of the
flat metric on $U$. Fortunately, the fact that these two metrics are
conformally equivalent implies that their Hodge operators are pointwise
proportionnal, i.e. one is deduced from the other by the multiplication by a
positive bounded smooth function. Thus, it is sufficient to prove that
$y\mapsto \int_{\varphi^{-1}(c)} *d\varphi^*G^1_y$ is bounded, the Hodge
operator being now that of the flat metric. As already noticed, we can remove
a smooth part of $G^1$ and keep only the part
$${1\over {4\pi}} \int_{\varphi^{-1}(c)} *d\log d(y,\cdot).$$
A short
computation shows that $1/4\pi *d\log d(y,\cdot)$ is nothing but $1/2\pi$
times the angle form $d\theta$ of the polar coordinates centered at $y$. This
allows us to estimate very easily the integral of this form along a path.

For example, we know that the integral of this form along a simple loop is
bounded by 1. On the other hand, it is obvious and easy to prove by a direct
computation that the integral of this form along a straight segment is bounded
by $1/2$.

Consider the path $\varphi^{-1}(c)$. It is injective, hence it is possible to
transform it into a simple loop by concatenating it with a finite number of
segments. So, we can make the function that we want to estimate to be bounded
by adding to it a finite number of bounded functions.  This gives the result.
\qed


\pagevide
\chapter{Small scale structure in the semi-simple case}

The theorem \ref{abelian_main} shows that it is possible to construct the
Yang-Mills measure in a short and quite pleasant way when $G=U(1)$, using a
white noise on $M$ as main ingredient. Is it possible to do something similar
in general? The works of Sengupta and Driver \cite{Driver_lassos, Sengupta_S2,
  Sengupta_AMS} lead to an ambiguous answer to that question. Indeed, in these
works, the authors have constructed random holonomies, starting from a Lie
algebra-valued white noise on $M$. Nevertheless, the family of loops along
which they are able to define the holonomy is strongly dependent of a
particular choice of coordinates on $M$, as we explained in the introduction.
We think that this is more than a simple technical problem. Although there
might exist some generalization of the construction made in section
\ref{rh_as_wn_functional}, we will show that a white noise is probably not the
right object to start with.

Our idea is the following. Was it possible to realize the random holonomy
using a white noise, it would be possible to find a lot of information by
looking at the random holonomy at small scale, i.e. along very small
loops. For example, the theorem \ref{construct_wn} basically says that when
$G=U(1)$, almost all the information about the holonomy along homologically
trivial loops is available at infinitesimally small scale. We prove that, when
$G$ is semi-simple, there is no information at all available at
infinitesimally small scale, at least when one looks at it in the same way
that we did in the Abelian case.

\section{Statement of a zero-one law}

We begin by stating the main result. The surface is $(M,\sigma)$ as usual.
We assume that $G$ is a compact connected semi-simple Lie group, for
example $SU(2)$. We choose $x_1,\ldots,x_p$ in $G$ and consider the
probability space $(\Omega_M,\mu_M(x_1,\ldots,x_p))$. 

Let $L$ be a simple loop on $M$ which is the boundary of an open set
$D$ diffeomorphic to a disk. For each $n\geq 0$, consider a graph on
$D$ which has exactly $n$ faces $F_{1,n},\ldots,F_{n,n}$ such that
$\sigma(F_{i,n})={{\sigma(D)}\over n}$ for each $i$. This is very
similar to the situation described in the section \ref{small_scale}.

Was $G$ Abelian, we would have the equality of cycles $L=\partial F_{1,n} +
\ldots +\partial F_{n,n}$, provided orientations are well
chosen. This would imply $\h_L=\h_{\partial F_{1,n}} \ldots \h_{\partial
  F_{n,n}}$ and for any function $f$ continuous on $G/\Ad=G$,
$$E[f(\h_L)|\h_{\partial F_{1,n }},\ldots,\h_{\partial
F_{n,n}}]=f(\h_L).$$

When $G$ is semi-simple, the situtation is the opposite.

\begin{theorem} \label{01}
For any function $f$ continuous on $G/\Ad$, the following convergence
holds:
$$E[f(\h_L)|\h_{\partial F_{1,n }},\ldots,\h_{\partial
F_{n,n}}] \build{\lra}_{n\to\infty}^{L^2} E f(\h_L).$$
\end{theorem}

\section{Proof of the zero-one law}
\subsection{Computation of the conditional expectation}

In this section, we will compute the conditional expectation
appearing in the statement of the theorem, keeping $n$ fixed. We
abbreviate $F_{i,n}$ in $F_i$. 

For each $F_i$, consider a sequence $(L_{i,k})_{k\geq 0}$ of simple loops
whose image is inside the interior of $F_i$ and such that $L_{i,k}
\build{\lra}_{k\to\infty}^{d_1} \partial F_i$. The proposition \ref{cont_5}
shows that the following convergence holds in probability:
$$(\h_{L_{1,k}},\ldots,\h_{L_{n,k}}) \build{\lra}_{k\to\infty}^{P}
(\h_{\partial F_1},\ldots,\h_{\partial F_n}).$$
Let $f$ and $f_1$ be
continuous functions on $G/\Ad$ and $(G/\Ad)^n$ respectively. This convergence
implies the following one:
$$E[f(\h_L)f_1(\h_{\partial F_1},\ldots,\h_{\partial
F_n})]=\lim_{k\to\infty}
E[f(\h_L)f_1(\h_{L_{1,k}},\ldots,\h_{L_{n,k}})].$$
We are led to the computation of the second expectation, keeping $k$
fixed. We abreviate temporarily $L_{i,k}$ by $L_i$. 

\begin{figure}[h]
\begin{center}
\input{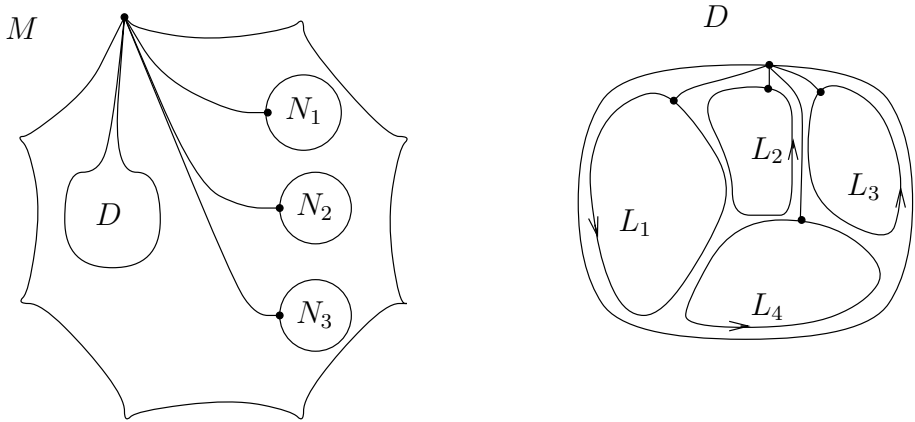}
\end{center}
\caption{Aspect of $\Gamma$ outside and inside $D$.
\label{cond}}
\end{figure}

We construct a particular graph on $M$ such that $L,L_1,\ldots,L_n \in
\Gamma^*$ (see fig. \ref{cond}). Outside $D$, it has only one face. Its
support contains the components $N_1,\ldots,N_p$ of $\partial M$, paths
$c_1,\ldots,c_p$ joining $L(0)$ to the $N_i$'s and simple loops
$a_1,\ldots,a_g,b_1,\ldots,b_g$ that represent a basis of the $H_1$ of a
minimal closure of $M$. 

The boundary of the unique face of $\Gamma$ outside $D$ is $L^{-1} c_1 N_1
c_1^{-1} \ldots c_p N_p c_p^{-1}\times$ $\times[a_1,b_1] \ldots [a_g,b_g]$,
where $[a,b]$ denotes the commutator $aba^{-1}b^{-1}$ of $a$ and $b$. This
notation is the same as that of joint conjugacy classes, but the context will
always make our meaning clear.  Inside $D$, the support of $\Gamma$ contains
$n$ paths $d_1,\ldots,d_n$ joining $L(0)$ to the $L_i(0)$. These paths meet
pairwise only at $L(0)$. Inside $D$, this graph has $n+1$ faces: $n$ whose
boundaries are the $L_i$'s and the last with boundary $Ld_n L_n^{-1} d_n^{-1}
\ldots d_1 L_1^{-1} d_1^{-1}$. We denote by $t$ the surface of this last face
and, for each $i=1,\ldots,n$, by $t_i$ the surface of the face bounded by
$L_i$.

We compute $E[f(\h_L) f_1(\h_{L_1},\ldots,\h_{L_n})]$ in this graph. It is
equal to
\begin{eqnarray*}
&& \hskip -.8cm{1\over
{Z(x_1,\ldots,x_p)}} \int_{G^\Gamma} f([h_L])
f_1([h_{L_1}],\ldots,[h_{L_n}]) p_{t_1}(h_{L_1})\ldots
p_{t_n}(h_{L_n}) \\
&& \hskip 3cm p_t(h_{d_1}^{-1} h_{L_1}^{-1} h_{d_1} \ldots h_{d_n}^{-1}
h_{L_n}^{-1} 
h_{d_n} h_L)\\
&&  p_{\sigma(D^c)}([h_{b_g}^{-1},h_{a_g}^{-1}] \ldots
[h_{b_1}^{-1},h_{a_1}^{-1}] h_{c_p}^{-1} h_{N_p} h_{c_p} \ldots
h_{c_1}^{-1} h_{N_1} h_{c_1} h_L^{-1}) \;  d\nu_{x_1} \ldots d\nu_{x_p} dg'\\
&=& {1\over{Z(x_1,\ldots,x_p)}} \int_{G^{2(g+n)+p+1}} f([g])
f_1([g_1],\ldots,[g_n]) p_{t_1}(g_1) \ldots p_{t_n}(g_n)\\
&& \hskip 0cm p_t(y_1^{-1} g_1^{-1} y_1 \ldots y_n^{-1} g_n^{-1} y_n g)
p_{\sigma(D^c)}([a_g,b_g] \ldots [a_1,b_1] z_p^{-1} x_p z_p
\ldots z_1^{-1} x_1 z_1 g^{-1}) \\
&& \hskip 2.8cm dg \; dg_1 \ldots dg_n \; dy_1 \ldots dy_n \; dz_1 \ldots dz_p \; da_1
\ldots da_g \; db_1 \ldots db_g.
\end{eqnarray*}
We used the fact that under $\nu_{x_1}\ldots \nu_{x_p} dg'$, $h_{N_i}=x_i$
a.s. and $h_L$ and all $h_{L_i}$, $h_{a_i}$, $h_{b_i}$, $h_{c_i}$,$h_{d_i}$
are uniform and independent. When $k$ tends to infinity, each $t_i$ tends to
${{\sigma(D)}\over n}$ and $t$ tends to zero. So, according to
\ref{conditioned_Haar_measure}, we can drop $g$ and replace it by $y_n^{-1}
g_n y_n \ldots y_1^{-1} g_1 y_1$. This terminates the proof of the following
proposition:

\begin{proposition} \label{comput}
The following equality holds:
\begin{eqnarray*}
&& \hskip -.5cm E[f(\h_L)|\h_{\partial F_{1,n }},\ldots,\h_{\partial F_{n,n}}]
=\\
&& \hskip -.7cm \left(\int p_{\sigma(D^c)}([a_g,b_g] \ldots [a_1,b_1] z_p^{-1} x_p z_p
\ldots z_1^{-1} x_1 z_1 y_1^{-1} K_1 y_1 \ldots y_n^{-1} K_n y_n)\;
dy_i dz_i da_i db_i
\right)^{-1}\!\!\! \cdot \\
&& \int f(y_n^{-1} K_n y_n \ldots y_1^{-1} K_1 y_1)\\
&& \hskip -.5cm p_{\sigma(D^c)}([a_g,b_g] \ldots [a_1,b_1] z_p^{-1} x_p z_p
\ldots z_1^{-1} x_1 z_1 y_1^{-1} K_1^{-1} y_1 \ldots y_n^{-1} K_n^{-1} y_n)\;
dy_i dz_i da_i db_i 
\end{eqnarray*}
where $K_1,\ldots,K_n$ are arbitrary representatives of $\h_{\partial
F_{1,n}},\ldots,\h_{\partial F_{n,n}}$. 
\end{proposition}

Since we are dealing with functions on $G/\Ad$, it is natural to use the
theory of characters on $G$. We give a very short account of the results that
we will use. For a detailed presentation of the subject, see for example
\cite{Broecker_tom_Dieck,Simon}.

\subsection{Characters of a semi-simple Lie group}
\label{carac}

A representation of $G$ is a smooth morphism of groups $\rho$ from $G$ into
the linear group of some $\C^n$. The integer $n$ is the dimension of $\rho$.
Since $G$ is compact, we may always assume that $\rho(G)$ is included in the
unitary group. A representation is said to be irreducible if there are no
subspaces of $\C^n$ invariant by all $\rho(g)$ except $\C^n$ and $0$. Two
representations $\rho_1$ and $\rho_2$ of same dimension $n$ are said to be
equivalent if there exists a linear isomorphism $\phi$ of $\C^n$ such that
$\rho_1 \circ \phi = \phi \circ \rho_2$.  . The character of $\rho$ is the
$\C$-valued function $\chi_\rho$ defined on $G$ by $\chi_\rho(g)=\tr \rho(g)$.
Two equivalent representations have the same character. In fact, this is also
a sufficient condition of equivalence.  The usual properties of the trace
imply that it is a central function on $G$, that is, $\chi_\rho(g)$ depends
only on the conjugacy class of $g$. Since all representations are unitary, the
relation $\chi_\rho(g^{-1})=\chi_\rho(g)^*$ holds, the star denoting the
complex conjugation. Note that $\chi_\rho(1)=\dim \rho$. The main theorem is
the following:

\begin{theorem}[Peter-Weyl theorem] The set of characters of all equivalence
  classes of irreducible representations of $G$ is an orthonormal basis of the
  space of central square-integrable functions on $(G,dg)$. Moreover, the
  algebra generated by this set is dense in the set of continuous $\C$-valued
  central functions on $G$ endowed with the uniform norm.
\end{theorem}

According to this theorem, any continuous function can be approximated by
linear combinations of products of characters. But it is a fact that such
combinations can always be written as linear combinations of characters of
irreducible representations. Thus it is sufficient to prove our theorem when
$f$ is the character of an irreducible representation.

Characters satisfy orthogonality relations that give rise to useful
formulas. We will mainly use two of them.

\begin{proposition} For any $x,y,z\in G$ and for any irreducible
representation $\alpha$ of $G$,
\begin{eqnarray} \label{rel_1}
\int_G \chi_\alpha (xyx^{-1}z) \; dx \! &=&\! {1\over {\dim \alpha}}
\chi_\alpha(x) \chi_\alpha(y),\\
 \label{rel_2}
\int_{G^2} \chi_\alpha([a,b]x)\; da db\! &=&\! {1\over {(\dim \alpha)^2}} 
\chi_\alpha(x).
\end{eqnarray}
\end{proposition} 

Let us endow $G$ with its biinvariant metric normalized to have total volume
equal to $1$. This metric gives rise to a Laplace operator on $G$. A
remarkable property of the characters is that they are eigenfunctions for this
operator. More precisely, for any irreducible representation $\alpha$, there
exists a positive real number $c_2(\alpha)$ such that
$$\Delta \chi_\alpha = -c_2(\alpha) \chi_\alpha.$$

A nice application of these properties is the computation of the character
expansion of the heat kernel on $G$. Let us denote by $\widehat G$ the set
of classes of irreducible representations of $G$.

\begin{proposition} \label{heat_kernel_expansion}
The following equality holds in $L^2(G,dg)$ and
also pointwise on $G$:
$$p_t =\sum_{\beta\in\widehat G} \dim \beta\; e^{-{{c_2(\beta)}\over 2} t}
\chi_\beta.$$
\end{proposition}

\pf. We first prove the $L^2$ convergence. For any $t>0$, the function
$p_t$ is a central $L^2$ function on $G$. Thus it admits a
decomposition 
$$p_t=\sum_{\beta\in\widehat G} C_\beta(t) \chi_\beta.$$
The differential equation $({1\over 2}\Delta - \partial_t) p_t=0$
implies that $({1\over 2}c_2(\beta)+ \partial_t) C_\beta(t) =0.$
Thus, $C_\beta(t)=C_\beta(0) e^{-c_2(\beta) {t\over 2}}$.
The constants $C_\beta(0)$ are determined using the fact that $p_t$
tends to $\delta_1$ as $t$ tends to $0$. It is easily checked that
$C_\beta(0)=\dim \beta$ is a convenient choice. Formally, it amounts
to check that $\sum_\beta \dim \beta \chi_\beta = \delta_1$. 

To see that the convergence holds pointwise, note that the expansion of $p_t$
is a series of continuous functions that converges normally. Since $p_t$ is
continuous, it is not only the sum of this series in the $L^2$ sense, but also
in the sense of the uniform convergence.\qed

\subsection{Character computations}

We go back to the big expression obtained in proposition \ref{comput}. From
now on, we fix an irreducible representation $\alpha$ and put $f=\chi_\alpha$.
We compute the numerator $N_n$ of the conditional expectation. The
computations using characters presented in this section and the next one are
very close to those done by Witten in \cite{Witten}, when he expands
explicitely partition functions in order to compute the symplectic volume of
the moduli space of flat connections. 

We begin by developing the heat kernel $p_{\sigma(D^c)}$ using proposition
\ref{heat_kernel_expansion}. We get a sum over $\beta$ of integrals of
$\chi_\alpha$ of something times $\chi_\beta$ of something else. We integrate
over the variables that appear only as arguments of $\chi_\beta$, first $a_1$
and $b_1$, and so on until $a_g$ and $b_g$, using (\ref{rel_2}). Each
integration against $a_i$ and $b_i$ produces a factor ${1\over {(\dim
    \beta)^2}}$. Then we integrate against $z_1,\ldots,z_n$ using
(\ref{rel_1}). Each integration gives a factor $\chi_\beta(x_i) / \dim\beta$.
At this stage, the arguments of $\chi_\alpha$ and $\chi_\beta$ are inverse of
each other. Using the relation $\chi_\beta(g^{-1})=\chi_\beta^*(g)$, we
obtain:
\begin{eqnarray*}
&& \hskip 0cm N_n=\sum_{\beta\in\widehat G} (\dim\beta)^{1-2g-p} \;
    e^{-{{c_2(\beta)}\over 2} \sigma(D^c)} \prod_{i=1}^p \chi_\beta(x_i)
    \cdot\\ 
&& \hskip 4cm \int_{G^n} \chi_\alpha
\chi_\beta^*(y_n^{-1} K_n y_n \ldots y_1^{-1} K_1 y_1)\; dy_1 \ldots dy_n.
\end{eqnarray*}
We use the formal expansion of the central function $\chi_\alpha
\chi_\beta^* = \sum_{\gamma\in\widehat G} (\chi_\alpha \chi_\beta^*,
\chi_\gamma)_{L^2} \; \chi_\gamma.$ We have now a sum over $\beta$ and
$\gamma$ but only $\chi_\gamma$ under the integral. We integrate against
$y_1,\ldots,y_n$ using (\ref{rel_1}). Note that the factor ${1\over
  {\dim\gamma}}$ produced by the last integration cancels out with the
remaining $\chi_\gamma(1)$ . We find
\begin{eqnarray}
&& \hskip -2cm  N_n=\sum_{\beta,\gamma\in\widehat G} (\dim\beta)^{1-2g-p}
(\dim\gamma)^{-(n-1)} 
\; e^{-{{c_2(\beta)}\over 2} \sigma(D^c)} (\chi_\alpha \chi_\beta^*,
\chi_\gamma)_{L^2} \nonumber \\
&& \hskip 6cm \prod_{i=1}^p \chi_\beta(x_i) \prod_{i=1}^n
\chi_\gamma(K_n),
\end{eqnarray} 
with $(\chi_\alpha \chi_\beta^*,\chi_\gamma)_{L^2}=\int_G
\chi_\alpha(g) \chi_\beta^*(g) \chi_\gamma^*(g) \; dg = (\chi_\alpha
,\chi_\beta \chi_\gamma)_{L^2}$. 

A similar and simpler computation leads to the following expression
for the denominator $D_n$ of \ref{comput}:
\begin{equation}
D_n=\sum_{\beta\in\widehat G} (\dim\beta)^{2-2g-p-n}
e^{-{{c_2(\beta)}\over 2} \sigma(D^c)} \prod_{i=1}^p \chi_\beta(x_p)
\prod_{i=1}^n \chi_\beta(K_i).
\end{equation}

Remark that it is equivalent to evaluate a character at $K_i$ or on
$\h_{\partial F_i}$. In order to prove the convergence, we need to know the
asymptotic behaviour of an expression like
$${1\over {(\dim\beta)^{n-1}}} \prod_{i=1}^n \chi_\beta(\h_{\partial
F_i}).$$

The following result will be proved in the next section:

\begin{proposition} \label{key01}
For any $\beta\in\widehat G$, the following
convergence holds:
$${1\over {(\dim\beta)^{n-1}}} \prod_{i=1}^n \chi_\beta(\h_{\partial
F_i}) \build{\lra}_{n\to\infty}^{L^2} \dim\beta\; e^{-{{c_2(\beta)}\over
2} \sigma(D)}.$$
\end{proposition}

The inequality $|\chi_\beta(g)|\leq \dim\beta$ shows that the sequence is
uniformly bounded by $\dim \beta$. This allows us to permute this convergence
with the summation over $\beta$. We get the following $L^2$-limit for the
numerator:
$$\sum_{\beta,\gamma\in\widehat G} (\dim\beta)^{1-2g-p} e^{-{{c_2(\beta)}\over
    2} \sigma(D^c)} \prod_{i=1}^p \chi_\beta(x_p) \; (\chi_\alpha,\chi_\beta
\chi_\gamma)_{L^2} \dim\gamma e^{-{{c_2(\gamma)}\over 2} \sigma(D)}.$$
By the
same kind of arguments that we used to derive the expressions of $D_n$ and
$N_n$, it is easy to check that this expression is equal to the following:
\begin{equation} \label{num01}
\int_{G^{p+2g+1}} \!\!\!\!\!\!\!\! \chi_\alpha(g) p_{\sigma(D)}(g^{-1})
p_{\sigma(D^c)}(g^{-1} z_1 x_1 z_1^{-1} \ldots z_p x_p
z_p^{-1}[a_1,b_1]\ldots [a_g,b_g])\; dg dz_i da_i db_i.
\end{equation}

For the denominator, we find the following $L^2$-limit:
$$\sum_{\beta\in\widehat G} (\dim\beta)^{2-2g-p} e^{-{{c_2(\beta)}\over 2}
(\sigma(D^c)+\sigma(D))} \prod_{i=1}^p \chi_\beta(x_p),$$
which is equal to:
\begin{equation} \label{den01}
\int_{G^{p+2g}} p_{\sigma(M)}(z_1 x_1 z_1^{-1} \ldots z_p x_p
z_p^{-1} [a_1,b_1] \ldots [a_g,b_g]) \; dz_i da_i db_i \; = \;
Z_M(x_1,\ldots,x_p). 
\end{equation}
Using the fact that $p_{\sigma(D)}(g^{-1})=p_{\sigma(D)}(g)$, we see
that the quotient of (\ref{num01}) by (\ref{den01}) is equal to $E
\chi_\alpha(\h_L)$. This proves the theorem, up to the proposition \ref{key01}
that was admitted.

\subsection{Zero-one law on the plane}

In order to prove proposition \ref{key01}, we begin with the following
result, which can be seen as a reformulation of the zero-one law when
the manifold $M$ is the plane $\R^2$.

\begin{proposition} \label{01_plane}
Let $(B^n_t)_{t\in\sR_+,n\geq 0}$ be a sequence of independent Brownian
motions on 
$G$. For any irreducible representation $\beta$ of $G$ and any
positive real number $T$, the following convergence holds in
distribution:
$${1\over {(\dim\beta)^{n-1}}} \prod_{i=1}^n \chi_\beta(B^i_{T\over
n}) \build{\lra}_{n\to\infty}^{{\rm law}} \dim\beta \; e^{-{{c_2(\beta)}\over
2} T}.$$
\end{proposition}

This result is really the center of the whole proof of the theorem. It is the
place where the fact that $G$ is semi-simple will be used, in the following
way. Given a representation $\rho$ of $G$, the differential at $1$ of $\rho$
is a linear map from $T_1G$, the Lie algebra of $G$, into $\l(\C^n)$. The
following statement is proved in Bourbaki ({\it Lie}, chap. I, $\S$ 6, N$^o$2,
corollary of the th.  1) \cite{Bourbaki_Lie_1}:

\begin{proposition} \label{sl}
Let $G$ be a semi-simple group and $\rho$ a
representation of $G$ of dimension $n$. Then
$$d_1\rho(T_1 G) \subset {\hbox{\got sl}}_n(\C),$$
where ${\hbox{\got
    sl}}_n(\C)$ denotes the set of endomorphisms of $\C^n$ whose trace is
equal to zero.
\end{proposition}

\pf \textsl{of proposition \ref{01_plane}}: We consider a Brownian motion
$(B_t)$ on $G$ and study the process $\chi_\alpha(B_t)$. There is a convenient
way to represent $(B_t)$, as a solution of a Stratonovich stochastic
differential equation \cite{Ikeda_Watanabe}. Recall that the data of a
biinvariant metric on $G$ is equivalent to that of a scalar product on the Lie
algebra of $G$, invariant by adjunction. Let $(X_1,\ldots,X_{\dim G})$ be a
basis of $T_1G$ orthonormal for this scalar product. Each $X_i$ is seen as a
left-invariant vector field on $G$. Let $W^1,\ldots,W^{\dim G}$ be independent
real Brownian motions. Then the Brownian motion on $G$ satisfies:
\begin{equation} \label{eds}
\cases{ dB_t= \sum_{i=1}^{\dim G} X_i \circ dW^i_t \cr 
        B_0 = 1}
\end{equation}
The meaning of this notation is that, for any continuous function $f$
on $G$,
$$f(B_t)=f(1) + \sum_{i=1}^{\dim G} \int_0^t X_i f(B_s) \;
dW^i_s + {1 \over 2} \int_0^t \Delta f(B_s)\; ds.$$
We apply this relation to $f=\chi_\alpha$. Using
$\chi_\alpha(1)=\dim\alpha$ and $\Delta \chi_\alpha=-c_2(\alpha)
\chi_\alpha$, it becomes:
\begin{eqnarray*}
&&{{\chi_\alpha(B_t)}\over {\dim\alpha}} = 1-{{c_2(\alpha) t}\over 2}
+ {1\over {\dim\alpha}} \sum_{i=1}^{\dim G} \int_0^t X_i \chi_\alpha
(B_s) \; dW^i_s \\
&& \hskip 6cm -{{c_2(\alpha)}\over {2 \dim\alpha}} \int_0^t
(\chi_\alpha(B_s) -\dim\alpha)\; ds.
\end{eqnarray*}
Set
$$Y_t=- {{c_2(\alpha)}\over {2 \dim\alpha}} \int_0^t (\chi_\alpha(B_s)
-\dim\alpha)\; ds,$$
$$Z_t={1\over {\dim\alpha}} \sum_{i=1}^{\dim G} \int_0^t X_i \chi_\alpha (B_s)
\; dW^i_s.$$
We keep the notation $\rho(g)$ for the biinvariant distance
$d(1,g)$ when $g\in G$. This $\rho$ has nothing to do with a representation
of $G$ !

For any $X\in T_1G$, we have:
$$d_1\chi_\alpha(X)={d\over {dt}}\biggl|_{t=0} \tr \alpha(\exp t X)=\tr
d_1\alpha(X)=0,$$
according to \ref{sl}. Thus the differential of
$\chi_\alpha$ at $1$ is zero. This implies that
$|\chi_\alpha(g)-\dim\alpha|=O(\rho(g)^2)$ in a neighbourhood of $1$. Using
the lemma \ref{estimate_heat_kernel},we get:
\begin{eqnarray}
E |Y_t|^2 &=& C E \left| \int_0^t (\chi_\alpha(B_s)-\dim\alpha)\; ds
\right|^2 \nonumber\\
&\leq & C t \; E \int_0^t |\chi_\alpha(B_s)-\dim\alpha |^2 \; ds
\nonumber\\
&\leq & C t \; E \int_0^t \rho(B_s)^4 \; ds \nonumber\\
&\leq & C t \int_0^t s^2 \; ds \nonumber\\
&\leq & C t^4. \label{Y_t}
\end{eqnarray}

For each $i=1,\ldots,\dim G$, the function $X_i\chi_\alpha$ is smooth and
$X_i\chi_\alpha(1)=0$. Thus, $|X_i\chi_\alpha(g)|=O(\rho(g))$ in a
neighbourhood of $1$. Thus,
\begin{eqnarray}
E |Z_t^2| &=& C E \left| \sum_{i=1}^{\dim G} \int_0^t X_i
\chi_\alpha(B_s)\; dW^i_s \right|^2 \nonumber\\
&=& C \sum_{i=1}^{\dim G} E \left| \int_0^t X_i \chi_\alpha(B_s)\;
dW^i_s \right|^2 \nonumber\\
&=& C \sum_{i=1}^{\dim G} E \int_0^t |X_i\chi_\alpha(B_s)|^2 \; ds
\nonumber\\
&\leq & C \dim G \int_0^t E \rho(B_s)^2 \; ds \nonumber\\
&\leq & C t^2. \label{Z_t}
\end{eqnarray}

The preliminary study of the process $\chi_\alpha(B_t)$ is now
finished. We consider a sequence $(B^n_t)$ as in the statement of
\ref{01_plane}, defined on a probability space $(\Omega,P)$. We look
at the product 
$$X_n=\prod_{i=1}^n (1-{{c_2(\alpha) T}\over {2n}} + Z^i_{T/ n}+
Y^i_{T/ n}),$$
where the random variables with different exponents are
independent. We would like to take the logarithm of this product. This
requires some precautions. Set
$$\Omega_n=\left\{ \left|Z^i_{{T/ n}} + Y^i_{T/n}\right| < {1\over 3}
\;\;\forall i=1\ldots n\right\}.$$
A Chebishev inequality gives
$$P\left( \left| Z_t + Y_t \right| < {1\over 3} \right) \geq 1 - 9 E
|Z_t + Y_t|^2 \geq 1 - Ct^2,$$
implying
$$P(\Omega_n) \geq \left( 1 -{{CT^2}\over n^2} \right)^n
\build{\lra}_{n\to\infty}^{} 1.$$
We do not change any convergence in distribution on $\Omega$ if we
replace $X_n$ by $1$ outside $\Omega_n$. So we set
$$\widetilde X_n=X_n \1_{\Omega_n} + \1_{\Omega_n^c}.$$
Then $\Log \widetilde X_n$ is
well defined, $\Log$ being the principal determination of the complex
logarithm. In fact, we have more than that. If $n$ is such that ${{c_2(\alpha)
    T}\over {2n}}$ is smaller than ${1\over 6}$, then each factor of $\widetilde
X_n$ is of the form $(1-z)$ with $|z|<{1\over 2}$. For such a $z$,we have
$|\log(1-z) +z| \leq |z|^2.$ Thus, the equality $$\Log(\widetilde
X_n) =\1_{\Omega_n} \sum_{i=1}^n \log(1-{{c_2(\alpha) T}\over {2n}} + Z^i_{T/
  n}+Y^i_{T/ n})$$
 implies
\begin{eqnarray}
\hskip 0cm \left| \Log\widetilde X_n + \sum_{i=1}^n \left({{c_2(\alpha)
T}\over {2n}} - Z^i_{T/ n} - Y^i_{T/ n}\right) \right| &\leq& \1_{\Omega_n}
\sum_{i=1}^n \left| {{c_2(\alpha) T}\over {2n}} - Z^i_{T/ n} - Y^i_{T/
n}\right|^2 + \nonumber\\
&& \hskip -1cm\1_{\Omega_n^c} \left| \sum_{i=1}^n {{c_2(\alpha)
T}\over {2n}} - Z^i_{T/ n} - Y^i_{T/ n} \right|. \label{log}
\end{eqnarray}

The last term tends to $0$ in probability because $P(\Omega_n^c)$
tends to $0$. Using (\ref{Y_t}) and (\ref{Z_t}), we find
$$E \left| {{c_2(\alpha) T}\over {2n}} - Z^i_{T/ n} -
Y^i_{T/n}\right|^2  \leq {{C T^2}\over {n^2}},$$
so that the first term of the right hand side of (\ref{log}) tends to $0$
in $L^1$ norm.

Now let us study the left hand side. On one hand, (\ref{Y_t}) implies $\sum_i
Y^i_{T/n} \build{\lra}_{}^{L^1} 0$, because
$$E \left| \sum_{i=1}^n Y^i_{T/n} \right| \leq \sum_{i=1}^n E
|Y^i_{T/n}| \leq {C\over n}.$$

On the other hand, using \ref{Z_t} and the fact that $E Z^i_{T/n}=0$,
we also have $\sum_i Z^i_{T/n} \build{\lra}_{}^{L^1} 0$, because
$$E \left| \sum_{i=1}^n Z^i_{T/n} \right|^2 = \sum_{i=1}^n E
|Z^i_{T/n}|^2 \leq {C\over n}.$$

We deduce that $\Log \widetilde X_n$ converges in probability to
$-{{c_2(\alpha) T}\over 2}$. This implies that $\widetilde X_n$, and so
$X_n$, converge in probability to $e^{-{{c_2(\alpha) T}\over
2}}$. Finally, we get
$$\dim\alpha \; X_n \build{\lra}_{n\to\infty}^{\rm law} \dim\alpha
\; e^{-{{c_2(\alpha) T}\over 2}},$$
proving the proposition.\qed

The proof of the theorem is almost finished. It remains to prove that
the proposition \ref{01_plane} implies the proposition \ref{key01}. \\

\pf \textsc{of \ref{key01}}: We set $T=\sigma(D)$. Let $F$ be a continuous
function on $\R$.
\begin{eqnarray*}
&& \hskip -.7cm E\left[ F\left({1\over {(\dim\beta)^{n-1}}} \prod_{i=1}^n
    \chi_\beta(\h_{\partial 
F_i})\right)\right] = {1\over {Z(x_1,\ldots,x_p)}} \int F\left({1\over
{(\dim\beta)^{n-1}}} \prod_{i=1}^n \chi_\beta(g_i)\right) \\
&&  \hskip .5cm p_{\sigma(D^c)}([a_g,b_g] \ldots [a_1,b_1] z_p^{-1} x_p z_p
\ldots z_1^{-1} x_1 z_1 y_1^{-1} g_1^{-1} y_1 \ldots y_n^{-1} g_n^{-1}
y_n) \\
&& \hskip 5.5cm p_{T/n}(g_1) \ldots p_{T/n}(g_n)  \; dg_i da_i db_i dy_i dz_i.
\end{eqnarray*}
We develop the heat kernel $p_{\sigma(D^c)}$ and integrate against all
variables except $g_1,\ldots,g_n$. We find
\begin{eqnarray*}
&&\hskip -.5cm {1\over {Z(x_1,\ldots,x_p)}} \sum_{\gamma\in\widehat G}
(\dim\gamma)^{1-2g-p} e^{-{{c_2(\gamma)}\over 2} \sigma(D^c)}
\prod_{i=1}^p \chi_\gamma(x_p)\\
&& \hskip 2cm \int_{G^n} F\left({1\over {(\dim\beta)^{n-1}}}
\prod_{i=1}^n \chi_\beta(g_i)\right) {1\over {(\dim\gamma)^{n-1}}}
\prod_{i=1}^n \chi_\gamma(g_i) \\
&& \hskip 7cm  p_{T/n}(g_1) \ldots p_{T/n}(g_n)\;
dg_1 \ldots dg_n.
\end{eqnarray*}
The proposition \ref{01_plane} says exactly that this last integral converges
to 
$$F(\dim\beta e^{-{{c_2(\beta)}\over 2} \sigma(D)}) \dim\gamma
e^{-{{c_2(\gamma)}\over 2} \sigma(D)}.$$ 
Thus,
\begin{eqnarray*}
&& \hskip -.5 cm E \left[F\left({1\over {(\dim\beta)^{n-1}}} \prod_{i=1}^n \chi_\beta(\h_{\partial
F_i})\right)\right] \build{\lra}_{n\to\infty}^{} {1\over {Z(x_1,\ldots,x_p)}}
\times \\
&&  \sum_{\gamma\in\widehat G}
(\dim\gamma)^{2-2g-p} e^{-{{c_2(\gamma)}\over 2} \sigma(M)} {1\over
{(\dim\gamma)^{p-1}}} \prod_{i=1}^p \chi_\gamma(x_i) \; F(\dim\beta
e^{-{{c_2(\beta)}\over 2} \sigma(D)}) \\
&=& {{Z(x_1,\ldots,x_p)}\over {Z(x_1,\ldots,x_p)}} F(\dim\beta
e^{-{{c_2(\beta)}\over 2} \sigma(D)}).
\end{eqnarray*}

This proves that the announced convergence holds in distribution. We already
noted that the function $\chi_\beta$ is bounded on $G$ by $\dim\beta$. So, the
sequence that we study is uniformly bounded by $\dim\beta$. The result
follows, since the convergence in distribution of a uniformly bounded sequence
to a deterministic limit implies its $L^2$ convergence. \qed


\pagevide

\chapter{Surgery of surfaces}

\section{Markov property of the Yang-Mills field}
\label{cut_1}

Consider a surface $M$ which is cut into two pieces $M_1$ and $M_2$. Our aim
in the first part of this chapter is to understand the relationships between
the Yang-Mills measures on $M_1$, $M_2$ and $M$. We already met a question of
this kind when we constructed the random holonomy on a surface with boundary
starting from that on a minimal closure of this surface and this led us to
prove in a particular case the Markov property which is the object of the
theorem \ref{Markov}.

Let $M_1$ and $M_2$ be two oriented surfaces with boundary such that $\partial
M_1$ and $\partial M_2$ have at least $p$ boundary components, where $p>0$ is
arbitrary. Pick $p$ components $N_1,\ldots,N_p$ on $\partial M_1$ and $p$
others $N'_1,\ldots,N'_p$ on $\partial M_2$. For each $i=1,\ldots,p$, consider
an orientation-reversing diffeomorphism $\psi_i : N_i \lra N'_i$ and call
$M$ the result of the gluing of $M_1$ and $M_2$ along $\psi_1,\ldots,\psi_p$. 
Denote by $L_1,\ldots,L_p$ $p$ loops on $M$ whose images are the components of
the common boundary of $M_1$ and $M_2$, oriented as boundary components of
$M_1$.  

For the statement of the theorem \ref{Markov}, we consider the Yang-Mills
measures on $M$, $M_1$ and $M_2$ measures on the spaces $\f(LM,G)$,
$\f(LM_1,G)$ and $\f(LM_2,G)$ endowed with the $\sigma$-algebras $\a$, $\a_1$
and $\a_2$ generated by the variables of the form $\h_{l_1,\ldots,l_n}$. 

There are two natural sub-$\sigma$-algebras on the space $(\f(LM,G),\a)$,
namely $\widetilde \a_i= \sigma(\h_{l_1,\ldots,l_n}, l_k \in LM_i)$, $i=1,2$.
For any $i=1,2$, a function $f_i$ on $\f(LM_i,G)$ gives rise to a function
$\widetilde f_i$ on $\f(LM,G)$ and it is equivalent to say that $f_i$ is
$\a_i$-measurable or to say that $\widetilde f_i$ is $\widetilde \a_i$-
measurable, so that we identify $\a_i$ and $\widetilde \a_i$, as well as $f_i$
and $\widetilde f_i$.

\begin{theorem} \label{Markov} The $\sigma$-algebras $\a_1$ and $\a_2$ are
  independent on $\f(LM,G)$ under $\mu_M$ conditionally to the random variable
  $(\h_{L_1},\ldots,\h_{L_p})$. Moreover, let $f_1$ and $f_2$ be two
  measurable functions on $(\f(LM_1,G),\a_1)$ and $(\f(LM_2,G),\a_2)$
  respectively. Then the product $f_1 f_2$ can be seen as a $\a$-measurable
  function on $\f(LM,G)$ and for any $t_1,\ldots,t_p \in G/\Ad$, the following
  equality holds:
  \begin{equation}\label{indcondf1f2}
  \mu_M(t_1,\ldots,t_p )(f_1f_2) = \mu_{M_1}(t_1,\ldots,t_p)(f_1)
  \mu_{M_2}(t_1^{-1},\ldots,t_p^{-1})(f_2). 
  \end{equation}
  Finally, these properties remain true if we condition further the Yang-Mills
  measures with respect to the random holonomy along some boundary components
  or some interior loops on $M_1$ or $M_2$.
\end{theorem}

This theorem says two things. It says that the random holonomy on $M_1$ is
independent of that on $M_2$ conditionally to the holonomy along the common
boundary of $M_1$ and $M_2$ and it also says that the restriction to $\a_i$ of
the measure $\mu_M(t_1,\ldots,t_p)$ is equal via the identification mentioned
above to the measure $\mu_{M_i}(t_1,\ldots,t_p)$, that is, the Yang-Mills
measure on $M_i$. 


We prove first a discrete result, which is due, in a slightly different form,
to Becker and Sengupta \cite{Becker_Sengupta}. 

\begin{proposition} \label{discMarkov} Let $\Gamma$ be a graph on $M$. Let
  $\Gamma_1$ and $\Gamma_2$ be the graphs on $M_1$ and $M_2$ induced by
  $\Gamma$. Let $f$ be a continuous function on $G^{\Gamma_1}\times
  G^{\Gamma_2}$ invariant by both the gauge transformations on $\Gamma_1$ and
  $\Gamma_2$. Then $f$ gives rise to a gauge-invariant function on $G^\Gamma$,
  still denoted by $f$ and for any $x\in G^p$,
  $$\int_{G^{\Gamma_1}\times G^{\Gamma_2}} f \; dP^{\Gamma_1}(y_1,x)
  dP^{\Gamma_2}(x^{-1},y_2) = \int_{G^\Gamma} f \; dP^\Gamma (y_1,x,y_2),$$
  where $y_1$ and $y_2$ represent the values of the holonomy imposed along the
  remaining components of $\partial M_1$ and $\partial M_2$ and possibly other
  loops inside $M_1$ and $M_2$.
\end{proposition}

\pf. By a standard density argument, it is enough to prove this proposition
when $f$ is a product $f_1f_2$ of two gauge-invariant functions on
$G^{\Gamma_1}$ and $G^{\Gamma_2}$. In order to shorten the expressions, we
abbreviate expressions like $d\nu_{x_1} \ldots d\nu_{x_p}$ into $d\nu_x$ if
$x=(x_1,\ldots,x_p)$. So, if $f=f_1f_2$, 
\begin{eqnarray*}
\int_{G^{\Gamma_1}\times G^{\Gamma_2}}\!\!\! f_1f_2\;
dP^{\Gamma_1}(y_1,x)dP^{\Gamma_2}(x^{-1},y_2) = \! \int_{G^{\Gamma_1}} f_1 \;
dP^{\Gamma_1}(y_1,x) 
\! \int_{G^{\Gamma_2}} f_2 \; dP^{\Gamma_2}(x^{-1},y_2) \\
&& \hskip -13cm = {1\over {Z_{M_1}(y_1,x)}} \int_{G^{\Gamma_1}} f_1
D^{\Gamma_1} \; d\nu_{y_1} d\nu_x dg'^{\Gamma_1} {1\over
  {Z_{M_2}(x^{-1},y_2)}} \int_{G^{\Gamma_2}} f_2 
D^{\Gamma_2} \; d\nu_{y_2} d\nu_{x^{-1}} dg'^{\Gamma_2}. 
\end{eqnarray*}
\begin{equation} \label{bip} \end{equation}
The two integrals in the last expression are very similar and it is enough to
study the first one in order to study both of them. Denote by
$\Gamma_{1,\partial}$ the graph on $\partial M_1$ induced by $\Gamma_1$. The
measure $d\nu_x$ concerns precisely the variables associated with the edges of
$\Gamma_{1,\partial}$. Thus,
$$\int_{G^{\Gamma_1}} f_1 D^{\Gamma_1} \; d\nu_{y_1} d\nu_x dg'^{\Gamma_1} =
\int_{G^{\Gamma_{1,\partial}}} d\nu_x \int_{G^{\Gamma_1 \backslash
    \Gamma_{1,\partial}}} f_1 D^{\Gamma_1} \; d\nu_{y_1} dg'^{\Gamma_1}.$$

The function $I_1= \int_{G^{\Gamma_1 \backslash\Gamma_{1,\partial}}} f_1
D^{\Gamma_1} \; d\nu_{y_1} dg'^{\Gamma_1}$ is a function on
$G^{\Gamma_{\partial,1}}$ and we claim that it depends only on the values of
the functions $\h_{L_1},\ldots,\h_{L_p}$, which are also well defined on this
space.
  
This depends on the fact that $f_1 D^{\Gamma_1}$ is a gauge-invariant function
on $G^{\Gamma_1}$ and on the biinvariance of the Haar measure. Indeed,
consider a gauge transformation $j$ on $\Gamma_1$, which is equal to $1$,
except on vertices located on $\Gamma_{1,\partial}$. This gauge transformation
changes the value of the holonomies along the edges of $\Gamma_{1,\partial}$,
and $j$ can be choosed in such a way that these holonomies take any
prescribed values, provided the conjugacy classes of the holonomies along
$L_1,\ldots,L_p$ remain unchanged. So, we just need to prove that the value of
$I_1$ is preserved by the action of such a gauge transformation. But this
action also affects the holonomy along the other edges that meet
$\Gamma_{1,\partial}$, so that the gauge invariance of $f_1$ and
$D^{\Gamma_1}$ is not sufficient to conclude.  It is the biinvariance of the
Haar measure that allows at this point to forget about the effect of $j$ on
the edges outside $\Gamma_{1,\partial}$.

Let us write this formally. Denote
$\Gamma_1=\{a_1,\ldots,a_k,a_{k+1},\ldots,a_l,a_{l+1},\ldots,a_m\}$, where
$\Gamma_{1,\partial}=\{a_1,\ldots,a_k\}$ and $a_{k+1},\ldots,a_l$ correspond
to the loops along which the holonomy is imposed by the choice of $y_1$. The
point for these edges is that they do not meet $\partial M_1$, so that their
discrete holonomy is not affected by $j$. For any $i=1,\ldots,m$, denote
$g_i^j =j(a_i(1))^{-1} g_i j(a_i(0))$. We have
\begin{eqnarray*}
&& \hskip -.5cm I_1(g_1,\ldots,g_k) = \int_{G^{\Gamma_1\backslash
    \Gamma_{1,\partial}}} f_1 
. D^{\Gamma_1}(g_1,\ldots,g_m) \; d\nu_{y_1}(g_{k+1},\ldots,g_l)
dg_{l+1}\ldots dg_m \\
&& \hskip .7cm = \int_{G^{\Gamma_1\backslash \Gamma_{1,\partial}}} f_1
. D^{\Gamma_1}(g_1^j,\ldots,g_m^j) \; d\nu_{y_1}
dg_{l+1}\ldots dg_m \\
&& \hskip .7cm = \int_{G^{\Gamma_1\backslash \Gamma_{1,\partial}}} f_1
. D^{\Gamma_1}(g_1^j,\ldots,g_k^j,g_{k+1},\ldots,g_{l},g_{l+1}^j,\ldots,g_m^j)
\; d\nu_{y_1} dg_{l+1}\ldots dg_m \\ 
&& \hskip .7cm = \int_{G^{\Gamma_1\backslash \Gamma_{1,\partial}}} f_1
. D^{\Gamma_1}(g_1^j,\ldots,g_k^j,g_{k+1},\ldots\ldots,g_m) \;
d\nu_{y_1} dg_{l+1}\ldots dg_m \\
&& \hskip .7cm =  I_1(g_1^j,\ldots,g_k^j).
\end{eqnarray*}

Now we go back to the computation (\ref{bip}). Denoting
$\Gamma_{2,\partial}=\{g'_1,\ldots,g'_k\}$, in such a way that $g_i$ and
$g'_i$ correspond to same edge on $M$, we have
\begin{eqnarray*}
\int_{G^{\Gamma_1}\times G^{\Gamma_2}}\!\!\! f_1f_2\;
dP^{\Gamma_1}(y_1,x)dP^{\Gamma_2}(x^{-1},y') = && \\
&& \hskip -6cm  = {1\over{Z_{M_1}(y_1,x)Z_{M_2}(x^{-1},y_2)}}
\int_{G^{\Gamma_{1,\partial}}\times G^{\Gamma_{2,\partial}}} 
I_1(g_1,\ldots,g_k) I_2(g'_1,\ldots,g'_k) \; d\nu_x
d\nu_{x^{-1}} \\
&& \hskip -6cm = {1\over{Z_{M_1}(y_1,x)Z_{M_2}(x^{-1},y_2)}} 
\int_{G^p} I_1.I_2(g_1,\ldots,g_k)\; d\nu_x(g_1,\ldots,g_p)\\
&& \hskip -6cm = {1\over{Z_{M_1}(y_1,x)Z_{M_2}(x^{-1},y_2)}} 
\int_{G^{\Gamma_\partial}} d\nu_x \int_{G^{\Gamma \backslash \Gamma_\partial}}
f_1 D^{\Gamma_1} f_2 D^{\Gamma_2} \; d\nu_{y_1} d\nu_{y_2} dg'^{\Gamma} \\
&& \hskip -6cm = {{Z_M(y_1,x,y_2)}\over{Z_{M_1}(y_1,x)Z_{M_2}(x^{-1},y_2)}}
\int_{G^\Gamma} f_1 f_2 dP^\Gamma(y_1,x,y_2). 
\end{eqnarray*}

The equality of the normalization constants follows from the case where
$f_1f_2$ is identically equal to 1, and this finishes the proof. \qed

The relation between conditional partition functions that was just established
deserves to be stated separately. We shall discuss this result and another one
of the same kind in the second part of this chapter.

\begin{proposition} \label{Z_1}
  For all $x\in G^p$, $y_1\in G^{p_1}$ and
  $y_2\in G^{p_2}$, the following relation holds:
  $$Z_M(y_1,x,y_2)=Z_{M_1}(y_1,x)Z_{M_2}(x^{-1},y_2).$$
\end{proposition}

We can use the proposition \ref{discMarkov} to prove the theorem
\ref{Markov}.\\ 

\pf \textsc{ of} \ref{Markov}. First of all, the relation (\ref{indcondf1f2})
implies the conditional independence stated in the first part of the
theorem. Now, by the continuity of the random holonomy, it is
sufficient to prove that this relation holds for functions $f_1$ and $f_2$
that depend on the holonomies along loops that can be put into a graph. So we
consider two such functions $f_1$ and $f_2$, and choose a graph $\Gamma$ on
$M$ such that all the loops which we need to consider are in
$\Gamma^*$. Choose $y_1$ and $y_2$ as in the proposition \ref{discMarkov}. Let
$\Gamma_1$, $\Gamma_2$ and $\Gamma_\partial$ be the graphs on $M_1$, $M_2$ and
$M_1 \cap M_2$ induced by $\Gamma$. The functions $f_1$ and $f_2$ may be
regarded as gauge-invariant functions on $G^{\Gamma_2}$ and $G^{\Gamma_1}$
respectively. Pick also a gauge-invariant function $f_\partial$ on
$G^{\Gamma_\partial}$. This function depends only on the conjugacy classes of
the holonomies along $L_1,\ldots,L_p$. So, we have:
\begin{eqnarray*}
&& \hskip -.5cm \int_{G^\Gamma} f_1 f_2 f_\partial \; dP(y,y') = {1\over {Z_M(y,y')}}
\int_{G^\Gamma} f_1 f_2 f_\partial D^\Gamma \; d\nu_y d\nu_{y'} dg' \\
&& \hskip -0cm = {1\over {Z_M(y,y')}} \int_{G^p} Z(y,x,y') \; dx\;
{1\over {Z(y,x,y')}} \int_{G^\Gamma} f_1 f_2 f_\partial D^\Gamma\;
d\nu_{x} d\nu_y d\nu_{y'} dg' \\
&& \hskip -0cm ={1\over {Z_M(y,y')}} \int_{G^p} dx \;f_\partial(x) Z(y,x,y')
\int_{G^\Gamma} f_1f_2 \; dP^\Gamma(y_1,x,y_2) \\
&&\hskip -0cm = {1\over {Z_M(y,y')}} \int_{G^p} dx \;f_\partial(x) Z(y,x,y')
\int_{G^{\Gamma_1}} f_1 \; dP^{\Gamma_1}(y_1,x) \int_{G^{\Gamma_2}} f_2 \;
dP^{\Gamma_2}(x^{-1},y_2),
\end{eqnarray*}
thanks to the proposition \ref{discMarkov}. This last equality, which is true
for every function $f_\partial$ shows that
\begin{eqnarray*}
E[f_1f_2|\a_\partial] &=& \int_{G^{\Gamma_1}} f_1 \; dP^{\Gamma_1}(y_1,x)
\int_{G^{\Gamma_2}} f_2 \; dP^{\Gamma_2}(x^{-1},y_2) \\
&=& \mu_{M_1}(u_1,t)(f_1) \mu_{M_2}(t^{-1},u_2)(f_2) \\
&=& E[f_1|\a_\partial] E[f_2|\a_\partial],
\end{eqnarray*}
where $\a_\partial=\sigma(\h_{L_1},\ldots,\h_{L_p})$ and $t$, $u_1$, $u_2$ are
the conjugacy classes corresponding to $x$, $y_1$ and $y_2$. These last
equalities prove the theorem. \qed

\section{Study of an example}
\label{examp}

Consider a closed surface $M$ of genus two realized as the connected sum of
two tori. Let $M_1$ and $M_2$ denote the two halves of $M$. Set $L=\partial
M_1 =-\partial M_2$. Let $f_1$ be a function on $(\f(LM_1,G),\a_1)$ and $f_2$
a function on $(\f(LM_2,G),\a_2)$.  Then $f_1 f_2$ can be seen as a function
on $(\f(LM,G),\a)$. We just proved that, for all $t \in G/\Ad$, $\mu_M(t) (f_1
f_2)=\mu_{M_1}(t)(f_1) \mu_{M_2}(t^{-1})(f_2).$ Since $\mu_M=\int_{G/\Ad}
\mu_M(t) Z_M(t) \; dt$ (see proposition \ref{disint}), we get:
\begin{eqnarray*}
\mu_M(f_1 f_2) &=& {1\over Z_M} \int_{G/\Ad} Z_M(t) \mu_M(t)(f_1 f_2)  \; dt
\\ 
&=& {1\over Z_M} \int_{G/\Ad} Z_{M_1}(t) \mu_{M_1}(t)(f_1) \; Z_{M_2}(t^{-1})
\mu_{M_2}(t^{-1})(f_2) \; dt.
\end{eqnarray*}
We rewrite this last equality in a more symmetric form:
$$Z_M\mu_M(f)=\int_{G/\Ad} Z_{M_1}(t) \mu_{M_1}(t)(f_1) \; Z_{M_2}(t^{-1})
\mu_{M_2}(t^{-1})(f_2) \; dt.$$

The point here is that the analytic objects that glue together in a simple way
are not the probability measures, but the measures with their natural weights.

Now, let us study further the relationships between the three
$\sigma$-algebras $\a_1$, $\a_2$ and $\a$. For this, it is convenient to
choose a base point $m$ on $L$, say $m=L(0)$ and to consider only loops based
at $m$. According to the proposition \ref{based_loops_gen_t}, the random
holonomies along these based loops generate the whole $\sigma$-algebra $\a$.
Thus, we consider the probability spaces $\Omega_{M_i}=\f(L_mM_i,G)/\Ad$ and
$\Omega_M=\f(L_mM,G)/\Ad$ endowed with the $\sigma$-algebras
$\a_i=\sigma(\h_{l_1,\ldots,l_n}, l\in L_mM_i)$, $i=1,2$ and
$\a=\sigma(\h_{l_1,\ldots,l_n}, l\in L_mM)$. 

The theorem \ref{Markov} says that the probability spaces
$(\Omega_{M_i},\a_i,\mu_{M_i}(t^{\pm 1}))$, $i=1,2$ are naturally isomorphic
with two independent subspaces of $(\Omega_M,\a,\mu_M(t))$. It is natural to
ask whether the inclusion $\a_1 \vee \a_2 \subset \a$ is an equality or not. 

In order to answer this question, choose $l_1 \in L_mM_1$ and $l_2\in L_mM_2$.
The random variable $\h_{l_1,L,l_2}$ is $\a$-measurable, but is it
$\a_1\vee\a_2$-measurable? The random variables $\h_{l_1,L}$ and $\h_{L,l_2}$
are respectively $\a_1$- and $\a_2$-measurable and it seems reasonable to
believe that they provide the whole information about $\h_{l_1,L,l_2}$
available in $\a_1\vee\a_2$.

If $G$ is Abelian, then $\h_{l_1,L}=(H_{l_1},H_L)$, $\h_{L,l_2}=(H_L,H_{l_2})$
and $\h_{l_1,L,l_2}=(H_{l_1},H_L,H_{l_2})$ so that $\h_{l_1,L,l_2}$ is
certainly $\a_1\vee\a_2$-measurable. We shall prove in this case that the
equality $(\a_1\vee \a_2)^\sim =\a$ holds, where the tilde denotes the
completion with respect to $\mu_M(t)$.

But if $G$ is not Abelian, this is not true anymore. Let us consider the
example $G=SO(3)$ and describe the conjugacy and joint conjugacy classes in
$SO(3)$ and $SO(3)^n$.

If $r\in SO(3)$ is neither the identity nor a symmetry, it has an angle and an
axis, which can be oriented in such a way that the angle is an element of
$(0,\pi)$. Thus, $r$ has a half-axis and an angle. This angle characterizes
the conjugacy class of $r$ (this is still true for the identity and the
symmetries). Now, consider $(r_1,\ldots,r_n)$ and $(r'_1,\ldots,r'_n)$ two
$n$-uples of rotations that have half-axes and angles $(u_i,\theta_i)$,
$(u'_j,\theta'_j)$. They belong to the same joint conjugacy class if and only
if $\theta_i=\theta'_i$ for all $i=1\ldots n$ and if there exists a rotation
$R\in SO(3)$ such that $u'_i=R(u_i)$.

\begin{figure}[h]
\begin{center}
\input{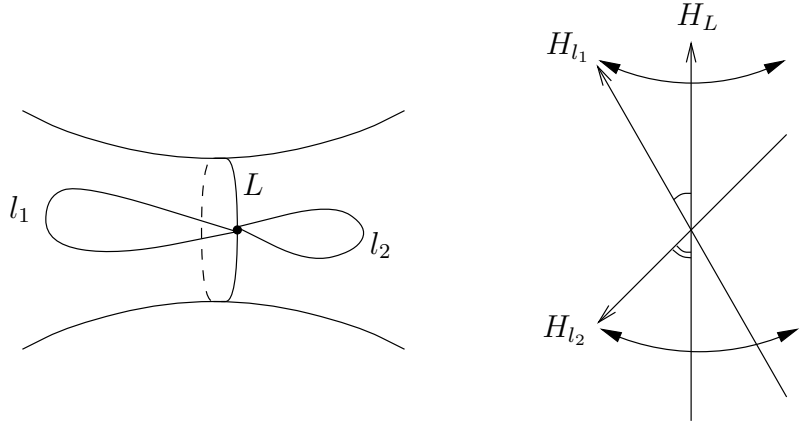}
\end{center}
\caption{Lack of information about the joint class $[H_L,H_{l_1},H_{l_2}]$
  when $G=SO(3)$.
\label{axes}}
\end{figure}

The random variables $\h_{l_1,L}$ and $\h_{L,l_2}$ determine the angles of the
three rotations $H_{l_1}, H_L, H_{l_2}$ and the relative position of the
half-axes of $H_{l_1}$ and $H_L$ on one hand and $H_L$ and $H_{l_2}$ on the
other hand. But this is not enough to determine the relative position of the
three half-axes (see figure \ref{axes}). There remains an undetermined
rotation around the axis of $H_L$. Such a rotation is an element of $SO(3)$
which commutes with $H_L$, in other words an element of the centralizer of
$H_L$. 

This informal argument suggests that $\a_1\vee \a_2$ is really a smaller
$\sigma$-algebra than $\a$. On the other hand, if the distributions of
$H_{l_1}$ and $H_{l_2}$ are not singular with respect to the Haar measure,
which is the case if $l_1$ and $l_2$ are not constant loops, it is possible to
extract from the random variable $\h_{l_1,L,l_2}$ the amount of information
that it is necessary to add to $\a_1\vee \a_2$ to recover $\a$.

Recall that $t$ is an element of $SO(3)/\Ad$ and that $\h_L=t$
$\mu_M(t)$-almost surely. Pick an element $x \in t$. Suppose that we are
given a measurable section $\tau : (t\times SO(3))/\Ad \lra \{x\} \times
SO(3)$ of the canonical projection. A realization of
$\h_{l_1,L,l_2}=[H_{l_1},H_L,H_{l_2}]$ determines $\mu_M(t)$-a.s. two classes
of $(t\times SO(3))/\Ad$, namely $[H_L,H_{l_1}]$ and $[H_L,H_{l_2}]$. Let
$(x,x_1)$ and $(x,x_2)$ be the images by $\tau$ of these classes. It is easy
to check that there exists an element $z$ of the centralizer $C(x)$ of $x$
such that $\h_{l_1,L,l_2}=[x_1,x,z^{-1}x_2 z]$ and that this $z$ is determined
up to right and left multiplication by an element of $C(x_1)\cap C(x_2)$,
which is almost surely equal to the center $Z(G)$ of $G$, so that the class
$zZ(G)\in C(x)/Z(G)$ is well defined. 

Thus, we outlined the construction of a random variable $Z$ with values in
$C(x)/Z(G)$, which can be seen as the missing angle of the figure \ref{axes}.
We will show that it has a uniform distribution, that it is independent of
$\a_1\vee \a_2$ and that $(\a_1\vee\a_2\vee\sigma(Z))^\sim=\a^\sim$.
Note that if $G$ is Abelian, the space $C(x)/Z(G)$ is just a point, so this
is not in contradiction with the fact that $(\a_1\vee\a_2)^\sim=\a^\sim$.

\section{Sewing of two surfaces along one circle}

We consider the general situation described at the beginning of the chapter,
with $p=1$, i.e. when the two surfaces are glued together along only one
circle $L$. Pick $x\in t$ and choose elements
$u_1$ and $u_2$ in $(G\Ad)^{p_1}$ and $(G/\Ad)^{p_2}$, representing the
holonomies along the other components of the boundaries of $M_1$ and $M_2$ and
maybe other loops inside $M_1$ and $M_2$. 

\begin{theorem} \label{tribuz1} 1. If $G$ is Abelian, then the completions of
  the $\sigma$-algebras $\a_1\vee \a_2$ and $\a$ with respect to the measure
  $\mu_M(u_1,t,u_2)$ are equal. \\
  \indent 2. If $G$ is not Abelian, there exists a non-trivial
  sub-$\sigma$-algebra $\z$ of $\a$ which is independent of $\a_1\vee \a_2$
  and such that the completions of $\a_1\vee\a_2\vee\z$ and $\a$ with respect
  to $\mu_M(u_1,t,u_2)$ are equal.
\end{theorem} 

By a non-trivial $\sigma$-algebra, we mean a $\sigma$-algebra which is not
independent of itself, or equivalently, that contains sets of probability
different from $0$ and $1$. 

The proof of this theorem requires a certain amount of preliminary work, that
leads to some interesting results. Let us choose on $M$ a Riemannian metric
such that $\partial M$ and $L$ are geodesic. Recall that $m=L(0)$. The set of
piecewise geodesic loops on $M$ based at $m$ is denoted by $PGL_mM$.

\begin{lemma}\label{gen12}
  On the space $\Omega_M$, the $\sigma$-algebra $\a$ is contained in
  $\a_{12}$, which is the completion of
  $\sigma(\h_{\lambda_1,\ldots,\lambda_n}, \lambda_1,\ldots,\lambda_n \in
  PGL_mM_1 \cup PGL_mM_2)$ with respect to $\mu_M(u_1,t,u_2)$. In particular,
  $\a^\sim=\a_{12}$. 
\end{lemma}

\pf. By the multiplicativity of the random holonomy described at the end of
section \ref{remarkable}, the holonomies along loops that are finite products
of loops of $L_mM_1$ and $L_mM_2$ are $\a_{12}$-measurable.  According to the
continuity property stated in \ref{reg_YM} and to the fact that the holonomy
depends only on the equivalence class of the loops, the result depends on the
fact that any loop of $L_mM$ can be approximated by loops that are equivalent
to finite products of loops of $PGL_mM_1$ and $PGL_mM_2$. Consider a piecewise
geodesic loop of $L_mM$. This loop cuts $L$ transversally at most a finite
number of times, hence it is equivalent to a finite product of loops of
$PGL_mM_1$ and $PGL_mM_2$.  Since any loop of $L_mM$ can be approximated by
piecewise geodesic loops, the result is proved. \qed

According to this lemma, we restrict our attention to the loops based at $m$
and whose image is either in $M_1$ or in $M_2$. This is enough to prove the
first part of the theorem: if $G$ is Abelian, a random variable of the form
$\h_{\lambda_1,\ldots,\lambda_n}$ with $\lambda_1,\ldots,\lambda_n \in
PGL_mM_1 \cup PGL_mM_2$ is just $(H_{\lambda_1},\ldots,H_{\lambda_n})$ and it
is certainly measurable with respect to $\a_1\vee\a_2$, since each variable
$H_{\lambda_i}$ is either $\a_1$- or $\a_2$-measurable.

If $G$ is not Abelian, we need to compute the conditional distribution of a
random variable like $\h_{\lambda_1,\ldots,\lambda_n}$ with respect to
$\a_1\vee\a_2$. For this, it is convenient to introduce the following
notation. 

Recall that $x$ is a fixed element in $t\in G/\Ad$. Consider an element of
$(t\times G^n)/\Ad$. It is a class that can be written under the form
$[x,x_1,\ldots,x_n]$, where the $n$-uple $(x_1,\ldots,x_n)$ is defined up to
conjugation by an element of the centralizer $C(x)$. The measure
$$\int_{C(x)} \delta_{(z^{-1}x_1z,\ldots,z^{-1}x_n z)} \; dz $$
is well
defined on $G^n$, depending only on the class $[x,x_1,\ldots,x_n]$ and on the
choice of $x$. We denote it by $\pi_x^n([x,x_1,\ldots,x_n])$. Note that if
$x'=y^{-1}xy$, then for any function $f$ on $G^n$,
\begin{equation}
\label{transfo_pi}
\int_{G^n} f \; d\pi_{x'}^n([x,x_1,\ldots,x_n]) = \int_{G^n} f\circ \Ad(y) \;
d\pi_x^n([x,x_1,\ldots,x_n]).
\end{equation}

\begin{proposition} \label{conda1a2}
  Let $\lambda_1,\ldots,\lambda_n$ be $n$ loops of $PGL_mM_1$ and
  $\lambda'_1,\ldots,\lambda'_{n'}$ be $n'$ loops of $PGL_mM_2$. For any
  continuous function $f$ on $G^{n+n'}$ invariant by diagonal adjunction, we
  have
  \begin{eqnarray*}  
  E_{\mu_M(u_1,t,u_2)}[f(\h_{\lambda_1,\ldots,\lambda_n,\lambda'_1,\ldots,
    \lambda'_{n'}})| \a_1\vee\a_2] &=& \\
  && \hskip -2.3cm \int_{G^{n+n'}} f \;
  d\pi_x^n(\h_{L,\lambda_1,\ldots,\lambda_n}) \otimes
  d\pi_x^{n'}(\h_{L,\lambda'_1,\ldots,\lambda'_{n'}}) \;\; {\rm a.s.}
  \end{eqnarray*}
\end{proposition}

The right hand side does not depend on $x$, only on $t$, because of the
invariance of $f$ and the relation (\ref{transfo_pi}). \\

\pf. Choose $l_1,\ldots,l_m \in PGL_mM_1$ and $l'_1,\ldots,l'_{m'} \in
PGL_mM_2$. Let $f_1$ be a continuous function on $G^m/\Ad \times G^{m'}/\Ad$.
Let $\Gamma$ be a graph on $M$ such that $L$, the $l_i$'s, the $l'_i$'s, the
$\lambda_i$'s and the $\lambda'_i$'s are in $\Gamma^*$. As in the proof of the
theorem \ref{Markov}, denote by $\Gamma_1$ and $\Gamma_2$ the graphs induced
by $\Gamma$ on $M_1$ and $M_2$.  Let $f$ be a continuous function on
$G^{n+n'}$ invariant by diagonal adjunction. Using the proposition
\ref{discMarkov}, we get
\begin{eqnarray*}
&&
E_{\mu_M(u_1,t,u_2)}[f(\h_{\lambda_1,\ldots,\lambda_n,\lambda'_1,\ldots,\lambda'_{n'}}) 
f_1(\h_{l_1,\ldots,l_m},\h_{l'_1,\ldots,l'_{m'}})] = \\
&& \hskip -.7cm =\int_{G^\Gamma} f(h_{\lambda_1},\ldots,h_{\lambda_n},
h_{\lambda'_1},\ldots,h_{\lambda'_{n'}}) f_1([h_{l_1},\ldots,h_{l_m}],
[h_{l'_1},\ldots,h_{l'_{m'}}]) \; dP^{\Gamma}(y,x,y') \\
&& \hskip -.7cm =\int_{G^{\Gamma_1}\times G^{\Gamma_2}}
f(h_{\lambda_1},\ldots,h_{\lambda_n}, 
h_{\lambda'_1},\ldots,h_{\lambda'_{n'}}) f_1([h_{l_1},\ldots,h_{l_m}],
[h_{l'_1},\ldots,h_{l'_{m'}}]) \\
&& \hskip 7.9cm dP^{\Gamma_1}(y,x) dP^{\Gamma_2}(x^{-1},y').
\end{eqnarray*}
For each element $z\in C(x)$ and each $i=1,2$, the gauge transformation equal
to $1$ at each vertex of $\Gamma_i$ except at $m$ where it takes the value $z$
leaves $P^{\Gamma_i}$ invariant. Thus, the last integral is equal to
\begin{eqnarray*}
&& \int_{G^{\Gamma_1}\times G^{\Gamma_2}}
\int_{C(x)^2} f(z_1^{-1}h_{\lambda_1}z_1,
\ldots,z_1^{-1}h_{\lambda_n}z_1,  
z_2^{-1} h_{\lambda'_1} z_2,\ldots,z_2^{-1} h_{\lambda'_{n'}}z_2) \; dz_1
dz_2\\ 
&& \hskip 3cm f_1([h_{l_1},\ldots,h_{l_m}],[h_{l'_1},\ldots,h_{l'_{m'}}])
\; dP^{\Gamma_1}(y,x) dP^{\Gamma_2}(x^{-1},y') \\
&=& \int_{G^{\Gamma_1}\times G^{\Gamma_2}}
\int_{G^{n+n'}} f \;
d\pi_x^n([h_L,h_{\lambda_1},\ldots,h_{\lambda_n}]) \otimes
d\pi_x^{n'}([h_L,h_{\lambda'_1},\ldots,h_{\lambda'_{n'}}]) \\
&& \hskip 3cm f_1([h_{l_1},\ldots,h_{l_m}],[h_{l'_1},\ldots,h_{l'_{m'}}]) 
\; dP^{\Gamma_1}(y,x) dP^{\Gamma_2}(x^{-1},y') \\
&=& E\left[\int_{G^{n+n'}} f \; 
d\pi_x^n([\h_{L,\lambda_1,\ldots,\lambda_n}]) \otimes
d\pi_x^{n'}([\h_{L,\lambda'_1,\ldots,\lambda'_{n'}}]) f_1(\h_{l_1,\ldots,l_m},\h_{l'_1,\ldots,l'_{m'}})\right],
\end{eqnarray*}
which proves the proposition. \qed

Remark that the ``reasonable'' intuition that the variables
$\h_{L,\lambda_1,\ldots,\lambda_n}$ and
$\h_{L,\lambda'_1,\ldots,\lambda'_{n'}}$ contain all the information about
$\h_{\lambda_1,\ldots,\lambda_n,\lambda'_1,\ldots,\lambda'_{n'}}$ available in
$\a_1\vee\a_2$ is confirmed.

The next step is the construction of the variable $Z$ in a general context. We
use the same method as in the example treated informally in section
\ref{examp} but this requires some technical results. 

\begin{lemma} \label{SN} Let $G$ be a compact connected Lie group of dimension
  $n$ and rank $k$. Set $N=n-k+1$. The set $S_N$ of all $N$-uples
  $(g_1,\ldots,g_N)$ such that the closed subgroup of $G$ generated by
  $g_1,\ldots,g_N$ is $G$ itself has full Haar measure in $G^N$. Moreover,
  $S_N$ is stable by diagonal adjunction, i.e. $(g_1,\ldots,g_N)$ belongs to
  $S_N$ if and only if $(h^{-1}g_1,h,\ldots,h^{-1}g_Nh)$ does for any $h$ in
  $G$.
\end{lemma}

\pf. The key of this result is that almost every element of $G$ is regular,
i.e. generates a maximal torus \cite{Simon}. Let $g_1$ be such a regular
element. It generates a subgroup $G_1=T_1$ of $G$ of dimension $k$. If $G_1$
is a proper subgroup of $G$, then $\dim G_1 < \dim G$ since $G$ is connected.
Thus, the complementary of $G_1$ has full measure as well as the set of
regular elements outside $G_1$. Let $g_2$ be such an element and $T_2$ the
torus that it generates. Denote by $G_2$ the subgroup generated by
$\{g_1,g_2\}$. Denote by {\got{g}}$_1$, {\got{g}}$_2$ and {\got{t}}$_2$ the
Lie algebras of $G_1$, $G_2$ and $T_2$ respectively. The fact that
{\got{t}}$_2 \not\subset$ {\got{g}}$_1$ and {\got{g}}$_2 \supset $
{\got{g}}$_1 + ${\got{t}}$_2$ shows that $\dim G_2 > \dim G_1$. Repeating this
procedure $N=n-k+1$ times, we get a subgroup $G_N$ which is equal to $G$. It
is clear from this construction that the set $S_N$ of convenient $N$-uples has
full Haar measure in $G^N$. 

The last statement depends on the fact that two conjugate $N$-uples generate
two conjugate subgroups of $G$. \qed

In particular, note that if $(g_1,\ldots,g_N)\in S_N$, then $C(g_1)\cap \ldots
\cap C(g_N)=Z(G)$.

Before to prove the next lemma, let us recall that, by a classical result of
Kuratowski \cite{Kuratowski}, a one-to-one measurable map between two Polish
spaces sends Borel subsets to Borel subsets.

\begin{lemma} \label{meas-tau}  There exists a measurable section $\tau:
  (t\times G^N)/\Ad \lra \{x\} \times G^N$.
\end{lemma}

\pf. We use a theorem of Bourbaki ({\it Topologie}, chap. IX, $\S$ 6, N$^o$9,
th.  5) \cite{Bourbaki_Topologie_9}, which says that there exists a Borel
subset $R$ of $t\times G^N$ that meets once and only once each orbit of the
action of $G$. This subset $R$ allows to define a section $\tau$. To prove
that $\tau$ is measurable, consider a Borel subset $B \subset t\times G^N$.
The fact that $\tau^{-1}(B)=p(B\cap R)$, where $p: t\times G^N \lra (t\times
G^N)/G$ is the natural projection, together with the result of Kuratowski
mentionned above, shows that $\tau^{-1}(B)$ is a Borel subset. Thus, $\tau$ is
measurable. \qed

Let us choose now $2N$ loops , $L^1_1,\ldots,L^1_N \in PGL_mM_1$ and
$L^2_1,\ldots,L^2_N \in PGL_mM_2$, such that the distribution of
$(H_{L^1_1},\ldots,H_{L^1_N},H_{L^2_1},\ldots,H_{L^2_N})$ has a density with
respect to the Haar measure on $G^{2N}$. This is a weak condition on the
$L^i_j$'s which ensures that the set where $(H_{L^i_1},\ldots,H_{L^i_N})$
belongs to $S_N$, $i=1,2$ has full Yang-Mills measure. Let us fix $\omega$
inside this set. Set
$\tau(\h_{L,L^1_1,\ldots,L^1_N}(\omega))=(x,x^1_1,\ldots,x^1_N)$ and
$\tau(\h_{L,L^2_1,\ldots,L^2_N}(\omega))=(x,x^2_1,\ldots,x^2_N)$, where $\tau$
is a measurable section given by the preceding lemma. Choose an element
$(y^1_1,\ldots,y^1_N,x,y^2_1,\ldots,y^2_N)$ inside the class
$\h_{L^1_1,\ldots,L^1_N,L,L^2_1,\ldots,L^2_N}(\omega)$. For $i=1,2$, the
equality $[x,x^i_1,\ldots,x^i_N]=[x,y^i_1,\ldots,y^i_N]$ implies the existence
of $z_i \in C(x)$ such that $y^i_1=z_i^{-1}x^i_1 z_i, \ldots,
y^i_N=z_i^{-1}x^i_N z_i.$ Thus there exists $z\in C(x)$ such that
$$[y^1_1,\ldots,y^1_N,x,y^2_1,\ldots,y^2_N]=[x^1_1,\ldots,x^1_N,x,
z^{-1}x^2_1z,\ldots,z^{-1}x^2_Nz],$$ for example $z=z_2 z_1^{-1}$. If another
element $z'\in C(x)$ satisfies the same relation, then it is easily checked
that there exist two elements $z'' \in C(x^1_1)\cap \ldots \cap
C(x^1_N)=Z(G)$ and $z'''\in  C(x^2_1)\cap \ldots \cap
C(x^2_N)=Z(G)$ such that $z'=z''z z'''=z z''z'''$. Thus, the class $zZ(G) \in
C(x)/Z(G)$ is well defined and we define the value at $\omega$ of the random
variable $Z$ by 
$$Z(\omega)=z Z(G).$$

Note that $Z$ depends on $\omega$ only through
$\h_{L^1_1,\ldots,L^1_N,L,L^2_1,\ldots,L^2_N}(\omega)$ so that we may write
$Z=\widetilde Z(\h_{L^1_1,\ldots\ldots,L^2_N})$ for some function $\widetilde
Z : G^{2N+1}/\Ad \lra C(x)/Z(G)$. Moreover, for any $z_1,z_2$ in $C(x)$ and
any $x^1,x^2$ in $G^N$, it is easily checked that
$$\widetilde Z([z_1^{-1}x^1z_1,x,z_2^{-1}x^2z_2])=z_1z_2^{-1}\widetilde
Z([x^1,x,x^2]).$$

\begin{proposition} 
The random variable $Z$ is uniformly distributed on the Lie group
  $C(x)/Z(G)$ and independent of $\a_1\vee\a_2$ under the measure
  $\mu_M(u_1,t,u_2)$. 
\end{proposition}

Since $Z(G)$ is a closed normal subgroup of $C(x)$, the quotient $C(x)/Z(G)$
is still a Lie group, so that it makes sense to speak about uniform
distribution. \\

\pf. Let $f$ be a continuous function on $C(x)/Z(G)$. It is more convenient
for the notations to consider $\widetilde Z$ as a function of the variable
$\h_{L^1_1,\ldots,L^1_N,L,L,L^2_1,\ldots,L^2_N}$, with twice the loop $L$. 
Using the proposition \ref{conda1a2}, we find
\begin{eqnarray*}
E[f(Z)|\a_1\vee\a_2] &=& E_{\mu_M(u_1,t,u_2)}[f(\widetilde
Z(\h_{L^1_1,\ldots,L^1_N,L,L,L^2_1,\ldots,L^2_N}))|\a_1\vee\a_2] \\
&=& \int_{G^{2n+2}} f\circ \widetilde Z \; d\pi_x(\h_{L,L^1_1,\ldots,L^1_N})
\otimes \delta_{(x,x)} \otimes \pi_x(\h_{L,L^2_1,\ldots,L^2_N}).
\end{eqnarray*}

Let us fix $\omega$ and set $\tau(\h_{L,L^1_1,\ldots,L^1_N}(\omega))=(x,x^1)$,
$\tau(\h_{L,L^2_1,\ldots,L^2_N}(\omega))=(x,x^2)$, where $x^1$ and $x^2$ are
elements of $G^N$. Then
\begin{eqnarray*}
E[f(Z)|\a_1\vee\a_2](\omega) &=& \int_{C(x)^2} f(\widetilde Z(z_1^{-1}x^1
z_1,x,x,z_2^{-1}x^2 z_2))\; dz_1 dz_2 \\
&=& \int_{C(x)^2} f(z_1 z_2^{-1} \widetilde Z(\h_{L^1_1,\ldots,L^1_N,L,
  L,L^2_1,\ldots,L^2_N})) \; dz_1 dz_2 \\
&=& \int_{C(x)} f(z) \; dz.
\end{eqnarray*}

Since this conditional expectation is a constant for any $f$, $Z$ is
independent of $\a_1 \vee \a_2$. \qed

The proof of the theorem is almost finished, there only remains to prove that
$\a_1\vee\a_2\vee\sigma(Z)$ contains $\a_{12}$. 

\begin{lemma}
The random variable $\h_{L^1_1,\ldots,L^1_N,L,L^2_1,\ldots,L^2_N}$ is
measurable with respect to the $\sigma$-algebra $\a_1\vee\a_2\vee\sigma(Z)$.
\end{lemma}

\pf. Let us denote by $\tau' : (t\times G^N)/\Ad \lra G^N$ the composition of
$\tau$ with the projection on the $N$ last factors. Then, by construction of
$Z$, the following equality holds almost surely:
$$\h_{L^1_1,\ldots,L^1_N,L,L^2_1,\ldots,L^2_N}
=[\tau'(\h_{L,L^1_1,\ldots,L^1_N}),x,Z^{-1}
\tau'(\h_{L,L^2_1,\ldots,L^2_N})Z],$$ 
proving the lemma. \qed

The next proposition will finish the proof:

\begin{proposition} \label{generate}
  The $\sigma$-algebra
  $\a_1\vee\a_2\vee\sigma(\h_{L^1_1,\ldots,L^1_N,L,L^2_1,\ldots,L^2_N})$
  contains $\a_{12}$.
\end{proposition}

\pf.  Pick $p$ loops $l_1,\ldots,l_p$ in $L_mM_1$ and $q$ loops
$l'_1,\ldots,l'_q$ in $L_mM_2$. We abbreviate by $l,l',L^1$ and $L^2$ the
corresponding families of loops.

We will show that the random variables $\h_{l,L^1}$, $\h_{L^2,l'}$ and
$\h_{L^1,L^2}$ determine $\h_{l,l'}$ if the values of $\h_{L^1}$ and
$\h_{L^2}$ are in $S_N/\Ad$. We do not even restrict ourselves to piecewise
geodesic loops. Since $\h_{l,l'}$ is just a continuous projection of
$\h_{l,L^1,L^2,l'}$, it is sufficient to write this last variable as a
function of the three given variables. For this, we construct a map:
$$(G^p\times S_N)/\Ad \times_{S_N} (S_N \times S_N)/\Ad \times_{S_N} (S_N
\times G^q)/\Ad \build{\lra}_{}^{\kappa} (G^p\times S_N \times S_N \times
G^q)/\Ad .$$
The symbols $\times_{S_N}$ above mean that the map $\kappa$ is only
defined on the set of elements of the form
$([u_1,\ldots,u_p,g_1,\ldots,g_N],[g'_1,\ldots,g'_N,h'_1,\ldots,h'_N],
[h_1,\ldots,h_N,v_1,\ldots,v_q])$  
such that $[g_1,\ldots,g_N]=[g'_1,\ldots,g'_N]$ and
$[h_1,\ldots,h_N]=[h'_1,\ldots,h'_N]$. 

We claim that it makes sense to construct $\kappa$ such that the image of such
a triple is the unique element
$[x_1,\ldots,x_p,r_1,\ldots,r_N,s_1,\ldots,s_N,y_1,\ldots,y_q]$ such that,
with compact notations, $[x,r]=[u,g]$, $[r,s]=[g',h']$ and $[s,y]=[h,v]$.  It
is not difficult to see that such an element exists: if $z_g$ and $z_h$ are
two elements of $G$ such that $\Ad(z_g)g=g'$ and $\Ad(z_h)h=h'$, then
$[\Ad(z_g) u, g',h',\Ad(z_h)v]$ is a possible choice.

Suppose that $[x,r,s,y]$ and $[x',r',s',y']$ are two candidates. Since
$[x,r]=[x',r']$, there exists $z_r\in G$ such that $\Ad(z_r) x=x'$ and
$\Ad(z_r)r =r'$. Similarly, there exists $z_s$ such that $\Ad(z_s)y=y'$ and
$\Ad(z_s)s=s'$. Now, $[r,s]=[r',s']$ implies $[r,s]=[\Ad(z_r) r, \Ad(z_s) s]=
[\Ad(z_r z_s^{-1})r,s]$. This forces $z_r z_s^{-1}$ to be an element of
$Z(G)$, so that $\Ad(z_r)=\Ad(z_s)$ and the two candidates are equal.

We have $\mu_M(t)$-almost surely
$\h_{l,L^1,L^2,l'}=\kappa(\h_{l,L^1},\h_{L^1,L^2},\h_{L^2,l'})$, so that it
remains only to prove that $\kappa$ is measurable. To see this, remark that
$\kappa^{-1}$ is easier to define than $\kappa$: it is a restriction of three
continuous projections defined on $(G^p\times G^N \times G^N \times G^q)/\Ad$.
Moreover, since $\kappa$ is well defined, $\kappa^{-1}$ is injective. Thus, by
the result of Kuratowski mentionned above, $\kappa^{-1}$ sends Borel subsets
to Borel subsets, which means exactly that $\kappa$ is measurable. \qed


\section{Cutting and gluing handles}
\label{cut_2}

An analysis very similar to that of the three preceding sections can be done
in the case of a surface $M$ obtained by sewing together two components of the
boundary of a surface $M_1$. There is no more conditional independence in this
situation, but it is still possible to study the relationships between the
measures $\mu_{M_1}$ and $\mu_M$.

Let $M_1$ be a surface with boundary such that $\partial M_1$ has at least two
components $N_1$ and $N_2$. Let $\psi:N_1 \lra N_2$ be an
orientation-reversing diffeomorphism and let $M$ be obtained by gluing $M_1$
along $\psi$. Let $L_1$ and $L_2$ be two loops on $M_1$ whose images are $N_1$
and $-N_2$ respectively, and call $L$ the corresponding loop on $M$. Set
$m=L(0)$, $m_1=L_1(0)$ and $m_2=L_2(0)$. Note that, in contrast to the
preceding situation, $M_1$ is not embedded in $M$, it is only immersed.
Nevertheless, this immersion allows us to map $LM_1$ into $LM$. So, a function
$f$ on $\f(LM_1,G)$ can be seen as a function $\widetilde f$ on $\f(LM,G)$ and
$f$ is measurable with respect to $\a_1$ if and only if $\widetilde f$ is
measurable with respect to $\widetilde \a_1=\sigma(\h_{l_1,\ldots,l_n},l\in
LM_1)$. As before, we identify $f$ and $\a_1$ with $\widetilde f$ and
$\widetilde \a_1$.

\begin{theorem} \label{cutting_law}
  Let $f$ be a measurable function on $(\f(LM_1,G),\a_1)$. Then for any
  $u_1,\ldots,u_p,t \in G/\Ad$, the following equality holds:
  $$\mu_M(u_1,\ldots,u_p,t)(f)=\mu_{M_1}(u_1,\ldots,u_p,t,t^{-1})(f),$$
  where
  the $u_p$ are elements of $G/\Ad$ that stand for the value of the holonomy
  along other boundary components of $\partial M_1$ and maybe some other loops
  inside $M_1$.
\end{theorem}

As for the proof of the theorem \ref{Markov}, we begin with a discrete result.
Let $i$ denote the immersion of $M_1$ into $M$. Given a graph $\Gamma$ on $M$,
there exists a graph $\Gamma_1$ on $M_1$ such that $\Gamma_1$ is mapped by $i$
to $\Gamma$ and such that the edges of $\Gamma_1$ lying on $N_1$ and
$N_2$ respectively are in one-to-one correspondence  via the diffeomorphism
$\psi$. We call $\Gamma_1$ the graph on $M_1$ induced by $\Gamma$. Recall that
$m_1=L_1(0)$ and $m_2=L_2(0)$.

\begin{proposition} \label{disccutting}
  Let $\Gamma$ be a graph on $M$ and $\Gamma_1$ the graph induced on $M_1$ by
  $\Gamma$. Let $f$ be a continuous function on $G^{\Gamma_1}$ which is
  invariant by the gauge transformations $j_1$ such that
  $j_1(m_1)=j_1(m_2)=1$.  Then $f$ gives rise to a function on $G^\Gamma$
  invariant by the gauge transformations $j$ such that $j(m)=1$ and for any
  $x\in G$,
  $$\int_{G^{\Gamma_1}} f \; dP^{\Gamma_1}(y,x,x^{-1}) = \int_{G^\Gamma} f \;
  dP^\Gamma(y,x),$$
  where $x$ is the holonomy along $N_1$ and $-N_2$ and $y$
  stands for the values of the holonomies along the other components of
  $\partial M_1$ and maybe also along other loops inside $M_1$.
\end{proposition}

The proof of this proposition is simlar to the proof of \ref{discMarkov}, with
a small difference due to the consideration of particular gauge
transformations. Note that \ref{disccutting} implies the corresponding
statement without any restrictions on the gauge transformations, which is in
fact enough to prove the theorem \ref{cutting_law}. We will need this
refinement in a forthcoming proof. \\

\pf \textsc{ of} \ref{disccutting}. We have
\begin{eqnarray*}
\int_{G^{\Gamma_1}} f \; dP^{\Gamma_1}(y,x,x^{-1}) &=& {1\over
  {Z_{M_1}(y,x,x^{-1})}} \int_{G^{\Gamma_1}} f D^{\Gamma_1} \; d\nu_y d\nu_x
  d\nu_{x^{-1}} dg'\\
&=& {1\over{Z_{M_1}(y,x,x^{-1})}} \int_{G^{\Gamma_{1,\partial}}} d\nu_x
  d\nu_{x^{-1}} \int_{G^{\Gamma_1\backslash\Gamma_{1,\partial}}} f
  D^{\Gamma_1} \; d\nu_y dg'.
\end{eqnarray*}
The last integral is a function on $G^{\Gamma_{1,\partial}}$ and the same
argument as in the proof of \ref{discMarkov} together with the fact that we
consider only gauge transformations such that $j(m_1)=j(m_2)=1$ shows that
this function depends only on the values of $h_{L_1}$ and $h_{L_2}$. Thus, we
can drop the integration against $d\nu_{x^{-1}}$ and replace the variables
associated with the edges lying on $N_2$ by the variables associated to the
corresponding edges on $N_1$. Thus,
\begin{eqnarray*}
\int_{G^{\Gamma_1}} f \; dP^{\Gamma_1}(y,x,x^{-1}) &=& 
{1\over{Z_{M_1}(y,x,x^{-1})}} \int_{G^{\Gamma_\partial}} d\nu_x
\int_{G^{\Gamma \backslash \Gamma_\partial}} f D^\Gamma \; d\nu_y dg' \\
&=& {{Z_M(y,x)}\over {Z_{M_1}(y,x,x^{-1})}} \int_{G^\Gamma} f \;
dP^\Gamma(y,x).
\end{eqnarray*}
The equality of the normalization constants follows from the case $f=1$ and
finishes the proof. \qed

Once again, we state separately the property of the conditional partition
functions that we just established:

\begin{proposition} \label{Z_2}
  For any $t_1,\ldots,t_p,t \in G/\Ad$, the following equality holds:
  $$Z_M(t_1,\ldots,t_p,t,t^{-1})=Z_{M_1}(t_1,\ldots,t_p,t).$$
\end{proposition}

\pf \textsc{ of} \ref{cutting_law}. The proof is similar to that of the
theorem \ref{Markov}, but simpler because there is essentially nothing more to
prove than what is already stated in the proposition \ref{disccutting}. By
continuity of the random holonomy, it is enough to consider the case of a
cylindrical function $f$ that depends on the holonomy along loops that can be
put into a graph. So we choose such a function and consider a graph $\Gamma$
on $M$ such that $\Gamma^*$ contains all the loops with which we are working.
Let $\Gamma_1$ be the graph on $M_1$ induced by $\Gamma$. The function $f$
induces a gauge-invariant function $f$ on $G^{\Gamma_1}$, and the equality
which we need to prove about $f$ is exactly that given by the proposition
\ref{disccutting}. \qed

As in the preceding situation, we want to know more about the inclusion of
$\sigma$-algebras $\a_1 \subset \a$ and for this it is technically convenient
to consider based loops: set $\Omega_M=\f(L_mM,G)/\Ad$ and
$\Omega_{M_1}=\f(L_mM_1,G)/\Ad$, where $m=L(0)$. We identify $\a_1$ with the
$\sigma$-algebra on $\Omega_M$ generated by the random holonomies along the
loops of $i(L_mM)$.

Let us begin by a short informal discussion. Choose a path $c$ from $m_1$ to
$m_2$. Although $c$ is an open path in $M_1$, $i(c)$ is a loop on $M$, so that
the holonomy along $i(c)$ is well defined on $\Omega_M$ but certainly not
$\a_1$-measurable. Once again, there is some information about $\h_c$ in
$\a_1$ but there is also something missing, which is closely related to the
centralizer of $H_L$. To see this, consider the random variable
$\h_{L_1,cL_2c^{-1}}$. Since the holonomies along $L_1$ and $L_2$ must be
equal after the gluing procedure, we choose arbitrarily an element $x\in t$
and compute as if $H_{L_1}=H_{L_2}=x$. The variable $\h_{L_1,cL_2c^{-1}}$,
which is $\a_1$-measurable, determines to some extent the holonomy along $c$:
we interpret its value as that of $[x,H_{c}^{-1} x H_c]$ and this determines
the value of $H_c$ up to left and right multiplications by elements of the
centralizer $C(x)$.

\begin{figure}[h]
\begin{center}
\input{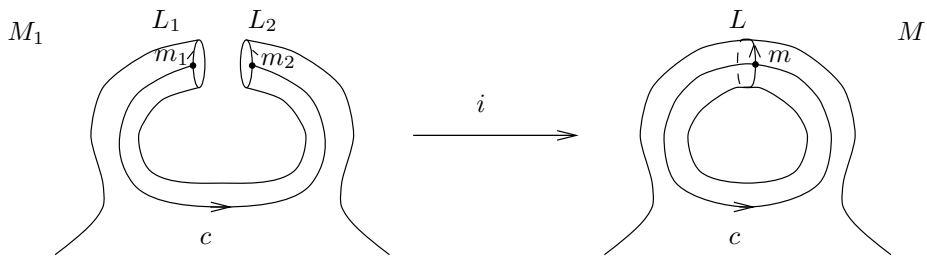}
\end{center}
\caption{The new handle of $M$.
\label{hand}}
\end{figure}

The main result is similar to the theorem \ref{tribuz1}:

\begin{theorem} \label{tribuz2} There exists a non-trivial
  sub-$\sigma$-algebra $\z'$ in $\a$ which is independent of $\a_1$ and such
  that the completions of $\a_1\vee\z'$ and $\a$ with respect to $\mu_M(u,t)$
  are equal.
\end{theorem}

Note that the conclusion of this theorem does not depend on the fact that $G$
is Abelian or not. We will even see that, in some vague sense, $\z'$ is the
biggest when $G$ is Abelian.

The proof of this theorem has the same structure as the proof of
\ref{tribuz1}. We begin with a lemma which allows us to consider a restricted
class of loops. We endow $M$ with a Riemannian metric such that $L$ is
geodesic. This metric on $M$ induces by pull-back by $i$ a metric on $M_1$.
We fix once for all a piecewise geodesic path $c$ in $M$ joining $m_1$ to
$m_2$ and meeting $\partial M_1$ only at $m_1$ and $m_2$.

\begin{lemma} On the space $\Omega_M$, the $\sigma$-algebra $\a$ is contained
  in the completion of the $\sigma$-algebra
  $\sigma(\h_{\lambda_1,\ldots,\lambda_n}, \lambda_1,\ldots,\lambda_n \in
  i(PGL_{m_1}M_1 \cup \{c\}))$ with respect to the measure $\mu_M(u,t)$.
\end{lemma}

\pf. Just as in the proof of \ref{gen12}, the point is to prove that any loop
of $L_mM$ can be approximated by loops that are equivalent to finite products
of loops of $i(PGL_{m_1}M_1\cup \{c\})$. Consider a piecewise geodesic loop.
Since it cuts at most a finite number of times $L$ transversally, it is
equivalent to a finite product of loops of $PGL_mM$ that are the images by $i$
of loops based at $m_1$ or at $m_2$ or of paths with endpoints $m_1$ and
$m_2$. Conjugation by $c$ or left or right multiplication by $c$ transform all
these paths on $M_1$ into loops based at $m_1$, hence the holonomy along our
piecewise geodesic loop can be expressed in terms of the holonomies along the
path of $i(PGL_{m_1}M_1 \cup \{c\})$. Since any loop of $L_mM$ is a limit of
piecewise geodesic loops, the result is proved. \qed

Now, let us compute the conditional expectation of these variables with
respect to $\a_1$. For this, we introduce a new notation, similar to the
notation $\pi^n_x$ introduced in the preceding section.

Remark that the mapping $PGL_{m_1}M_1 \cup \{c\} \lra PGL_mM$ induced by $i$
is injective, so that we regard $PGL_{m_1}M_1 \cup \{c\}$ as a subset of
$PGL_mM$.

Recall that $x$ is a fixed element in $t$ and consider an element of $(t\times
t\times G^n)/\Ad$. It is a joint conjugacy class that can be written
$[x,x_c^{-1}xx_c,x_1,\ldots,x_n]$. Let us compare two such representations of
this class. If
$[x,x_c^{-1}xx_c,x_1,\ldots,x_n]=[x,x_c'^{-1}xx_c',x_1',\ldots,x_n']$, then
there exists $z_1\in C(x)$ such that $x'_i=z_1^{-1}x_iz_1$ for $i=1,\ldots,n$
and $x_c'^{-1}xx_c'=z_1^{-1}x_c^{-1}xx_cz_1$. This last equality implies that
$x_cz_1x_c'^{-1}$ is an element $z_2$ of $C(x)$, so that $x_c'=z_2^{-1}
x_cz_1$. Thus, for any function $f$ on $G^{n+1}/\Ad$, the following integral
is well defined:
$$\int_{C(x)^2} f(z_1^{-1}x_1z_1,\ldots,z_1^{-1}x_nz_1,z_2^{-1}x_cz_1)\; dz_1
dz_2 =\int_{C(x)} f(x_1,\ldots,x_n,zx_c)\; dz.$$
We denote by
$\alpha^{n+1}_x([x,x_c^{-1}xx_c,x_1,\ldots,x_n])$ the corresponding measure on
$G^{n+1}/\Ad$.

\begin{proposition} \label{conda1} Let $\lambda_1,\ldots,\lambda_n$ be loops
  of $PGL_{m_1}M_1$. For any continuous function $f$ on $G^{n+1}/Ad$, we have
  $$E[f(\h_{\lambda_1,\ldots,\lambda_n,c})|\a_1]=\int_{G^{n+1}} f \;
  d\alpha^{n+1}_x(\h_{L,cLc^{-1},\lambda_1,\ldots,\lambda_n}).$$
\end{proposition}

\pf. Choose $l_1,\ldots,l_m$ in $PGL_{m_1}M_1$ and a continuous function $f_1$
on $G^m/\Ad$. Let $\Gamma$ be a graph on $M$ such that all the loops we are
considering belong to $\Gamma^*$. Denote by $\Gamma_1$ the graph induced by
$\Gamma$ on $M_1$, in the same way as we did in the proof of the theorem
\ref{cutting_law}. 
\begin{eqnarray*}
&&\hskip -.5cm E[f(\h_{\lambda_1,\ldots,\lambda_n,c})
f_1(\h_{l_1,\ldots,l_m})] = 
\int_{G^\Gamma} f(h_{\lambda_1},\ldots,h_{\lambda_n},h_c)
f_1(h_{l_1},\ldots,h_{l_m}) \; dP(y,x) \\
&&\hskip -.2cm = \int_{G^{\Gamma_1}} f(h_{\lambda_1},\ldots,h_{\lambda_n},h_c)
f_1(h_{l_1},\ldots,h_{l_m}) \; dP^{\Gamma_1}(y,x,x^{-1}),
\end{eqnarray*}
thenks to \ref{disccutting}. Indeed, the integrand here is not gauge-invariant
because of the $h_c$, but it is precisely invariant by those gauge
transformations $j$ such that $j(m_1)=j(m_2)=0$. 

The end points of $h_c$ are now different and we can use the
invariance of $P^{\Gamma_1}(y,x,x^{-1})$ under a gauge transformation
identically equal to $1$, except at $m_2$ where it takes the value $z^{-1} \in
C(x)$. This transformation leaves the $h_{\lambda_i}$'s and $h_{l_i}$'s
invariant since they are based at $m_1$, and transforms $h_c$ into
$zh_c$. Thus,
\begin{eqnarray*}
&& E[f(\h_{\lambda_1,\ldots,\lambda_n,c})
f_1(\h_{l_1,\ldots,l_m})] = \\
&&\hskip -.8cm  =\int_{G^{\Gamma_1}}
\left[ \int_{C(x)} f(h_{\lambda_1},\ldots,h_{\lambda_n},z h_c) \; dz
\right] 
f_1(h_{l_1},\ldots,h_{l_m}) \; dP^{\Gamma_1}(y,x,x^{-1})\\
&&\hskip -.8cm  =\int_{G^{\Gamma_1}}
\left[\int_{G^{n+1}} f \;
  d\alpha^{n+1}_x([h_L,h_{cLc^{-1}},h_{\lambda_1},\ldots,h_{\lambda_n}])
\right] f_1(h_{l_1},\ldots,h_{l_m}) \; dP^{\Gamma_1}(y,x,x^{-1})\\
&& \hskip -.8cm =E\left[\left(\int_{G^{n+1}} f \;
  d\alpha^{n+1}_x(\h_{L,cLc^{-1},\lambda_1,\ldots,\lambda_n})\right) f_1(\h_{l_1,\ldots,l_m}) \right],
\end{eqnarray*}
which implies the result. \qed

It is time now to construct the random variable $Z'$ that will generate the
$\sigma$-algebra $\z'$. We will use the existence of the set $S_N$ and the
section $\tau$ defined by the lemmas \ref{SN} and \ref{meas-tau}. We need also
a new mapping $\sigma : t \lra G$ defined as follows. Let $g=h^{-1}xh$ be an
element of $t$. The element $h$ is not uniquely defined by $g$, but only its
class in $G/C(x)$. Then $\sigma$ is the composition of this well-defined map
$t \lra G/C(x)$ with a measurable section $G/C(x)\lra G$. The existence of
such a measurable section can be shown in the same way as we proved that
$\tau$ exists. 

We fix now $N$ loops $L^1_1,\ldots,L^1_N$ in $PGL_{m_1}M_1$ such that the
distribution of $(H_{L_1},\ldots,H_{L_N})$ has a density with respect to the
Haar measure on $G^N$. Then $\mu_M(u,t)$-almost surely, this $N$-uple takes
its values in $S_N$. Let $\omega\in \Omega_M$ be typical with this respect.
Set $\tau(\h_{L,L^1_1,\ldots,L^1_N})=(x,x_1,\ldots,x_N)$ and choose $x_c\in G$
such that $\h_{L,L^1_1,\ldots,L^1_N,c}=[x,x_1,\ldots,x_N,x_c]$. Then $x_c$ is
defined up to conjugation by an element of $C(x_1)\cap \ldots \cap
C(x_N)=Z(G)$, so that it is in fact uniquely defined. Now, set
$x'_c=\sigma(x_c^{-1}xx_c)$. The product $x_c{x'_c}^{-1}$ belongs to $C(x)$
and we set
$$Z'(\omega) =x_c{x'_c}^{-1}.$$

Note that, in a similar fashion to $Z$, $Z'$ depends on $\omega$ only through
$\h_{L,L^1_1,\ldots,L^1_N,c}(\omega)$, via a function that we denote
$\widetilde Z': G^{N+2}/\Ad \lra C(x)$. It is also easy to check that 
$$\widetilde Z'([x,x_1,\ldots,x_N,zx_c])=z\widetilde
Z'([x,x_1,\ldots,x_N,x_c])$$
holds for all $z$ in $C(x)$.

\begin{proposition}
The random variable $Z'$ is independent of $\a_1$ and uniformly distributed on
$C(x)$ under the measure $\mu_M(u,t)$.
\end{proposition}

In particular, $Z'$ is never a constant random variable, since the centralizer
of an element contains at least a one-parameter subgroup of $G$. If $G$ is
Abelian, then $C(x)=G$, so that in this case, it is necessary to bring a
uniform $G$-valued random variable in order to close a handle on $M$. This
should be compared to the properties of the random holonomy along
homologically non-trivial loops (see \ref{law_1}).\\

\pf. Let $f$ be a continuous function on $C(x)$. The function $Z'$ can be
seen as a function of the random variable $\h_{L,L^1_1,\ldots,L^1_N,c}$ and so
will we do. Using \ref{conda1}, we get
\begin{eqnarray*}
E[f(Z)|\a_1] &=& E[f(\widetilde Z'(\h_{L,L^1_1,\ldots,L^1_N,c}))] \\
&=& \int_{G^{n+2}} f\circ \widetilde Z' \;
d\alpha^{n+2}_x(\h_{L,cLc^{-1},L,L^1_1,\ldots,L^1_N}).
\end{eqnarray*}

Let us fix $\omega$ and set
$(x,x_1,\ldots,x_N)=\tau(\h_{L,L^1_1,\ldots,L^1_N}(\omega))$. We assume that
$C(x_1)\cap \ldots \cap C(x_N)=Z(G)$. Let $x_c$ be the unique element of $G$
such that $\h_{L,L^1_1,\ldots,L^1_N,c}(\omega)=[x,x_1,\ldots,x_N,x_c]$. Then in
particular,
$\h_{L,cLc^{-1},L,L^1_1,\ldots,L^1_N}=[x,x_c^{-1}xx_c,x_1,\ldots,x_N]$, so that by
definition of $\alpha^{n+2}_x$,
\begin{eqnarray*}
E[f(Z')|\a_1](\omega) &=& \int_{C(x)} f(\widetilde
Z'(x,x_1,\ldots,x_N,zx_c)) \; dz\\ 
&=& \int_{C(x)} f(z \widetilde Z'(x,x_1,\ldots,x_N,x_c)) \; dz \\
&=& \int_{C(x)} f(z) \; dz.
\end{eqnarray*}

This conditional expectation is a constant for all $f$, hence the random
variable $Z'$ is independent of $\a_1$. \qed

It remains to prove that the completion of $\a_1\vee \sigma(Z')$ contains
$\a$.

\begin{lemma} The random variable $\h_{L,L^1_1,\ldots,L^1_N,c}$ is measurable with
  respect to the $\sigma$-algebra $\a_1\vee \sigma(Z')$.
\end{lemma}

\pf. By construction of $Z'$, the random variable $\h_{L,L^1_1,\ldots,L^1_N,c}$ is
almost surely equal to
$$[\tau(\h_{L,L^1_1,\ldots,L^1_N}),Z'\sigma(\h_{cLc^{-1}})],$$ 
which is measurable with respect to $\a_1 \vee \sigma(Z')$. \qed

The next proposition finishes the proof of the theorem \ref{tribuz2}:

\begin{theorem} \label{generate2} The completion of the $\sigma$-algebra $\a_1
  \vee \sigma(\h_{L,L^1_1,\ldots,L^1_N,c})$ with respect to the measure $\mu_M(t)$
  contains $\a$.
\end{theorem}

\pf. Let $l_1,\ldots,l_p$ be $p$ loops of $i(PGL_{m_1}M_1\cup\{c\})$.
Since this proof is very close to that of \ref{generate}, we do not write all
details. The variables $\h_{c,L^1_1,\ldots,L^1_N}$ and
$\h_{L^1_1,\ldots,L^1_N,l_1,\ldots,l_p}$ determine $\h_{c,l_1,\ldots,l_p}$ because
they determine $\h_{c,L^1_1,\ldots,L^1_N,l_1,\ldots,l_p}$. This is proved using a
map
$$(G\times S_N)/\Ad \times_{S_N}(S_N \times G^p)/\Ad \build{\lra}_{}^{\sigma'}
(G\times S_N \times G^p)/\Ad.$$
It is easily checked that this map can be
defined in such a way that $\h_{c,L,l}=\sigma'(\h_{c,L},\h_{L,l})$. The fact
that $\sigma'$ is measurable terminates the proof. \qed

\section{Conditional partition functions}

The propositions \ref{Z_1} and \ref{Z_2} show that the conditional partition
functions deserve to be studied separately. We are interested in the
conditional partition functions with respect to the boundary components of a
surface. The importance of these functions had been pointed out by Witten
\cite{Witten}. He already proved their algebraic properties using character
expansions.

\subsection{Algebraic properties of the partition functions}

Let us summarize the properties of the conditional partition function that
were proved at different points in the preceding chapters.

Let $(M,\sigma)$ be a surface, with a boundary $\partial M=N_1\cup\ldots\cup
N_p$ or without boundary. For any graph $\Gamma$ on $M$ and any
$x_1,\ldots,x_p \in G$, the number $\int_{G^\Gamma} D^\Gamma \; d\nu_{x_1}
\ldots d\nu_{x_p} \; dg'$ is well defined (see section \ref{conditionning}).
By the lemma \ref{inv_part}, it depends only on the conjugacy classes of the
$x_i$'s: it is a central function of the $x_i$'s. By the proposition
\ref{inv_part_gr}, which is true on a surface with boundary by the lemma
\ref{valid_boundary}, this number does not depend on $\Gamma$. If
$t_1=[x_1],\ldots,t_p=[x_p]$, we denote this number by $Z_M(t_1,\ldots,t_p)$,
or just $Z_M$ if $M$ is closed.

Consider now an area-preserving diffeomorphism between $(M,\sigma)$ and
another surface $(M',\sigma')$, i.e. a diffeomorphism that sends $\sigma$ to
$\sigma'$. Then an expression like $\int_{G^\Gamma} D^\Gamma \;
d\nu_{x_1}\ldots d\nu_{x_p}\; dg'$ is obviously invariant by this
diffeomorphism. Thus, the function $Z_M$ depends on $M$ only through its class
modulo area-preserving diffeomorphisms. As a consequence of Moser's theorem,
that we extended to the case of surfaces with boundary in the proof of
\ref{metric_Moser}, this class is easily parametrized by a triple $(p,g,T)\in
\N^2\times \R^*_+$, where $p$ is the number of components of $\partial M$, $g$
the genus of $M$ and $T$ the total surface of $M$. Another consequence of this
invariance and of Moser's theorem is the symmetry of $Z_M$. Indeed, given any
two components of $\partial M$, there exists a diffeomorphism of $M$ that
permutes these components, hence an area-preserving diffeomorphism. Thus, for
any $i$ and $j$ such that $1\leq i<j\leq p$,
$Z_M(t_1,\ldots,t_i,\ldots,t_j,\ldots,t_p)=Z_M(t_1,\ldots,t_j,\ldots,t_i,
\ldots,t_p)$.

Let us give an expression of the function $Z_M$ that makes clear that it
depends on $M$ only through $p,g,T$. We consider a graph with only one face on
$M$, such that the boundary of this face is $[a_1,b_1]\ldots[a_g,b_g]
c_1^{-1}N_1c_1 \ldots c_p^{-1}N_pc_p$, where $a_i,b_i$ are the edges of a
polygonal fundamental domain in the universal covering of $M$ and each $c_i$
joins $N_i$ to a point on the boundary of this fundamental domain. We find
\begin{eqnarray}
&& Z_{p,g,T}(t_1,\ldots,t_n)= \label{integral} \\ 
&& \hskip -.7cm =\int_{G^{2g+p}} p_{\sigma(M)}(y_1^{-1} x_1 y_1 \ldots
y_p^{-1} x_p y_p [a_1,b_1] \ldots [a_g,b_g]) \; da_1 db_1 \ldots da_g db_g
dy_1 \ldots dy_p, \nonumber
\end{eqnarray} 
where the $x_i$'s are arbitrary representatives of the $t_i$'s, and
\begin{equation} \label{integral_closed}
Z_{0,g,T}=\int_{G^{2g+p}} p_{\sigma(M)}([a_1,b_1] \ldots [a_g,b_g]) \; da_1
db_1 
\ldots da_g db_g
\end{equation}
when $p=0$, i.e. when $M$ is closed. From now on, we index the
function $Z$ by the triple $(p,g,T)$ instead of the surface $M$.

The expression \ref{integral} shows also that $Z_{p,g,T}$ is a smooth central
function on $G^p$ and a continuous function on $(G/\Ad)^p$. On the other hand,
the symmetry of $Z_{p,g,T}$ is less obvious in this form. Using character
expansions, it is possible to give a manifestly symmetric expression of
$Z_{p,g,T}$.  The reader which is not familiar with the characters of a
compact Lie group should read the beginning of the section \ref{carac} before
to go further.  Using the expansion of the heat kernel proved in proposition
\ref{heat_kernel_expansion}, we transform (\ref{integral}) and
(\ref{integral_closed}) into:
\begin{equation} \label{Z_carac}
Z_{p,g,T}(t_1,\ldots,t_p)= \sum_{\alpha\in\widehat G} (\dim\alpha)^{2-2g}
e^{-{{c_2(\alpha)}\over 2} T} \prod_{i=1}^p {{\chi_\alpha(t_i)}\over
  {\dim\alpha}},
\end{equation}
\begin{equation}
Z_{0,g,T}=\sum_{\alpha\in\widehat G} (\dim\alpha)^{2-2g}
e^{-{{c_2(\alpha)}\over 2} T}.
\end{equation}

Before to state these results in a theorem, recall that $G/\Ad$ is endowed
with the image measure of the Haar measure by the canonical projection $G\lra
G/\Ad$. As we said at the beginning of this section, this theorem was
essentially already proved by Witten.

\begin{theorem} \label{mainZ} For each $(p,g,T)\in \N^2 \times \R^*_+$, the
  function $Z_{g,p,T}$ is a continuous symmetric function on $(G/\Ad)^p$.
  Moreover, for any $(p',g',T')$ and any $t_1,\ldots,t_p,t'_1,\ldots,t'_{p'}
  \in
  G/\Ad$, the following relations hold: \\
  \begin{eqnarray} 
  \int_{G/\Ad} Z_{p+1,g,T}(t_1,\ldots,t_p,t) Z_{p'+1,g',T'}
  (t^{-1},t'_1,\ldots,t'_{p'}) \; dt =&& \nonumber \\
  && \hskip -5cm
  Z_{p+p',g+g',T+T'}(t_1,\ldots,t_p,t'_1,\ldots,t'_{p'}),\label{ZZ_1} 
  \end{eqnarray} 
  \begin{equation} \label{ZZ_2}
  \int_{G/\Ad} Z_{p+2,g,T}(t_1,\ldots,t_p,t,t^{-1})\; dt =
  Z_{p,g+1,T}(t_1,\ldots,t_p).
  \end{equation}
\end{theorem}

\pf. The symmetry and continuity of $Z_{p,g,T}$ were already discussed. The
relation (\ref{ZZ_1}) is a consequence of the proposition \ref{Z_1}. Indeed,
in this proposition, the number of components of the boundary of $M$ is
$p_1+p_2-2$, where $p_1$ and $p_2$ are the number of components of $\partial
M_1$ and $\partial M_2$ and its genus and total surface are the sums of those
of $M_1$ and $M_2$. This gives:
\begin{eqnarray*}
\int_{G/\Ad} Z_{p+1,g,T}(t_1,\ldots,t_p,t) Z_{p'+1,g',T'}
(t^{-1},t'_1,\ldots,t'_{p'}) \; dt =&& \\ 
&& \hskip -3cm \int_{G/\Ad} Z_M(t_1,\ldots,t_p,t,t'_1,\ldots,t'_{p'}) \; dt.
\end{eqnarray*}
In this last partition
function, the variable $t$ corresponds to an interior loop of $M$, not to a
component of the boundary. If we compute this function using an expression
like $\int_{G^\Gamma} D^\Gamma \; d\nu_x \ldots$, where $[x]=t$, we see that
the conditioning with respect to this interior loop disappears if we integrate
over $x$, so that the last integral is exactly equal to
$Z_{p+p',g+g',T+T'}(t_1,\ldots,t_p,t'_1,\ldots,t'_{p'})$.

Similarly, the relation (\ref{ZZ_2}) is a consequence of the proposition
\ref{Z_2}. In this gluing operation, the surface $M_1$ had lost two components
of its boundary and gained one handle. Thus,
$$\int_{G/\Ad} Z_{p+2,g,T}(t_1,\ldots,t_p,t,t^{-1})\; dt = \int_{G/\Ad}
Z_M(t_1,\ldots,t_p,t) \; dt.$$
Just as above, $t$ disappears when we integrate
against $dt$ because it corresponds to an interior loop and we find
$Z_{p,g+1,T}(t_1,\ldots,t_p)$. \qed

\subsection{Building bricks of the theory}

\begin{figure}[h]
\begin{center}
\input{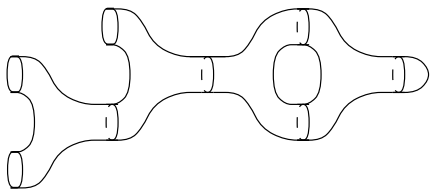}
\end{center}
\caption{An example of decomposition in three-holed spheres and disks.
\label{bric}}
\end{figure}

The two relations (\ref{ZZ_1}) and (\ref{ZZ_2}) are the analytic counterparts
of the behaviour of the Yang-Mills measure under the two basic surgery
operations. It is well known that a few elementary surfaces are enough to
build any surface by a sequence of these basic operations, namely a disk and a
three-holed sphere (see fig. \ref{bric}). It is not surprising that a
corresponding result holds for the conditional partition functions.

\begin{proposition} \label{detZ} The family of functions $Z_{p,g,T}$ is
  completely determined by the functions $Z_{1,0,T}$ and $Z_{3,0,T}$, $T>0$,
  and the relations (\ref{ZZ_1}) and (\ref{ZZ_2}).
\end{proposition}

\pf. We choose $g,p$ and construct the functions $Z_{g,p,T}$, $T>0$, starting
with the functions $Z_{1,0,T}$ and $Z_{3,0,T}$. 

Suppose first that $p+2g\geq 3$. In this case, repeated applications of
(\ref{ZZ_1}) to the function $Z_{3,0,T}$ allow to compute $Z_{p+2g,0,T}$ for
any $T>0$. Now, $g$ applications of (\ref{ZZ_2}) to $Z_{p+2g,0,T}$ give the
function $Z_{p,g,T}$.

The case $p+2g=2$ happens when $(p,g)=(0,1)$ or $(2,0)$. The first case is
that of a closed torus. Start with a three-holed sphere and glue two
components of its boundary. We get a torus with one hole. This corresponds to
(\ref{ZZ_2}) applied to $Z_{3,0,T}$ to get $Z_{1,1,T}$. Now, it remains to
glue a disk on the hole of the torus. In other words, the relation
(\ref{ZZ_1}) applied to $Z_{1,1,T}$ and $Z_{1,0,T}$ gives $Z_{0,1,T}$. The
case $(p,g)=(2,0)$ is that of a cylinder, which is obtained by gluing a disk
on a three-holed sphere. So, $Z_{2,0,T}$ is obtained by applying (\ref{ZZ_1})
to $Z_{3,0,T}$ and $Z_{1,0,T}$.

The case $p+2g=1$ happens only when $(p,g)=(1,0)$ and the corresponding
function is one of our building bricks.

Finally, $p=g=0$ is a closed sphere, which can be obtained by gluing two disks
together. So, (\ref{ZZ_1}) applied to $Z_{1,0,T}$ gives $Z_{0,0,T}$.\qed

The natural question arising from this result is to identify the elementary
functions $Z_{1,0,T}$ and $Z_{3,0,T}$.

\begin{proposition} \label{id_Z1} The function $Z_{1,0,T}$ is the projection
  on $G/\Ad$ of the heat kernel $p_T$ on $G$.
\end{proposition}

\pf. Any expression of $Z_{1,0,T}$, for example (\ref{integral}), proves this
assertion. \qed

The meaning of $Z_{3,0,T}$ is less obvious. Let us consider the functions
$Z_{p,g,T}$ as functions on $G^p$ invariant by adjunction on each variable.

\begin{lemma} \label{heateq} For any $(p,g)\in \N^2$, any $T,T'>0$, the
  following relation holds between central functions on $G$:
  $$e^{T{\Delta \over 2}}
  Z_{p,g,T'}(x_1,\ldots,x_{p-1},\cdot)(x)=Z_{p,g,T+T'}
  (x_1,\ldots,x_{p-1},x).$$ 
\end{lemma}

In other words, $Z_{p,g,T}$ is a solution of the heat equation in each of its
variables.\\

\pf. Given the fact that $\Delta \chi_\alpha=-c_2(\alpha) \chi_\alpha$ for any
irreducible representation $\alpha$, this assertion is a simple consequence of
\ref{Z_carac}. \qed

\begin{figure}[h]
\begin{center}
\input{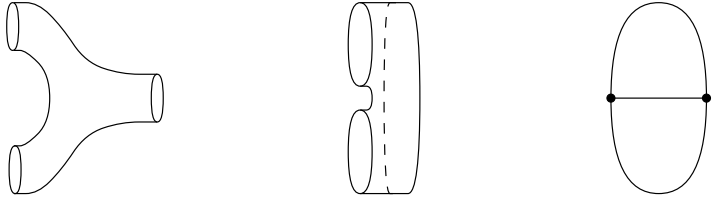}
\end{center}
\caption{A thiner and thiner three-holed sphere.
\label{pant}}
\end{figure}

This lemma shows that the algebraic meaning of $Z_{3,0,T}$, if there is one,
is contained in the formal limit $\lim_{T\to 0} Z_{3,0,T}$. Let us look at a
three-holed sphere with a very small surface (see fig. \ref{pant}).

At the $T\to 0$ limit, there remains only two adjacent circles that form a
graph. If we remind that the conditional partition function is the density of
the natural law of the holonomy along the boundary of a surface (see
\ref{natural_law}), we see in this case that $Z_{3,0,T}$ is closely related
to the multiplication in $G$.

Recall that the convolution product of two function $f,g\in L^2(G,dx)$ is
defined by $f*g(x)=\int_G f(y) g(y^{-1}x) \; dx$. It is also a
square-integrable function. Let us denote by $L^2(G)^G$ the space of central
$L^2$ functions on $G$. It is easily checked that the convolution product is a
commutative operation in $L^2(G)^G$. Indeed, let $f$ and $g$ be two central
functions.  Using the left invariance of the Haar measure and its invariance
by inversion, we get:
$$f*g(x)=\int_G f(y)g(y^{-1}x) \; dy = \int_G f(xy)g(y^{-1})\; dy =\int_G g(y)
f(y^{-1}x) \; dy = g*f(x).$$
It remains to check that $f*g$ is central.
\begin{eqnarray*}
f*g(xy)=\int_G f(z)g(z^{-1}xy)\; dz = \int_G f(z) g(y z^{-1}x) \; dz = &&\\
&& \hskip -2cm \int_G f(xzy) g(z^{-1})\; dz = g*f(yx)
\end{eqnarray*}
and the result follows by commutativity of $*$.

This product on $L^2(G)^G$ is what remains from the product on $G$ when one
considers conjugacy classes. The following result relates $Z_{3,0,T}$ to this
product.

\begin{proposition} \label{id_Z3}
  Let $f$ and $g$ be two functions of $L^2(G)^G$. Then the following equality
  holds in $L^2(G)^G$:
  $$\int_{G^2} f_1(x_1) f_2(x_2) Z_{3,0,T}(x_1,x_2,x)\; dx_1 dx_2 =
  \left[e^{-T{\Delta \over 2}} (f_1 * f_2)\right](x).$$
\end{proposition}

This gives the interpretation that we were looking for: formally, $Z_{3,0,0}$
is the distributional kernel of the operator $*:L^2(G)^G \otimes L^2(G)^G \lra
L^2(G)^G$. From this point of view, the commutativity of $*$ finds its
geometric counterpart in the fact that two holes of a three-holed sphere are
indistinguishable under area-preserving diffeomorphisms.\\

\pf. We use the fact that any central square-integrable function can be
expanded into a series of characters. Thus, it is sufficient to prove the
theorem when $f_1$ and $f_2$ are the characters of two irreducible
representations $\alpha$ and $\beta$. We use the expansion of $Z_{3,0,T}$
given by (\ref{Z_carac}).
\begin{eqnarray}
\int_{G^2} \chi_\alpha(x_1) \chi_\beta(x_2) Z_{3,0,T}(x_1,x_2,x) \; dx_1 dx_2=
&& \nonumber \\
&& \hskip -7cm
=\sum_{\gamma\in\widehat G} (\dim \gamma)^{-1} e^{-{{c_2(\gamma)}\over 2}
  T} 
\chi_\gamma(x) \int_G \chi_\alpha(x_1) \chi_\gamma(x_1) \; dx_1
  \int_G 
\chi_\beta(x_2) \chi_\gamma(x_2) \; dx_2 \nonumber\\
&& \hskip -7cm =\sum_{\gamma\in\widehat G} (\dim \gamma)^{-1}
e^{-{{c_2(\gamma)}\over 2}T}  
\chi_\gamma(x) \delta_{\alpha,\gamma} \delta_{\beta,\gamma} \nonumber\\
&& \hskip -7cm =\delta_{\alpha,\beta} e^{-{{c_2(\alpha)}\over 2} T}
{{\chi_\alpha(x)}\over {\dim\alpha}}. \label{calc}
\end{eqnarray}
On the other hand, the orthogonality relations between characters imply:
$$\chi_\alpha*\chi_\beta(x)=\int_G \chi_\alpha(y) \chi_\beta(y^{-1}x) \; dy =
\delta_{\alpha,\beta} {{\chi_\alpha(x)}\over {\dim \alpha}}.$$
Finally, the
fact that $\Delta \chi_\alpha = -c_2(\alpha) \chi_\alpha$ shows that the
expression (\ref{calc}) is exactly equal to $\exp(T{\Delta \over 2})
(\chi_\alpha * \chi_\beta)$. \qed

\subsection{Transition fonctions of the Markov field}

Consider the following very simple example. Take $M$ to be a cylinder
$S^1\times [0,1]$ endowed with the Riemannian volume of the standard product
metric, with total volume equal to 1. Pick two elements $t_0$ and $t_1$ in
$G/\Ad$ and consider the Yang-Mills measure $\mu_M(t_0^{-1},t_1)$.

For each
$s\in[0,1]$, let $l_s$ be a loop whose image is the slice $S^1\times \{s\}$ in
$M$. The family of random variables $\h_{l_s},s\in[0,1]$ is a $G/\Ad$-valued
process with index set $[0,1]$. The Markov property of the Yang-Mills field
(theorem \ref{Markov}) implies that this process is a Markov process. Let us
compute its transition functions. Choose $0<s_1<s_2<1$. Let $f_1$ and $f_2$ be
two continuous functions on $G/\Ad$. We know by the proposition
\ref{partition_law} that 
$$E_{\mu_M(t_0^{-1},t_1)}[f_1(\h_{l_{s_1}}) f_2(\h_{l_{s_2}})] =
\int_{(G/\Ad)^2} f_1(u_1) f_2(u_2)
{{Z_M(t_0^{-1},u_1,u_2,t_1)}\over{Z_{2,0,1}(t_0^{-1},t_1)}} \; du_1 du_2,$$
where the partition function in the integral is taken with respect to $l_0,
l_{s_1},l_{s_2},l_1$. Using the proposition \ref{Z_1}, we find that the
expectation $E_{\mu_M(t_0,t_1)}[f_1(\h_{l_{s_1}}) f_2(\h_{l_{s_2}})]$ is equal
to
$$ \int_{(G/\Ad)^2} f_1(u_1) f_2(u_2)
{{Z_{2,0,s_1}(t_0^{-1},u_1) Z_{2,0,s_2-s_1}(u_1^{-1},u_2)
Z_{2,0,1-s_2}(u_2^{-1},t_1)}\over{Z_{2,0,1}(t_0^{-1},t_1) }} \; du_1 du_2.$$
Thus, the conditional expectation $E[f_1(\h_{l_{s_1}})|\h_{l_{s_2}}]$ is equal
to
$$\frac{1}{Z_{2,0,1}(t_0^{-1},t_1)} \int_{G/\Ad} f_1(u)
Z_{2,0,s_1}(t_0^{-1},u) Z_{2,0,s_2-s_1}(u^{-1},\h_{l_{s_2}})
Z_{2,0,1-s_2}(\h_{l_{s_2}}^{-1},t_1) \; du.$$

\begin{figure}[h]
\begin{center}
\input{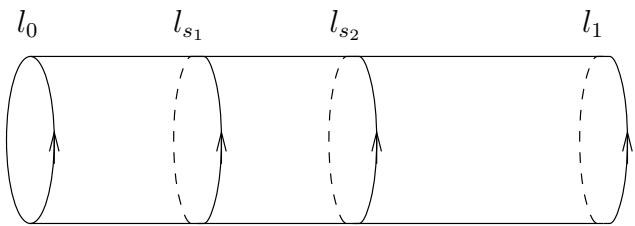}
\end{center}
\caption{Transition functions on a cylinder. \label{trans}}
\end{figure}

The transition functions of the process are exactly the functions $Z_{2,0,T}$.
This suggests that the functions $Z_{p,g,T}$ determine to some extent the law
of the random holonomy.  This was essentially the content of the proposition
\ref{partition_law}. More precisely, this proposition shows that it is
possible to write down the law of the holonomy along a family of disjoint
simple loops using only the partition functions. Using the continuity of the
random holonomy, we can extend this statement a little bit. Indeed, let
$l_1,\ldots,l_n$ be simple loops on $M$ that can be approximated by families
of disjoint simple loops in such a way that none of the components of $M$
delimited by these families has a surface tending to 0. Then the density of
the law of the variable $(\h_{l_1},\ldots,\h_{l_n})$ is the limit of the
densities of the holonomies along the approximating families. Since
$Z_{p,g,T}$ depends continuously on its parameters and also on $T$ provided
$T>0$, this limit density is also a combination of the functions $Z_{p,g,T}$.

Nevertheless, one cannot hope to express the law on $G^2/\Ad$ of a variable
like $\h_{l_1,l_2}$ using partition functions when $l_1$ and $l_2$ are based
at the same point. Indeed, partition functions are functions on the space
$(G/\Ad)^p$ which is much smaller than the space $G^p/\Ad$.


\pagevide

\end{document}